%% file: main.tex
\PassOptionsToPackage{x11names, table}{xcolor}
\documentclass[12pt]{article}
\usepackage{hyperref}
\usepackage{amsmath}
\usepackage{amssymb}
\usepackage{mathtools}
\usepackage{adjustbox}
\usepackage{amsthm}
\usepackage[utf8]{inputenc}
\usepackage{fullpage}
\usepackage{caption}
\usepackage{subcaption}
\usepackage{graphicx}
\usepackage{wrapfig}
\usepackage{natbib}
\setcitestyle{authoryear} 

\usepackage{abstract} 

\usepackage{titlesec}
\def\sectionfont{\sffamily\Large\bfseries\boldmath}
\def\subsectionfont{\sffamily\large\bfseries\boldmath}
\def\paragraphfont{\sffamily\normalsize\bfseries\boldmath}
\titleformat*{\section}{\sectionfont}
\titleformat*{\subsection}{\subsectionfont}
\titleformat*{\subsubsection}{\paragraphfont}
\titleformat*{\paragraph}{\paragraphfont}
\titleformat*{\subparagraph}{\paragraphfont}
\usepackage[labelfont={bf,sf}]{caption}  

\makeatletter
\renewcommand{\@maketitle}{
    \begin{center}
        {\sffamily\LARGE\bfseries\boldmath \@title \par}
        \vskip 1.5em
        {\large \@author \par}
        \vskip 1em
        {\large \@date \par}
    \end{center}
    \vskip -2em
}
\makeatother

\usepackage{tikz}
\usetikzlibrary{backgrounds}
\usepackage{xcolor}

\usepackage{cancel}
\usepackage{amsthm}
\usepackage{mathtools}
\newtheorem{assumption}{Assumption}
\newtheorem{theorem}{Theorem}
\newtheorem{lemma}{Lemma}
\newtheorem{definition}{Definition}

\newcommand{\purple}[1]{\textcolor{purple}{#1}}

\newcommand{\calC}{{\mathcal C}}
\newcommand{\calD}{{\mathcal D}}
\newcommand{\calE}{{\mathcal E}}
\newcommand{\calG}{{\mathcal G}}

\newcommand{\calN}{{\mathcal N}}
\newcommand{\calP}{{\mathcal P}}
\newcommand{\calS}{{\mathcal S}}
\newcommand{\calU}{{\mathcal U}}
\newcommand{\calV}{{\mathcal V}}

\newcommand{\calZ}{{\mathcal Z}}

\newcommand{\bmu}{{\bm \mu}}
\newcommand{\vnu}{{\bm \nu}}
\newcommand{\btheta}{{\bm \theta}}
\newcommand{\blambda}{{\bm \lambda}}
\newcommand{\bxi}{{\bm \xi}}
\newcommand{\bSigma}{{\bm \Sigma}}
\definecolor{myYellow}{RGB}{255,255,204}

\usepackage{booktabs}

\input{commands.tex}

\title{Deep Distributed Optimization for \\ Large-Scale Quadratic Programming}

\author{
    Augustinos D. Saravanos, 
    Hunter Kuperman, 
    Alex Oshin, \\
    Arshiya Taj Abdul,  
    Vincent Pacelli and 
    Evangelos A. Theodorou    
    \\[1em]
Autonomous Control and Decision Systems Laboratory
     \\
    Georgia Institute of Technology
    }
\date{December 2024}

\begin{document}
\maketitle
\begin{tikzpicture}[remember picture, overlay]
    \node[anchor=north west, yshift=-8mm, xshift=20mm] at (current page.north west)
    {\fontsize{8}{8}\selectfont Accepted at International Conference on Learning Representations (ICLR) 2025. Preprint version.};
\end{tikzpicture}

\begin{abstract}
Quadratic programming (QP) forms a crucial foundation in optimization, encompassing a broad spectrum of domains and serving as the basis for more advanced algorithms. Consequently, as the scale and complexity of modern applications continue to grow, the development of efficient and reliable QP algorithms is becoming increasingly vital. In this context, this paper introduces a novel deep learning-aided distributed optimization architecture designed for tackling large-scale QP problems.
First, we combine the state-of-the-art Operator Splitting QP (OSQP) method with a consensus approach to derive \textbf{\textsf{DistributedQP}}, a new method tailored for network-structured problems, with convergence guarantees to optimality. Subsequently, we unfold this optimizer into a deep learning framework, leading to \textbf{\textsf{DeepDistributedQP}}, which leverages learned policies to accelerate reaching to desired accuracy within a restricted amount of iterations. Our approach is also theoretically grounded through Probably Approximately Correct (PAC)-Bayes theory, providing generalization bounds on the expected optimality gap for unseen problems. The proposed framework, as well as its centralized version \textbf{\textsf{DeepQP}}, significantly outperform their standard optimization counterparts on a variety of tasks such as randomly generated problems, optimal control, linear regression, transportation networks and others. Notably, DeepDistributedQP demonstrates strong generalization by training on small problems and scaling to solve much larger ones (up to 50K variables and 150K constraints) using the same policy. Moreover, it achieves orders-of-magnitude improvements in wall-clock time compared to OSQP. The certifiable performance guarantees of our approach are also demonstrated, ensuring higher-quality solutions over traditional optimizers.
\end{abstract}

\input{sections/sec1_intro}

\input{sections/sec2_related_work}

\input{sections/sec3_distr_qp}

\input{sections/sec4_deep_distr_qp}

\input{sections/sec5_gen_bounds}
\input{sections/sec6_experiments}

\input{sections/sec7_conclusion}

\section*{Acknowledgements}
This work is supported by the National Aeronautics and Space Administration under ULI Grant 80NSSC22M0070 and the ARO Award \#W911NF2010151.
Augustinos Saravanos acknowledges financial support by the A. Onassis Foundation Scholarship. The authors also thank Alec Farid for helpful discussions on PAC-Bayes Theory.

\bibliography{references}

\appendix
\newpage
\input{sections/appendix}

\end{document}

%% file: commands.tex
\usepackage{amsmath,amsfonts,bm}

\def\eqref#1{equation~\ref{#1}}

\def\1{\bm{1}}

\def\vzero{{\bm{0}}}

\def\vb{{\bm{b}}}
\def\vc{{\bm{c}}}
\def\vd{{\bm{d}}}
\def\ve{{\bm{e}}}
\def\vf{{\bm{f}}}
\def\vg{{\bm{g}}}

\def\vp{{\bm{p}}}
\def\vq{{\bm{q}}}

\def\vs{{\bm{s}}}
\def\vt{{\bm{t}}}
\def\vu{{\bm{u}}}
\def\vv{{\bm{v}}}
\def\vw{{\bm{w}}}
\def\vx{{\bm{x}}}
\def\vy{{\bm{y}}}
\def\vz{{\bm{z}}}

\def\mA{{\bm{A}}}
\def\mB{{\bm{B}}}
\def\mC{{\bm{C}}}
\def\mD{{\bm{D}}}

\def\mF{{\bm{F}}}
\def\mG{{\bm{G}}}

\def\mI{{\bm{I}}}

\def\mL{{\bm{L}}}

\def\mP{{\bm{P}}}
\def\mQ{{\bm{Q}}}
\def\mR{{\bm{R}}}

\DeclareMathAlphabet{\mathsfit}{\encodingdefault}{\sfdefault}{m}{sl}
\SetMathAlphabet{\mathsfit}{bold}{\encodingdefault}{\sfdefault}{bx}{n}

\def\sD{{\mathbb{D}}}

\def\sR{{\mathbb{R}}}
\def\sS{{\mathbb{S}}}

\newcommand{\E}{\mathbb{E}}

\DeclareMathOperator*{\argmin}{arg\,min}


%% file: sections/sec1_intro.tex
\section{Introduction}

Quadratic programming (QP) serves as a fundamental cornerstone in optimization with a wide variety of applications in machine learning \citep{cortes1995support, tibshirani1996regression}, control and robotics \citep{garcia1989model, rawlings2017model}, signal processing \citep{mattingley2010real}, finance \citep{cornuejols2018optimization},  and transportation networks \citep{mota2014distributed} among other fields. Beyond its standalone applications, QP also acts as the core component of many advanced non-convex optimization algorithms such as sequential quadratic programming \citep{nocedal1999numerical}, trust-region methods \citep{conn2000trust}, augmented Lagrangian approaches \citep{houska2016augmented}, mixed-integer optimization \citep{belotti2013mixed}, etc. For these reasons, the pursuit of more efficient QP algorithms remains an ever-evolving area of research from active set \citep{wolfe1959simplex} and interior point methods \citep{nesterov1994interior} during the previous century to first-order methods such as the state-of-the-art Operator Splitting QP (OSQP) algorithm \citep{stellato2020osqp}.

\begin{wrapfigure}[15]{r}{0.45\textwidth}
\vskip -0.15in
\begin{tikzpicture}
\node[anchor=south west,inner sep=0] at (0,0){\includegraphics[width=0.45\textwidth, trim={0.5cm 0.0cm 0cm 0cm},clip]{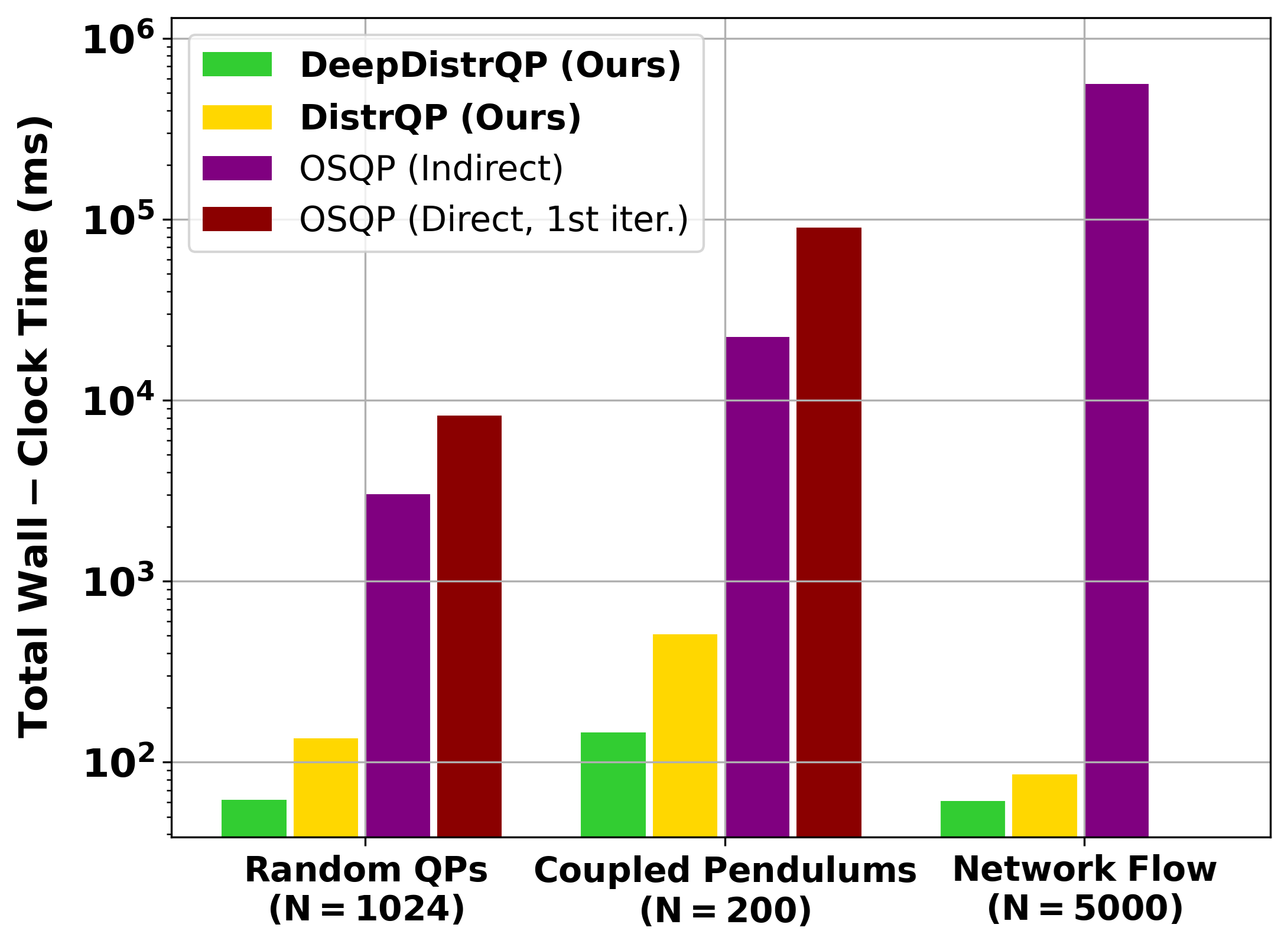}
};
\node[align=center, scale = 0.45, color = Green4] (c) at (1.4, 0.98) {\textbf{\textsf{62ms}}};
\node[align=center, scale = 0.45, color = DarkGoldenrod2] (c) at (1.78, 1.33) {\textbf{\textsf{129ms}}};
\node[align=center, scale = 0.45, color = DarkOrchid2] (c) at (2.26, 2.75) {\textbf{\textsf{3s}}};
\node[align=center, scale = 0.45, color = Firebrick2] (c) at (2.7, 3.22) {\textbf{\textsf{8.2s}}};
\node[align=center, scale = 0.45, color = Green4] (c) at (3.42, 1.37) {\textbf{\textsf{146ms}}};
\node[align=center, scale = 0.45, color = DarkGoldenrod2] (c) at (3.83, 1.95) {\textbf{\textsf{511ms}}};
\node[align=center, scale = 0.45, color = DarkOrchid2] (c) at (4.33, 3.66) {\textbf{\textsf{22.4s}}};
\node[align=center, scale = 0.45, color = Firebrick2] (c) at (4.77, 4.33) {\textbf{\textsf{90.1s}}};
\node[align=center, scale = 0.45, color = Green4] (c) at (5.57, 0.96) {\textbf{\textsf{61ms}}};
\node[align=center, scale = 0.45, color = DarkGoldenrod2] (c) at (6.0, 1.12) {\textbf{\textsf{86ms}}};
\node[align=center, scale = 0.45, color = DarkOrchid2] (c) at (6.43, 5.13) {\textbf{\textsf{9m 58s}}};
\node[align=center, scale = 0.45, color = Firebrick2] (c) at (6.87, 0.77) {\textbf{\textsf{N/A}}};
\end{tikzpicture}
\vspace{-0.8cm}
\captionof{figure}{\textbf{Wall-clock time comparison.} DeepDistributedQP, DistributedQP (ours) and OSQP on large-scale QPs.}
\label{fig: intro figure}
\end{wrapfigure}
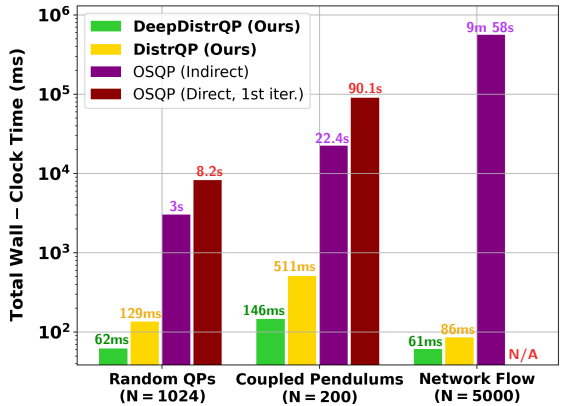
As the scale of modern decision-making applications rapidly increases, there is an emerging interest in developing effective optimization architectures for addressing high-dimensional problems. Given the fundamental role of QP in optimization, there is a clear demand for algorithms capable of solving large-scale QPs with thousands, and potentially much more, variables and constraints. Such problems arise in diverse applications including sparse linear regression \citep{mateos2010distributed} and support vector machines \citep{navia2006distributed} with decentralized data, multi-agent control \citep{van2017distributed}, resource allocation \citep{huang2014distribution}, network flow \citep{mota2014distributed}, power grids \citep{lin2012distributed} and image processing \citep{soheili2020dqpfs}. Traditional centralized optimization algorithms are inadequate for solving such problems at scale (see for example Fig. \ref{fig: intro figure}), prompting the development of distributed methods that leverage the underlying network/decentralized structure to parallelize computations. In this context, the Alternating Direction Method of Multipliers (ADMM) has gained widespread popularity as an effective approach for deriving distributed algorithms \citep{boyd2011distributed, mota2013d}. Nevertheless, as scale increases, such algorithms continue to face significant challenges such as their need for \textit{meticulous tuning}, lack of \textit{generalization guarantees} and restrictions on the \textit{allowed number of iterations} imposed by computational or communication limitations.



Learning-to-optimize has recently emerged as a methodology for enhancing existing optimizers or developing entirely new ones through training on sample problems \citep{chen2022learning, shlezinger2022model, amos2023tutorial}. A notable approach within this paradigm is \textit{deep unfolding}, which under the realistic assumption of computational budget restrictions, 
unrolls a fixed number of iterations as layers of a deep learning network and learns the optimal parameters for improving performance \citep{monga2021algorithm, shlezinger2022model}. Our key insight is that deep unfolding is particularly well-suited for overcoming the limitations of \textit{distributed constrained optimization}, as it can eliminate the need for extensive tuning, manage iteration restrictions and enhance generalization. However, to our best knowledge, its combination with distributed ADMM has only recently been explored in \cite{noah2024distributed}. While this framework shows promising initial results, it relies on a relatively simple setup that studies unconstrained problems, assumes local updates consisting of gradient steps, focuses solely on parameter tuning, and is not accompanied by any performance guarantees.

\paragraph{Contributions.} This paper introduces a novel deep learning-aided distributed optimization architecture for solving large-scale constrained QP problems. Our proposed approach relies on unfolding a newly introduced distributed QP algorithm as a supervised learning framework for a prescribed number of iterations. To our best knowledge, this is the first work to propose a deep unfolded architecture for distributed constrained optimization using ADMM, despite its widespread popularity. Our framework demonstrates remarkable performance and scalability, being trained exclusively on low-dimensional problems and then effectively applied to much higher-dimensional ones. Furthermore, its performance is theoretically supported by establishing guarantees based on generalization bounds from statistical learning theory. We believe that this work lays the foundation for developing learned distributed optimizers capable of handling large-scale constrained optimization problems without requiring training at such scales. Our specific contributions can be summarized as follows:
%
\begin{itemize}
\item First, we introduce \textbf{\textsf{DistributedQP}}, a new distributed quadratic optimization method that combines the well-established OSQP solver with a consensus approach, to achieve parallelizable computations. We further prove that DistributedQP is guaranteed to converge to optimality, even under local iteration-varying algorithm parameters.
\item Then, we propose \textbf{\textsf{DeepDistributedQP}}, a deep learning-aided distributed architecture that unrolls the iterations of DistributedQP in a supervised manner, learning feedback policies for the algorithm parameters. As a byproduct, we also present \textbf{\textsf{DeepQP}}, its centralized counterpart which corresponds to unfolding the standard OSQP solver.

\item To certify the performance of the learned solver, we establish generalization guarantees on the optimality gap of the final solution of DeepDistributedQP for unseen problems using Probably Approximately Correct (PAC)-Bayes theory.

\item Finally, we present an extensive experimental evaluation that validates the following:
\begin{itemize}
    \item For centralized QPs, DeepQP consistently outperforms OSQP requiring 1.5-3 times fewer iterations for achieving the desired accuracy.
    \item DeepDistributedQP successfully scales for \textit{high-dimensional} problems (up to 50K variables and 150K constraints), despite being \textit{trained exclusively} on much \textit{lower-dimensional} ones. Furthermore, both DeepDistributedQP and DistributedQP outperform OSQP in wall-clock time by orders of magnitude as dimension increases, which indicates their advantage against conventional centralized solvers.
    \item The proposed PAC bounds offer valuable guarantees on the quality of solutions produced by DeepDistributedQP for unseen problems from the same class.
\end{itemize}

\end{itemize}

%% file: sections/sec2_related_work.tex
\section{Related Work}

This section provides an overview of the existing related literature from the angles of distributed optimization and learning-to-optimize approaches.

\paragraph{Distributed optimization with ADMM.} Distributed ADMM algorithms have emerged as a scalable approach for addressing large-scale optimization problems \citep{boyd2011distributed, mota2013d}. Despite their significant applicability to machine learning \citep{mateos2010distributed}, robotics \citep{shorinwa2024distributed} and many other fields, their successful performance has been shown to be highly sensitive to the proper tuning of its underlying parameters \citep{xu2017adaptive, saravanos2023distributed}. Moreover, tuning parameters for large-scale problems is often tedious and time-consuming, making it desirable to develop effective \textit{learned} optimizers that can be trained on smaller problems instead. Furthermore, even if an distributed optimizer performs well for a specific problem instance, its generalization to new problems remains challenging to verify. These challenges constitute our main motivation for studying learning-aided distributed ADMM architectures. We also note that an ADMM-based distributed QP solver resembling a simpler version of DistributedQP was presented in \cite{Pereira-RSS-22}, but it focused on multi-robot control and lacked theoretical analysis.

\paragraph{Learning-to-optimize for distributed optimization.} The concept of integrating learning-to-optimize approaches into distributed optimization is particularly compelling, as algorithms of the latter class typically rely on a significant amount of designing and tuning by experts. Nevertheless, the area of distributed learning-to-optimize methods remains largely unexplored. For instance, although ADMM has achieved widespread success in distributed constrained optimization, its unfolded extension as a deep learning network has only been recently explored by \cite{noah2024distributed}. This framework demonstrates promising results, but it is limited to an unconstrained problem formulation, assumes gradient-based local updates, focuses solely on parameter tuning and lacks formal performance guarantees. \cite{biagioni2020learning} presented an ADMM framework which utilizes recurrent neural networks for predicting the converged values of the variables demonstrating substantial improvements in convergence speed. In \cite{zeng2022reinforcement}, a reinforcement learning (RL) approach for learning the optimal parameters of distributed ADMM was proposed, showing promising speed improvements, but requiring a substantial amount of training effort.

Beyond distributed ADMM, \cite{wang2021decentralized} proposed unrolling two decentralized first-order optimization algorithms (ProxDGD and PG-Extra) as graph neural networks (GNNs) for addressing the decentralized statistical inference problem. Similarly, \cite{hadou2023stochastic} presented an distributed gradient descent algorithm unrolled as a GNN focusing on the federated learning problem setup. From a different point of view, \cite{he2024mathematics} recently introduced a distributed gradient-based learning-to-optimize framework for unconstrained optimization which partially imposes structure on the learnable updates instead of unrolling predefined iterations. A deep RL approach for adapting the local updates of the approximate method of multipliers was recently proposed in \cite{zhu2023deep}. Finally, \cite{kishida2020deep} and \cite{ogawa2021deep} presented distributed learned optimization methods for tackling the average consensus problem.

\paragraph{Learning-to-optimize for (centralized) QP.} Recent works have focused on accelerating QP through learning; however these efforts have solely concentrated on a centralized setup. In particular, \cite{ichnowski2021accelerating} introduced an RL-based algorithm for accelerating OSQP demonstrating promising reductions in iterations, yet training this algorithm incurs significant computational costs. From a different perspective, \cite{sambharya2023end, sambharya2024learning} focused on learning-to-initialize fixed-point methods including OSQP, while maintaining constant parameters in the unrolled algorithm iterations. Concurrently with the development of the present work, \cite{sambharya2024learningb} presented a methodology for selecting the optimal algorithm parameters for various first-order optimization methods. Considering OSQP as the unrolled method coincides with the open-loop version of the proposed DeepQP framework without any notion of feedback policies.

\paragraph{Generalization guarantees for learning-to-optimize.} The works in \cite{sucker2023pac} and \cite{sucker2024learning} presented generalization bounds for learned optimizers, considering the update function as a gradient step or a multi-layer perceptron, respectively. \cite{sambharya2024data} recently also explored incorporating PAC-Bayes bounds in learning-to-optimize methods without assuming a specific underlying algorithm structure. However, our approach differs fundamentally, as their method employs a binary error function, whereas in this work we directly establish bounds based on the optimality gap of the final solution.

%% file: sections/sec3_distr_qp.tex
\section{Distributed Quadratic Programming}
\label{sec: distr qp}

\subsection{Problem Formulation}

\begin{wrapfigure}[12]{r}{0.4\textwidth}
    \vskip -0.35in
\begin{tikzpicture}
\node[anchor=south west,inner sep=0] at (0,0){\includegraphics[width=0.4\textwidth, trim={0.5cm 0.4cm 0.4cm 0.13cm},clip]{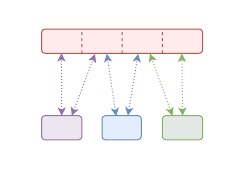}
};
\node[align=center, scale = 0.9] (c) at (3.25, 4.38) {\textbf{\textsf{Global variable components}}};
\node[align=center, scale = 0.9] (c) at (3.2, -0.15) {\textbf{\textsf{Local variables}}};
\node[align=center, scale = 1.0] (c) at (1.1, 3.6) {$\vw_1$};
\node[align=center, scale = 1.0] (c) at (2.5, 3.6) {$\vw_2$};
\node[align=center, scale = 1.0] (c) at (3.91, 3.6) {$\vw_3$};
\node[align=center, scale = 1.0] (c) at (5.32, 3.6) {$\vw_4$};
\node[align=center, scale = 0.8] (c) at (1.02, 2.2) {$[\vx_1]_1 = \vw_1$ \\ $[\vx_1]_2 = \vw_2$};
\node[align=center, scale = 0.8] (c) at (3.15, 2.2) {$[\vx_2]_1 = \vw_2$ \\ $[\vx_2]_2 = \vw_3$};
\node[align=center, scale = 0.8] (c) at (5.3, 2.2) {$[\vx_3]_1 = \vw_3$ \\ $[\vx_3]_2 = \vw_4$};
\node[align=center, scale = 1.0] (c) at (1.17, 0.58) {$\vx_1$};
\node[align=center, scale = 1.0] (c) at (3.2, 0.58) {$\vx_2$};
\node[align=center, scale = 1.0] (c) at (5.31, 0.58) {$\vx_3$};
\end{tikzpicture}
\vspace{-0.73cm}
\captionof{figure}{Example of consensus mapping $\calG$ in problem (\ref{eq: distr qp problem}).}
\label{fig: consensus figure}
\end{wrapfigure}
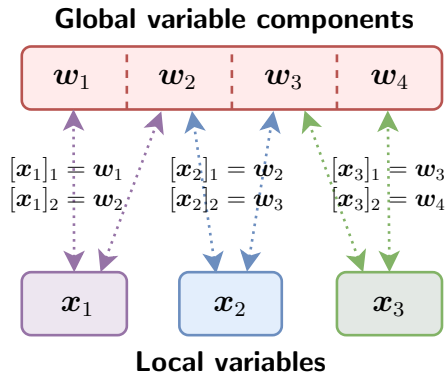
A convex QP problem is expressed in a general \textit{centralized} form as
\begin{equation}
\min ~ \frac{1}{2} \vx^\top \mQ \vx + \vq^\top \vx
\quad \mathrm{s.t.} 
\quad \mA \vx \leq \vb,
\label{eq: general qp problem}
\end{equation}
where $\vx \in \sR^n$ is the decision vector and 
$\zeta = \{ 
\mQ \in \sS_{++}^n, 
\vq \in \sR^n,
\mA \in \sR^{m \times n}, 
\vb \in \sR^m 
\}$ are the problem data.
%
\footnote{Note that equality constraints can also be captured as pairs of inequalities.} 
%
As the scale of such problems increases to higher dimensions, there is often an underlying networked/decentralized structure that could be leveraged for achieving distributed computations. This work aims to address problems characterized by such structures. Let $\vw \in \sR^n$ be the main global variable and $\vx_i \in \sR^{n_i}$ be local variables $i \in \mathcal{V} = \{1, \dots, N\}$. Then, assume a mapping $(i,j) \mapsto \mathcal{G}(i,j)$ from all index pairs $(i,j)$ of local variable components $[\vx_i]_j$ to indices $l = \mathcal{G}(i,j)$ of global components $\vw_l$\footnote{This formulation is adopted from the standard consensus ADMM framework \citep{boyd2011distributed}, wherein local variables are typically associated with their respective computational nodes.} - for an example see Fig. \ref{fig: consensus figure}. We consider QP problems of the following \textit{distributed consensus} form:
%
\begin{equation}
\min ~ \sum_{i \in \calV} \frac{1}{2} \vx_i^\top \mQ_i \vx_i + \vq_i^\top \vx_i
\quad \mathrm{s.t.} 
\quad \mA_i \vx_i \leq \vb_i,
\quad \vx_i = \tilde{\vw}_i, 
\quad i \in \calV,
\label{eq: distr qp problem}
\end{equation}
where the problem data are now given by $\zeta = \{ \zeta_i \}_{i = 1}^N$ with $\zeta_i = ( 
\mQ_i \in \sS_{++}^{n_i}, 
\vq_i \in \sR^{n_i},
\mA_i \in \sR^{m_i \times n_i}, 
\vb_i \in \sR^{m_i} 
).$
%
The vector $\vx = [\{ \vx_i \}_{i \in \calV} ]$ is the concatenation of all local variables, while $\tilde{\vw}_i \in \sR^{n_i}$, defined as $\tilde{\vw}_i = [ \{\vw_l \}_{l \in \calG(q,j) : q = i } ]$, is the selection of global variable components that correspond to the components of $\vx_i$. This form captures a wide variety of large-scale QPs found in machine learning \citep{mateos2010distributed, navia2006distributed}, optimal control \citep{van2017distributed}, transportation networks, \citep{mota2014distributed}, power grids \citep{lin2012distributed}, resource allocation \citep{huang2014distribution} and many other fields.

\subsection{DistributedQP: The Underlying Optimization Algorithm}
\label{sec: distr qp - the algo}

This section introduces a new distributed algorithm named \textbf{\textsf{DistributedQP}} for solving problems of the form (\ref{eq: distr qp problem}). The proposed method can be viewed as a combination of consensus ADMM \citep{boyd2011distributed} and OSQP using local iteration-varying penalty parameters.

Let us introduce the auxiliary variables $\vz_i, \vs_i \in \sR^{m_i}$, such that problem (\ref{eq: distr qp problem}) can be reformulated as
\begin{equation}
\min ~ \sum_{i \in \calV} \frac{1}{2} \vx_i^\top \mQ_i \vx_i + \vq_i^\top \vx_i
\quad \mathrm{s.t.} 
\quad \mA_i \vx_i = \vz_i,
~~ \vs_i \leq \vb_i, 
~~ \vz_i = \vs_i, 
~~ \vx_i = \tilde{\vw}_i, 
~~ i \in \calV.
\label{eq: distr qp problem - v2}
\end{equation}
The proposed DistributedQP algorithm is summarized below, where $k$ denotes iterations:
\begin{enumerate}
\item \textbf{\textsf{Local updates for $\vx_i, \vz_i$.}} For each node $i \in \calV$, solve in parallel:
\begin{equation} 
\begin{bmatrix} 
\mQ_i + \mu_i^k \mI & \mA_i^\top \\
\mA_i & - 1 / \rho_i^k \mI 
\end{bmatrix}
\begin{bmatrix}
\vx_i^{k+1} \\
\vnu_i^{k+1}
\end{bmatrix}
= 
\begin{bmatrix}
- \vq_i + \mu_i^k \tilde{\vw}_i - \vy_i
\\
\vz_i - 1 / \rho_i^k \blambda_i
\end{bmatrix},
\label{eq: DistrQP x update}
\end{equation}
and then update in parallel:
\begin{equation}
\vz_i^{k+1} =  \vs_i^k + 1 / \rho_i^k (\vnu_i^{k+1} - \blambda_i^k).
\label{eq: DistrQP z update}
\end{equation}
\item \textbf{\textsf{Local updates for $\vs_i$ and global update for $\vw$.}} For each $i \in \calV$, update in parallel:
\begin{equation}
\vs_i^{k+1} = \Pi_{\vs_i \leq \vb_i}
\left( \alpha^k \vz_i^{k+1} + (1 - \alpha^k) \vs_i^k + \blambda_i^k / \rho_i^k \right).
\label{eq: DistrQP s update}
\end{equation}
In addition, each global variable component $\vw_l$ is updated through:
\begin{equation}
\vw_l^{k+1} = \alpha^k \frac{\sum_{\mathcal{G}(i,j) = l} \mu_i^k [\vx_i]_j}{\sum_{\mathcal{G}(i,j) = l} \mu_i^k}
+ (1 - \alpha^k) \vw_l^k.
\label{eq: DistrQP global update}
\end{equation}
\item \textbf{\textsf{Local updates for dual variables $\blambda_i, \vy_i$.}} For each $i \in \calV$, update in parallel:
\begin{align}
\blambda_i^{k+1} & = \blambda_i^k + \rho_i^k 
( \alpha^k \vz_i^{k+1} 
+ (1 - \alpha^k) \vs_i^k
- \vs_i^{k+1}),
\label{eq: DistrQP lambda update}
\\
\vy_i^{k+1} & = \vy_i^k + \mu_i^k 
( \alpha^k \vx_i^{k+1} 
+ (1 - \alpha^k) \tilde{\vw}_i^k
- \tilde{\vw}_i^{k+1}).
\label{eq: DistrQP y update}
\end{align}
\end{enumerate}
The Lagrange multipliers $\vnu_i, \blambda_i$ and $\vy_i$ correspond to the equality constraints $\mA_i \vx_i = \vz_i$, $\vz_i = \vs_i$ and $\vx_i = \tilde{\vw}_i$, respectively. The penalty parameters $\rho_i, \mu_i > 0$ correspond to $\vz_i = \vs_i$ and $\vx_i = \tilde{\vw}_i$, while $\alpha^k \in [1,2)$ are over-relaxation parameters.
A complete derivation is provided in Appendix \ref{sec: appendix - distr_qp derivation}.

\subsection{Convergence Guarantees}

Prior to unrolling DistributedQP into a deep learning framework, it is particularly important to establish that the underlying optimization algorithm is well-behaved even for iteration-varying over-relaxation and local penalty parameters, i.e., it is expected to asymptotically converge to the optimal solution of problem. This property is especially important in deep unfolding where parameters are expected to be distinct between different iterations. 

In the simpler case of $\alpha^k = 1$, $\rho_i^k = \rho$, $\mu_i^k = \mu$, the standard convergence guarantees of two-block ADMM would apply directly \citep{deng2016global}; for a detailed discussion, see Appendix \ref{sec: appendix - distr_qp standard guarantees}. Nevertheless, the introduction of local iteration-varying penalty parameters $\rho_i^k, \mu_i^k$, as well as the over-relaxation with varying parameters $\alpha^k$ makes proving the convergence of this algorithm non-trivial. 

In the following, we prove that under mild assumptions on the asymptotic behavior of the penalty parameters, the DistributedQP algorithm is guaranteed to converge to optimality. We consider the following assumption on the penalty parameters.

\begin{assumption}
As $k \rightarrow \infty$, the parameters $\rho^k_i = \rho^{k-1}_i$, $\mu^k_i = \mu^{k-1}_i$, for all $i \in \calV$.
\label{assumption: eta}
\end{assumption}
The following theorem states the convergence guarantees of DistributedQP to optimality.
\begin{theorem}[Convergence guarantees for DistributedQP] If Assumption \ref{assumption: eta} holds and $\alpha^k \in [1,2)$, then the iterates $\vw^k$ generated by the DistributedQP algorithm converge to the optimal solution $\vw^*$ of problem (\ref{eq: distr qp problem}), as $k \rightarrow \infty$.
\label{theorem: DistrQP convergence}
\end{theorem}
The proof of Theorem \ref{theorem: DistrQP convergence} and all intermediate lemmas are provided in Appendix  \ref{sec: appendix - distr_qp convergence proof}.

%% file: sections/sec4_deep_distr_qp.tex
\section{The DeepDistributedQP Architecture}
\label{sec: distr deep qp}
The proposed DeepDistributedQP architecture emerges from unfolding the iterations of the DistributedQP optimizer into a deep learning framework. Section \ref{sec: distr deep qp - main arch} illustrates the main architecture, key aspects of our methodology, as well as the centralized version DeepQP. Section \ref{sec: distr deep qp - backprop} leverages implicit differentiation during backpropagation to facilitate the training of our framework.

\begin{figure}[t]
\centering
\begin{tikzpicture}
\node[anchor=south west,inner sep=0] at (0,0){\includegraphics[width=1.0\textwidth, trim={0cm 0cm 0cm 0cm},clip]{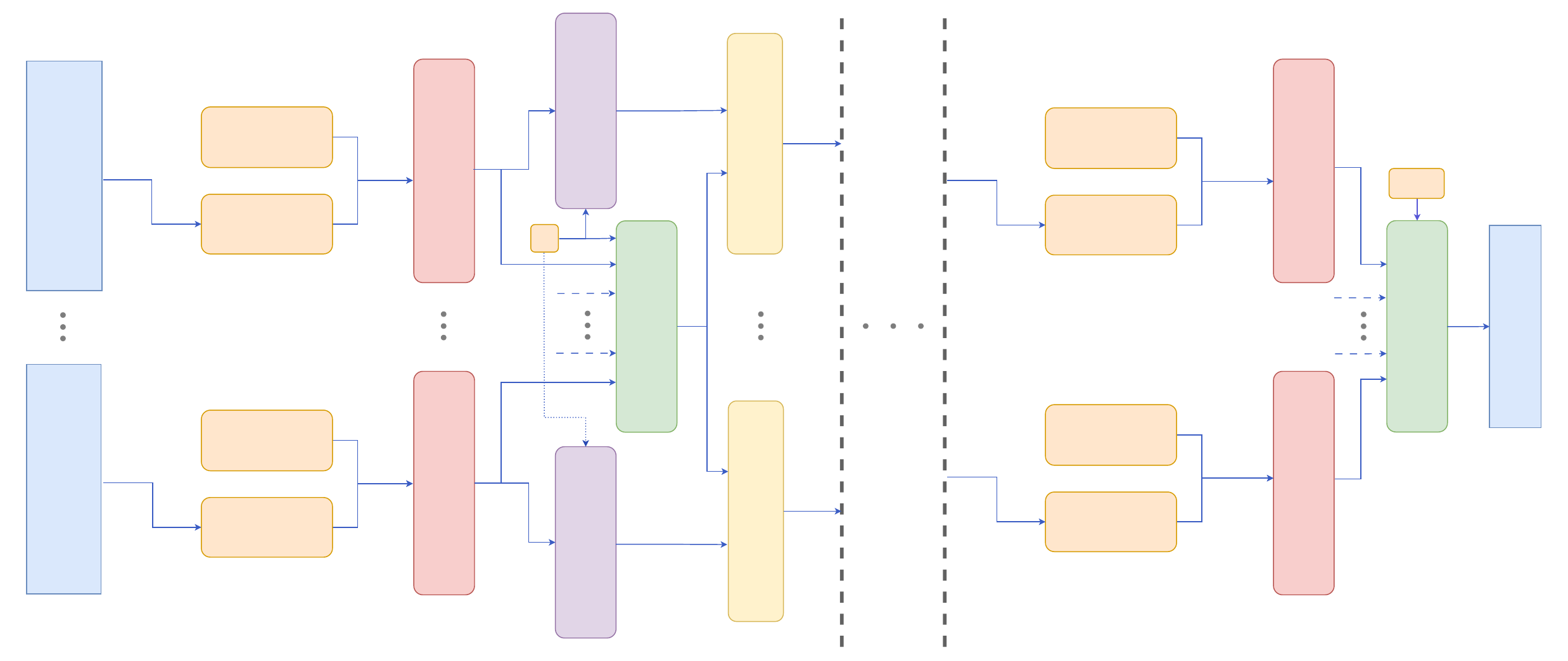}
};
\node[align=center, scale = 0.8] (c) at (5.0, -0.1) {\textbf{\textsf{Layer}} $k = 1$};
\node[align=center, scale = 0.8] (c) at (12.0, -0.1) {\textbf{\textsf{Layer}} $k = K$};
\node[align=center, fill = myYellow, scale = 0.75] (c) at (2.2, 6.7) {\textbf{\textsf{Node}} $i = 1$};
\node[align=center, fill = myYellow, scale = 0.75] (c) at (2.2, 0.35) {\textbf{\textsf{Node}} $i = N$};
\node[align=center, rotate=90, text = black, scale = 0.55] (c) at (0.68, 5.12) {\textbf{\textsf{Local Problem Data}} \\[-0.03cm]  \textbf{\textsf{\& Initializations}}};
\node[align=center, rotate=90, text = black, scale = 0.55] (c) at (0.68, 1.92) {\textbf{\textsf{Local Problem Data}} \\[-0.03cm]  \textbf{\textsf{\& Initializations}}};
\node[align=center, rotate=90, text = black, scale = 0.55] (c) at (4.69, 1.87) {\textbf{\textsf{Local update}} (4)-(5)};
\node[align=center, rotate=90, text = black, scale = 0.55] (c) at (4.69, 5.18) {\textbf{\textsf{Local update}} (4)-(5)};
\node[align=center, rotate=90, text = black, scale = 0.55] (c) at (6.18, 5.8) {\textbf{\textsf{Local update}} (6)};
\node[align=center, rotate=90, text = black, scale = 0.55] (c) at (6.18, 1.27) {\textbf{\textsf{Local update}} (6)};
\node[align=center, rotate=90, text = black, scale = 0.55] (c) at (6.8, 3.52) {\textbf{\textsf{Global update}} (7)};
\node[align=center, rotate=90, text = black, scale = 0.55] (c) at (7.96, 5.46) {\textbf{\textsf{Dual update}} (8)-(9)};
\node[align=center, rotate=90, text = black, scale = 0.55] (c) at (7.96, 1.53) {\textbf{\textsf{Dual update}} (8)-(9)};
\node[align=center, rotate=90, text = black, scale = 0.55] (c) at (13.74, 1.9) {\textbf{\textsf{Local update}} (4)-(5)};
\node[align=center, rotate=90, text = black, scale = 0.55] (c) at (13.74, 5.2) {\textbf{\textsf{Local update}} (4)-(5)};
\node[align=center, rotate=90, text = black, scale = 0.55] (c) at (14.93, 3.52) {\textbf{\textsf{Global update}} (7)};
\node[align=center, rotate=0, scale = 0.5] (c) at (4.07, 5.27) {$\rho_1^0, \mu_1^0$};
\node[align=center, rotate=0, scale = 0.45] (c) at (4.07, 2.1) {$\rho_N^0, \mu_N^0$};
\node[align=center, scale = 0.5] (c) at (2.82, 5.5) {\textbf{\textsf{Feed-forward}} \\
\textbf{\textsf{params}}};
\node[align=center, scale = 0.5] (c) at (2.82, 4.59) {\textbf{\textsf{Feedback}} \\ \textbf{\textsf{FC layer}}};
\node[align=center, scale = 0.5] (c) at (2.82, 2.3) {\textbf{\textsf{Feed-forward}} \\
\textbf{\textsf{params}}};
\node[align=center, scale = 0.5] (c) at (2.82, 1.41) {\textbf{\textsf{Feedback}} \\ \textbf{\textsf{FC layer}}};
\node[align=center, rotate=0, scale = 0.5] (c) at (5.75, 4.47) {$\bar{\alpha}^0$};
\node[align=center, rotate=0, scale = 0.5, text = RoyalBlue3] (c) at (1.6, 5.43) {\textbf{\textsf{Local}} \\ \textbf{\textsf{residuals}}};
\node[align=center, rotate=0, scale = 0.5, text = RoyalBlue3] (c) at (1.6, 2.2) {\textbf{\textsf{Local}} \\ \textbf{\textsf{residuals}}};
\node[align=center, rotate=0, scale = 0.47] (c) at (5.28, 5.37) {$\vx_1^1, \vz_1^1$};
\node[align=center, rotate=0, scale = 0.43] (c) at (5.3, 1.7) {$\vx_N^1, \vz_N^1$};
\node[align=center, rotate=0, scale = 0.47] (c) at (5.45, 3.88) {$\vx_2^1, \vz_2^1$};
\node[align=center, rotate=0, scale = 0.45] (c) at (5.36, 3.28) {$\vx_{N-1}^1, \vz_{N-1}^1$};
\node[align=center, rotate=0, scale = 0.47] (c) at (7.45, 6.0) {$\vs_1^1$};
\node[align=center, rotate=0, scale = 0.47] (c) at (7.45, 1.0) {$\vs_N^1$};
\node[align=center, rotate=0, scale = 0.5] (c) at (7.3, 3.7) {$\vw^1$};
\node[align=center, rotate=0, scale = 0.47] (c) at (7.45, 5.35) {$\tilde{\vw}_1^1$};
\node[align=center, rotate=0, scale = 0.47] (c) at (7.45, 1.76) {$\tilde{\vw}_N^1$};
\node[align=center, rotate=0, scale = 0.45] (c) at (8.53, 5.6) {$\vy_1^1, \blambda_1^1$};
\node[align=center, rotate=0, scale = 0.42] (c) at (8.55, 1.75) {$\vy_N^1, \blambda_N^1$};
\node[align=center, rotate=0, scale = 0.5] (c) at (13.03, 5.4) {$\rho_1^{K-1}$ \\[0.1cm] $\mu_1^{K-1}$};
\node[align=center, rotate=0, scale = 0.5] (c) at (13.03, 2.29) {$\rho_N^{K-1}$ \\[0.1cm] $\mu_N^{K-1}$};
\node[align=center, rotate=0, scale = 0.5, text = RoyalBlue3] (c) at (10.47, 5.41) {\textbf{\textsf{Local}} \\ \textbf{\textsf{residuals}}};
\node[align=center, rotate=0, scale = 0.5, text = RoyalBlue3] (c) at (10.47, 2.26) {\textbf{\textsf{Local}} \\ \textbf{\textsf{residuals}}};
\node[align=center, rotate=0, scale = 0.5] (c) at (14.93, 5.05) {$\bar{\alpha}^{K-1}$};
\node[align=center, rotate=0, scale = 0.47] (c) at (14.41, 5.4) {$\vx_1^K, \vz_1^K$};
\node[align=center, rotate=0, scale = 0.47] (c) at (14.42, 1.73) {$\vx_N^K, \vz_N^K$};
\node[align=center, rotate=0, scale = 0.47] (c) at (13.7, 3.83) {$\vx_2^K, \vz_2^K$};
\node[align=center, rotate=0, scale = 0.47] (c) at (13.52, 3.27) {$\vx_{N-1}^K, \vz_{N-1}^K$};
\node[align=center, scale = 0.5] (c) at (11.7, 5.5) {\textbf{\textsf{Feed-forward}} \\
\textbf{\textsf{params}}};
\node[align=center, scale = 0.5] (c) at (11.7, 4.6) {\textbf{\textsf{Feedback}} \\ \textbf{\textsf{FC layer}}};
\node[align=center, scale = 0.5] (c) at (11.7, 2.39) {\textbf{\textsf{Feed-forward}} \\
\textbf{\textsf{params}}};
\node[align=center, scale = 0.5] (c) at (11.7, 1.47) {\textbf{\textsf{Feedback}} \\ \textbf{\textsf{FC layer}}};
\node[align=center, rotate=0, scale = 0.5] (c) at (15.45, 3.7) {$\vw^K$};
\node[align=center, rotate=90, text = black, scale = 0.55] (c) at (15.95, 3.52) {\textbf{\textsf{Training Loss}}};
\end{tikzpicture}
\caption{\textbf{The DeepDistributedQP architecture.} The proposed framework relies on unrolling the DistributedQP optimizer as a supervised deep learning framework. In particular, we interpret its iterations (\ref{eq: DistrQP x update})-(\ref{eq: DistrQP y update}) as sequential network layers and introduce learnable components ({\color{DarkOrange1} orange blocks}) to facilitate reaching the desired accuracy after a predefined number of allowed iterations.}
\label{fig: deep distr qp arch}
\end{figure}
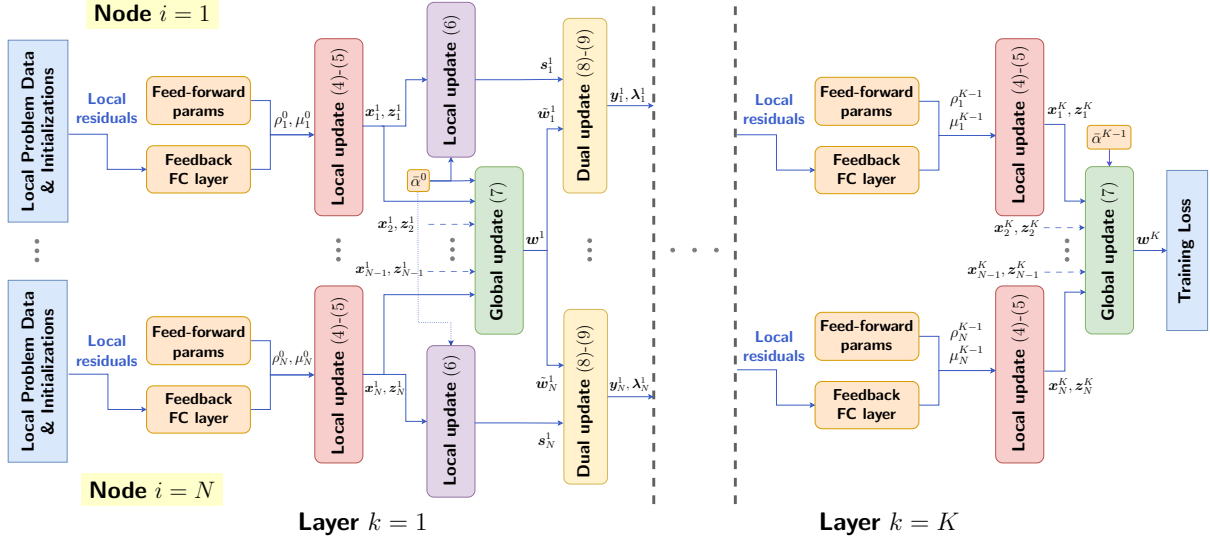

\subsection{Main Architecture}
\label{sec: distr deep qp - main arch}
\paragraph{Architecture overview.} 
The \textbf{\textsf{DeepDistributedQP}} architecture arises from unrolling the DistributedQP optimizer within the supervised learning paradigm. (Fig. \ref{fig: deep distr qp arch}). This is accomplished through treating the updates (\ref{eq: DistrQP x update})-(\ref{eq: DistrQP lambda update}) as blocks in sequential layers of a deep learning network. The number of layers is equal to the predefined number of allowed iterations $K$, with each layer corresponding to an iteration $k = 1,\dots,K$. The inputs of the network are the local problem data $\zeta_i$ and initializations $\vx_i^0$, $\vz_i^0$, $\tilde{\vw}_i^0$, $\vs_i^0$, $\blambda_i^0$ and $\vy_i^0$. These are initially passed to $N$ parallel local blocks corresponding to (\ref{eq: DistrQP x update})-(\ref{eq: DistrQP z update}), which output the new variables $\vx_i^1$ and $\vz_i^1$. Then, all $\vz_i^1$ are fed into $N$ new parallel local blocks (\ref{eq: DistrQP s update}), yielding the new iterates $\vs_i^1$. In the meantime, all $\vx_i^1$ are communicated to a central node that computes the new iterate $\vw^1$ through the weighted averaging step (\ref{eq: DistrQP global update}). Subsequently, the global variable components $\tilde{\vw}_i$ are communicated back to each local node $i$, to perform the updates (\ref{eq: DistrQP lambda update})-(\ref{eq: DistrQP y update}) which output the updated dual variables $\blambda_i, \vy_i$. This group of blocks is then repeated $K$ times, yielding the output of the network which is the final global variable iterate $\vw^K$. 
\paragraph{Learning feedback policies.} Standard deep unfolding typically leverages data to learn algorithm parameters tailored for a specific problem \citep{shlezinger2022model}. From a control theoretic point of view, this process can be interpreted as seeking \textit{open-loop} policies without the incorporating any feedback. In our setup, this would be equivalent with learning the optimal parameters $\bar{\rho}_i^k$, $\bar{\mu}_i^k$, $\bar{\alpha}^k$
\begin{align}
\rho_i^k = \mathrm{SoftPlus}({\bar{\rho}_i^k}), \quad 
& \mu_i^k = \mathrm{SoftPlus}(\bar{\mu}_i^k), \quad 
\alpha^k = \mathrm{Sigmoid}_{1,2}(\bar{\alpha}^k), 
\end{align}
for all $i = 1,\dots, N$ and $k = 1,\dots, K$, where the $\mathrm{SoftPlus}(\cdot)$ function is used  to guarantee the positivity of $\rho_i^k$, $\mu_i^k$, and the sigmoid function $\mathrm{Sigmoid}_{1,2}(\cdot)$ restricts each $\alpha^k$ between $(1,2)$.

In the meantime, the predominant practice for online adaptation of the ADMM penalty parameters relies on observing the primal and dual residuals every few iterations \citep{boyd2011distributed}. The widely-used rule suggests that if the ratio of primal-to-dual residuals is high, the penalty parameter $\rho$ should be increased; conversely, if the ratio is low, $\rho$ should be decreased. 
Despite its heuristic nature, this approach includes a notion of ``feedback'' since the current state of the optimizer is used to adapt the parameters, and as a result, it can be interpreted as a closed-loop policy. Based on this point of view, our goal is to learn the optimal \textit{closed-loop} policies for the local penalty parameters
\begin{equation}
\rho_i^k = \mathrm{SoftPlus}\Big( \bar{\rho}_i^k 
+ \underbrace{\pi_{i,\rho}^k (r_{i, \rho}^k, s_{i, \rho}^k; \theta_{i, \rho}^k)}_{\hat{\rho}_i^k} \Big), \quad 
\mu_i^k = \mathrm{SoftPlus}\Big( \bar{\mu}_i^k + \underbrace{\pi_{i,\mu}^k (r_{i, \mu}^k, s_{i, \mu}^k; \theta_{i, \mu}^k)}_{\hat{\mu}_i^k}
\Big),
\end{equation}
where $\hat{\rho}_i^k$, $\hat{\mu}_i^k$ are feedback components obtained from policies $\pi_{i, \cdot}^k (r_{i, \cdot}^k, s_{i, \cdot}^k; \theta_{i, \cdot}^k)$, 
parameterized by fully-connected neural network layers with inputs $r_{i, \cdot}^k, s_{i, \cdot}^k$ and weights $\theta_{i,\cdot}^k$. The terms $r_{i, \cdot}^k$ and $s_{i, \cdot}^k$ represent the local primal and dual residuals of node $i$ at layer $k$ and are detailed in Appendix \ref{sec: appendix - feedback residuals}.
\paragraph{Solving the local updates.}
The most computationally demanding block in DeepDistributedQP is solving the local updates (\ref{eq: DistrQP x update}), as this requires solving a linear system of size $n_i + m_i$.
Similar to OSQP \citep{stellato2020osqp}, this can be accomplished using either a direct or an indirect method.
The direct method factors the KKT matrix, solving the system via forward and backward substitution.
This approach is particularly efficient when penalty parameters remain fixed, as the same factorization can then be reused accross iterations. Nevertheless, at larger scales, this factorization might become impractical. In contrast, with the indirect method, we  eliminate $\vnu_i^{k + 1}$ to solve the linear system:
\begin{equation} \label{eq:indirect_method_linear_system}
    \underbrace{(\mQ_i + \mu_i^k \mI + \mA_i^\top \rho_i^k \mA_i)}_{\bar{\mQ}_i^k} \vx_i^{k + 1} = \underbrace{- \vq_i + \mu_i^k \tilde{\vw}_i - \vy_i + \mA_i^\top \rho_i^k \vz_i - \mA_i^\top \blambda_i}_{\bar{\vb}_i^k}.
\end{equation}
This new linear system is solved for $\vx_i^{k + 1}$ using an iterative scheme such as the conjugate gradient (CG) method.
We then substitute $\vnu_i^{k + 1} = \rho_i^k (\mA_i \vx_i^{k + 1} - \vz_i) + \blambda_i$. The indirect method has three important properties that make it particularly attractive in our setup.
First, its computational complexity scales better w.r.t. the dimension of the local problem, while no additional overhead is introduced by changing the penalty parameters.
Second, it can be warmstarted using the solution from the previous iteration, greatly reducing the number of iterations required to converge to a solution.
The final important property, which is critical for the scalability of the DeepDistributedQP, is that training with the indirect method can be much more memory efficient as shown in Section \ref{sec: distr deep qp - backprop}.
\paragraph{Training loss.} Let $\calS = \{ \zeta^j \}_{j = 1}^H$ be a dataset consisting of $H$ problem instances $\zeta^j = \{ (\mQ_i, \vq_i, \mA_i, \vb_i )_{i=1}^N, \vw^* \}_j$ subject to the known mapping $\calG$ of problem (\ref{eq: distr qp problem}). The loss we are using for training is the average of the $\gamma_k$-scaled distances of the global iterates $\vw_1, \dots, \vw_N$ from the known optimal solution $\vw^*$ of each problem instance $\zeta^j$, provided as
\begin{equation}
\ell(\calS; \theta) = \frac{1}{H} \sum_{j = 1}^H 
\sum_{k = 1}^K
\gamma_k \| \vw^k(\zeta^j; \theta) - \vw^*(\zeta^j) \|_2,
\end{equation}
where $\theta$ corresponds to the concatenation of all learnable parameters/weights.
\begin{figure}[t]
\centering
\begin{tikzpicture}
\node[anchor=south west,inner sep=0] at (0,0){\includegraphics[width=1.0\textwidth, trim={0.5cm 0.0cm 0cm 0cm},clip]{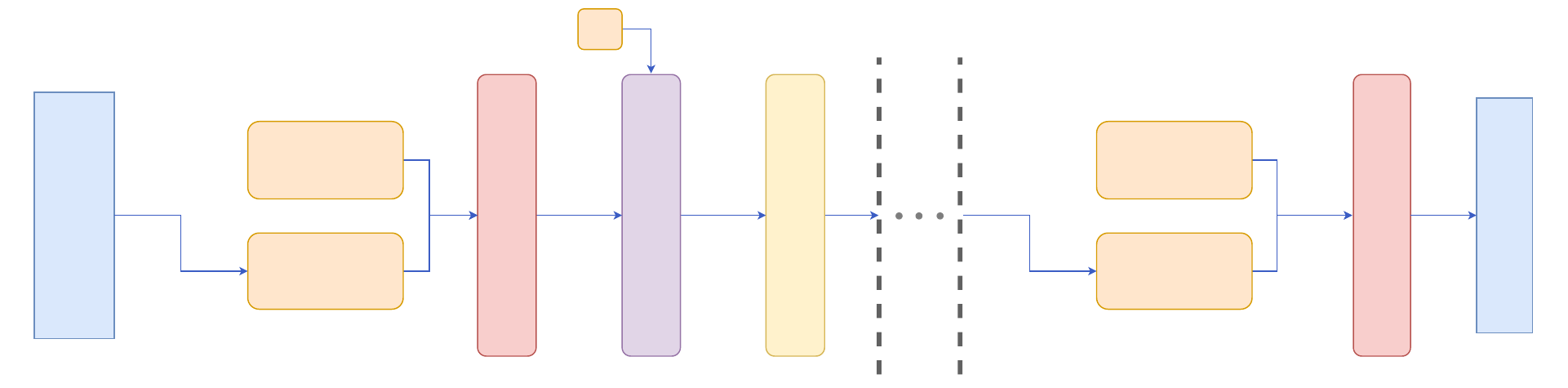}};
\node[align=center, scale = 0.7] (c) at (6.0, -0.1) {\textbf{\textsf{Layer}} $k = 1$};
\node[align=center, scale = 0.7] (c) at (12.6, -0.1) {\textbf{\textsf{Layer}} $k = K$};
\node[align=center, scale = 0.6] (c) at (3.23, 2.43) {\textbf{\textsf{Feed-forward}} \\
\textbf{\textsf{params}}};
\node[align=center, scale = 0.6] (c) at (3.23, 1.23) {\textbf{\textsf{Feedback}} \\ \textbf{\textsf{FC layer}}};
\node[align=center, rotate=0, scale = 0.6, text = RoyalBlue3] (c) at (1.63, 2.1) {\textbf{\textsf{Residuals}}};
\node[align=center, rotate=90, text = black, scale = 0.65] (c) at (5.19, 1.88) {\textbf{\textsf{First update}} (94)-(95)};
\node[align=center, rotate=90, text = black, scale = 0.65] (c) at (6.71, 1.88) {\textbf{\textsf{Second update}} (97)};
\node[align=center, rotate=90, text = black, scale = 0.65] (c) at (8.26, 1.88) {\textbf{\textsf{Dual update}} (98)};
\node[align=center, rotate=90, text = black, scale = 0.65] (c) at (14.54, 1.88) {\textbf{\textsf{First update}} (94)-(95)};
\node[align=center, rotate=0, scale = 0.7] (c) at (5.93, 2.1) {$\vx^1, \vz^1$};
\node[align=center, rotate=0, scale = 0.7] (c) at (7.44, 2.1) {$\vt^1, \vs^1$};
\node[align=center, rotate=0, scale = 0.7] (c) at (8.82, 2.1) {$\blambda^1$};
\node[align=center, rotate=90, text = black, scale = 0.65] (c) at (0.53, 1.85) {\textbf{\textsf{Problem Data}} \\  \textbf{\textsf{\& Initializations}}};
\node[align=center, rotate=90, text = black, scale = 0.65] (c) at (15.85, 1.78) {\textbf{\textsf{Training Loss}}};
\node[align=center, rotate=0, scale = 0.7] (c) at (6.17, 3.85) {$\bar{\alpha}^0$};
\node[align=center, rotate=0, scale = 0.7] (c) at (4.59, 2.1) {$\rho^{0}$};
\node[align=center, rotate=0, scale = 0.6, text = RoyalBlue3] (c) at (10.7, 2.1) {\textbf{\textsf{Residuals}}};
\node[align=center, scale = 0.6] (c) at (12.31, 2.43) {\textbf{\textsf{Feed-forward}} \\
\textbf{\textsf{params}}};
\node[align=center, scale = 0.6] (c) at (12.31, 1.23) {\textbf{\textsf{Feedback}} \\ \textbf{\textsf{FC layer}}};
\node[align=center, rotate=0, scale = 0.7] (c) at (15.15, 2.1) {$\vx^{K}$};
\node[align=center, rotate=0, scale = 0.7] (c) at (13.81, 2.1) {$\rho^{K-1}$};
\end{tikzpicture}
\caption{\textbf{DeepQP:} The centralized version of DeepDistributedQP which boils down to unfolding the standard OSQP method.}
\label{fig: DeepQP}
\end{figure}
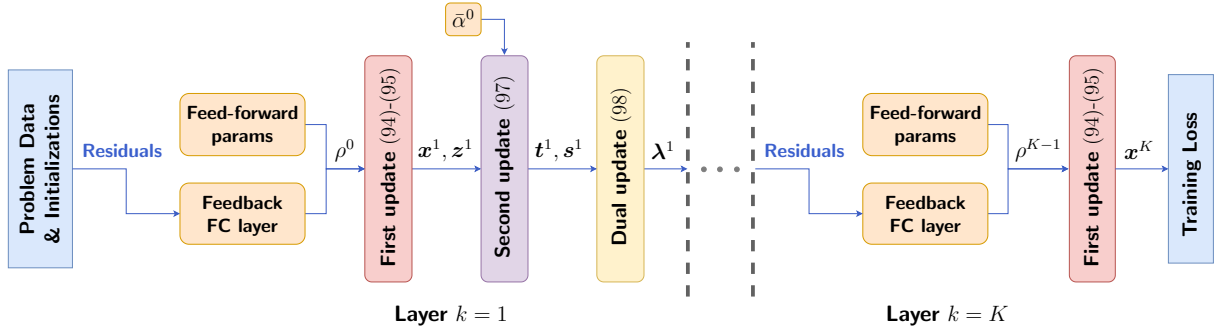

\paragraph{Centralized version.} While this work primarily focuses on distributed optimization, for completeness, we also introduce \textbf{\textsf{DeepQP}}, the centralized version of our framework, for addressing general QPs of the form (\ref{eq: general qp problem}). In the centralized case, our framework simplifies to $N = 1$, eliminating the need for distinguishing between local and global variables. Under this simplification, the DistributedQP optimizer coincides with OSQP. Hence, DeepQP consists of unfolding the OSQP updates (see Appendix \ref{sec: appendix - non distr deep qp}) and learning policies for adapting its penalty and over-relaxation parameters. The resulting framework is illustrated in Fig. \ref{fig: DeepQP}. Additional details on DeepQP are provided in Appendix \ref{sec: appendix - non distr deep qp}.
%

\label{sec: distr deep qp - feedback policies}

\subsection{Implicit Differentiation}
\label{sec: distr deep qp - backprop}

When solving for the local updates in (\ref{eq:indirect_method_linear_system}) using the indirect method, it is computationally intractable to backpropagate through all CG iterations.
This is especially important in the context of unfolding, as it would become necessary to unroll multiple inner CG optimization loops. To address this, we leverage the implicit function theorem (IFT) to express the solution of (\ref{eq:indirect_method_linear_system}) as an implicit function of the local problem data.
This allows us to compute gradients in a manner that avoids unrolling the CG iterations and requires solving a linear system with the same coefficient matrix, but with a new RHS, achieved by rerunning the CG method.
This result is formalized in the following theorem.
\begin{theorem}[Implicit Differentiation of Indirect Method] \label{thm:implicit_differentiation}
Let $\vx_i^{k + 1}$ be the unique solution to the linear system $\bar{\mQ}_i^k \vx_i^{k + 1} = \bar{\vb}_i^k$ in (\ref{eq:indirect_method_linear_system}). Let $\nabla_{\vx} L(\vx_i^{k + 1})$ be a backward pass vector computed through reverse-mode automatic differentiation of some loss function $L$. Then, the gradient of $L$ with respect to $\bar{\mQ}_i^k$ and $\bar{\vb}_i^k$ is given by
\begin{align*}
    \nabla_{\bar{\mQ}_i^k} L &= \frac{1}{2} (\vx_i^{k + 1} \otimes d\vx_i^{k + 1} + d\vx_i^{k + 1} \otimes \vx_i^{k + 1}), \\
    \nabla_{\bar{\vb}_i^k} L &= -d\vx_i^{k + 1},
\end{align*}
where $d\vx_i^{k + 1}$ is the unique solution to the linear system $\bar{\mQ}_i^k d\vx_i^{k + 1} = -\nabla_x L(\vx_i^{k + 1})$.
\end{theorem}
The proof is provided in Appendix \ref{sec: appendix - IFT proof} and is a straightforward application of the IFT, similar to the results established by \citet{amos2017optnet} and \citet{agrawal2019differentiable}.
%

%% file: sections/sec5_gen_bounds.tex
\section{Generalization Bounds}
\label{sec: gen bounds}

In this section, we establish guarantees on the expected performance of DeepDistributedQP. To achieve this, we leverage the PAC-Bayes framework \citep{alquier2024user}, a well-known statistical learning methodology for providing bounds on expected loss metrics that hold with high probability. In our case, we provide bounds on the \textit{expected progress} of the final iterate $\vw^K$ towards reaching the optimal solution $\vw^*$ for unseen problems drawn from the same distribution as the training dataset.
\vspace{-0.2cm}
\paragraph{Learning stochastic policies.} PAC-Bayes theory is applicable to frameworks that learn weight distributions rather than fixed weights. For this reason, in order to establish such guarantees, we switch to learning a Gaussian distribution of weights $\mathcal{P} = \calN(\mu_{\Theta}, \Sigma_{\Theta})$ based on a prior $\mathcal{P}_0 = \calN(\mu_{\Theta}^0, \Sigma_{\Theta}^0)$. 
This choice is motivated by the fact that PAC-Bayes bounds include Kullback–Leibler (KL) divergence terms which can be easily evaluated and optimized for Gaussian distributions. 
%
%
%
%
%
\vspace{-0.2cm}
\paragraph{Generalization bound for DeepDistributedQP.} To facilitate the exhibition of our performance guarantees, we provide necessary preliminaries on PAC-Bayes theory in Appendix \ref{ssec:pac bayes background}. To establish a generalization guarantee for DeepDistributedQP, a meaningful loss function must first be selected. This quantity will be denoted $q(\zeta; \theta)$ to differentiate from the loss used for training. To capture the progress the optimizer makes towards optimality, we propose the following \emph{progress metric}:
\begin{align}
q(\zeta; \theta) = \min \left\{ \frac{ \| \vw^K(\zeta; \theta) - \vw^*(\zeta) \|_2}{ \| \vw^0(\zeta) - \vw^*(\zeta)\|_2}, 1 \right\}.
\end{align}
This loss function measures progress by comparing the distance between the final iterate $\vw^K(\zeta; \theta)$ and problem solution $\vw^*(\zeta)$ with the distance between the initialization $\vw^0(\zeta; \theta)$ and the solution. This choice satisfies the requirement of being bounded between 0 and 1 while being more informative than the indicator losses used in prior work that simply determine whether the final iterate is within a specified neighborhood of the optimal solution \citep{sambharya2024data}. Moreover, this loss is invariant to the scale of the problem data since it is a relative measurement.

As in Section \ref{ssec:pac bayes background}, let $q_{\calD}(\calP)$ be the true expected loss and $q_{\calS}(\calP)$ the empirical expected loss. To evaluate the PAC-Bayes bounds in (\ref{eq:pac bayes}), the expectation $\E_{\theta \sim \mathcal{P}}[q(\zeta; \theta)]$ must be computed as part of the definition of $q_{\calS}(\calP)$. Since no closed-form solution is available, an empirical estimate using $M$ sampled weights $(\theta_i)_{i = 1}^M$ is required to upper bound $q_{\calS}(\calP)$ with high probability. We adopt a standard approach involving a sample convergence bound (\cite{majumdar2021pac}, \cite{dziugaite2017computing}, \cite{langford2001not}). Specifically, define the empirical estimate of $q_{\calS}(\calP)$ as:
\begin{align}
\hat{q}_{\calS}(\mathcal{P}; M) = \frac{1}{M H} \sum_{i = 1}^H \sum_{j = 1}^M q(\zeta_i; \theta_j). 
\end{align}
Then, the following sample convergence bound provides an upper bound on $q_{\calS}(\calP)$,
\begin{align}
q_{\calS}(\calP) \leq \bar{q}_{\calS}(\calP; M, \epsilon) \coloneqq \sD_{\textsc{KL}}\left(\hat{q}_{\calS}(\mathcal{P}; M)\ \|\ M^{-1}\log\left(2/\epsilon\right)\right),
\label{eq:sample convergence}
\end{align}
with probability $1 - \epsilon$.
The following theorem summarizes the PAC-Bayes bound we use to evaluate the generalization capabilities of our framework.

\begin{theorem}[Generalization bound for DeepDistributedQP]
\label{thm:gen bound}
For problems $\zeta \in \calZ$ drawn from distribution $\calD$, the true expected progress metric of DeepDistributedQP with policy $\calP$, i.e.,
\begin{equation}
q_{\calD}(\calP) = \E_{\zeta \sim \calD} ~ \E_{\theta \sim \mathcal{P}} 
\left[
\min \left \{ \frac{ \| \vw^K(\zeta; \theta) - \vw^*(\zeta) \|_2}{ \| \vw^0(\zeta) - \vw^*(\zeta) \|_2}, 1 \right\}
\right],
\end{equation}
is bounded with probability at least $1 - \delta - \epsilon$ by:
\begin{equation}
q_{\calD}(\calP) \leq
\sD_{\textsc{KL}}^{-1}
\left( 
\bar{q}_\calS(\calP; M, \epsilon)
\middle\|
\left(\sD_{\textsc{KL}}(\calP \| \calP_0) + \log (2 \sqrt{H} / \delta) \right)/H
\right),
\end{equation}
where $\bar{q}_{\calS}(\calP; M, \epsilon)$ is the estimate of $q_\calS(\calP; M, \epsilon)$ described in (\ref{eq:sample convergence}).
\end{theorem}

We explain in detail how we train for optimizing this bound in Appendix \ref{sec: opt for gen bounds}.
\vspace{-0.2cm}

%% file: sections/sec6_experiments.tex
\section{Experiments}
\begin{figure}[t]
\centering
\includegraphics[width=\textwidth, trim={0cm 0cm 0cm 0cm},clip]{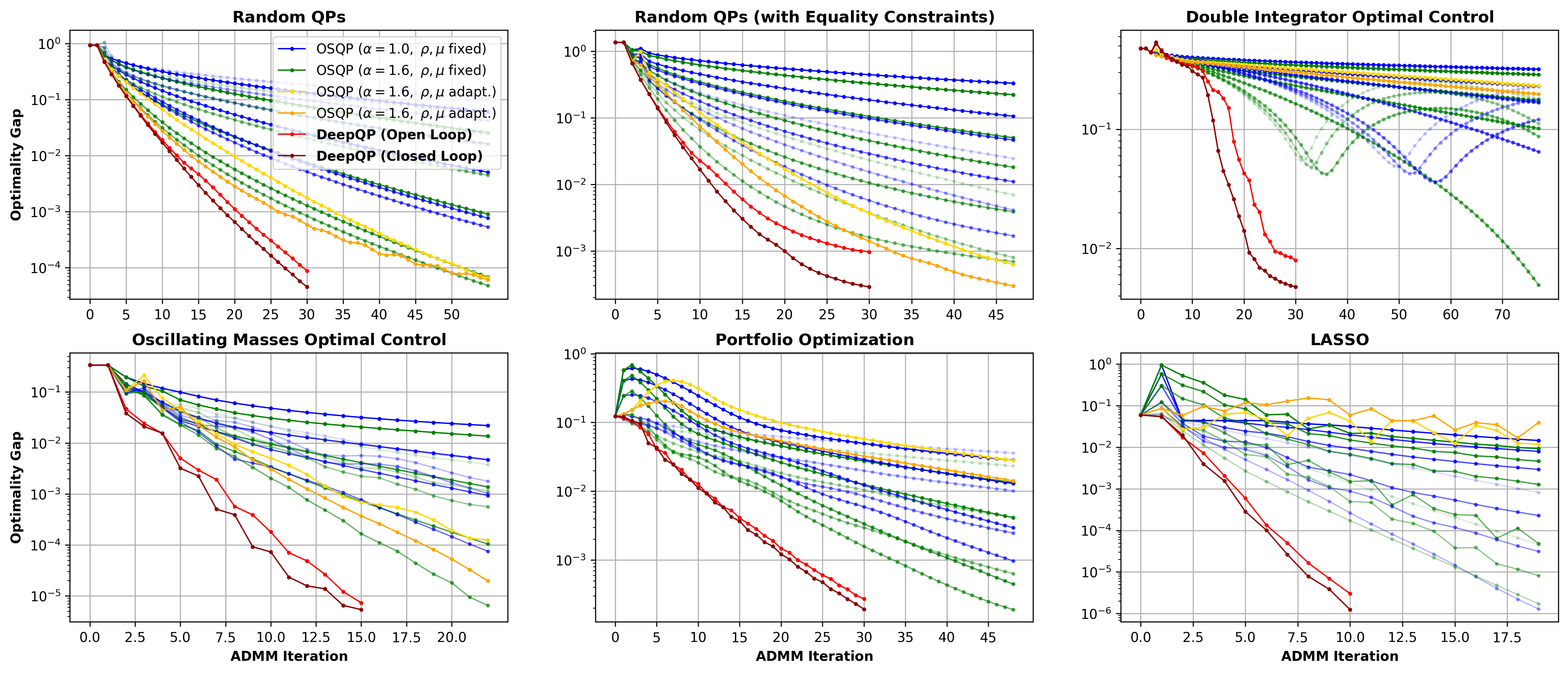}
\caption{\textbf{Small-scale centralized comparison between DeepQP and OSQP.} Across all tested problems, DeepQP consistently outperforms OSQP (same per-iteration complexity using the indirect method).}
\label{fig: results deep qp}
\end{figure}
We conduct extensive experiments to highlight the effectiveness, scalability and generalizability of the proposed methods. Section \ref{sec: results centr} shows the advantageous performance of DeepQP against OSQP on a variety of centralized QPs. In Section \ref{sec: results distr}, we address large-scale problems, showcasing the scalability of DeepDistributedQP despite being trained exclusively on much lower-dimensional instances. Additionally, we discuss the advantages of learning local policies over shared ones and evaluate the proposed generalization bounds, which provide guarantees for the performance of our framework on unseen problems.
An overall discussion and potential limitations are provided in Section \ref{sec: results disc}. All experiments were performed on an system with an RTX 4090 GPU 24GB, a 13th Gen Intel(R) Core(TM) i9-13900K and 64GB of RAM.

\subsection{Small-Scale Centralized Experiments: DeepQP vs OSQP}
\label{sec: results centr}

\paragraph{Setup.} We begin with comparing DeepQP against OSQP for solving centralized QPs (\ref{eq: general qp problem}). The following problems are considered: i,ii) random QPs without/with equality constraints, iii, iv) optimal control for double integrator and oscillating masses, v) portfolio optimization, and vi) LASSO regression. For all problems, we set a maximum allowed amount of iterations $K$ for DeepQP within $[10,30]$ and examine how many iterations OSQP requires to reach the same accuracy. We train DeepQP using both open-loop and closed-loop policies and with a dataset of size $H \in [500,2000]$. For OSQP, we consider both constant and adaptive penalty parameters $\rho$ and we set $\alpha$ to be either $1.0$ or  $1.6$. Additional details on DeepQP, OSQP and the problems can be found in Appendix \ref{sec: appendix - exp details}.
\paragraph{Performance comparison.} 
The comparison between DeepQP and OSQP is illustrated in Fig. \ref{fig: results deep qp}. Note that both methods share the same per-iteration complexity from solving (\ref{eq: osqp - first update - indirect}). We evaluate their performance by comparing the (normalized) optimality gap $\| \vx^k - \vx^* \|_2 / \sqrt{n}$. For all tested problems, DeepQP provides a consistent improvement over OSQP, requiring $1.5-3$ times fewer iterations to reach the desired accuracy. Furthermore, the advantage of incorporating feedback in the policies is shown, as closed-loop policies outperform open-loop ones in all cases. 
%
%
\subsection{Large-Scale Distributed Experiments: Scaling DeepDistributedQP}
\label{sec: results distr}
\paragraph{Setup.} The purpose of the following analysis is to compare the performance and scalability of DeepDistributedQP (ours), DistributedQP (ours) and OSQP for large-scale QPs of the form (\ref{eq: distr qp problem}). We consider the following six problems: i,ii) random networked QPs without/with equality constraints, iii, iv) multi-agent optimal control for coupled pendulums and oscillating masses, v) network flow, and vi) distributed LASSO. We select a maximum allowed number of iterations $K$ for DeepDistributedQP within $[20,50]$ and examine what is the computational effort required by DistributedQP and OSQP to achieve the same accuracy measured by the optimality gap $\| \vw^k - \vw^* \|_2 / \sqrt{n}$.  More details about our experimental setup are provided in Appendix \ref{sec: appendix - exp details}.
%
%
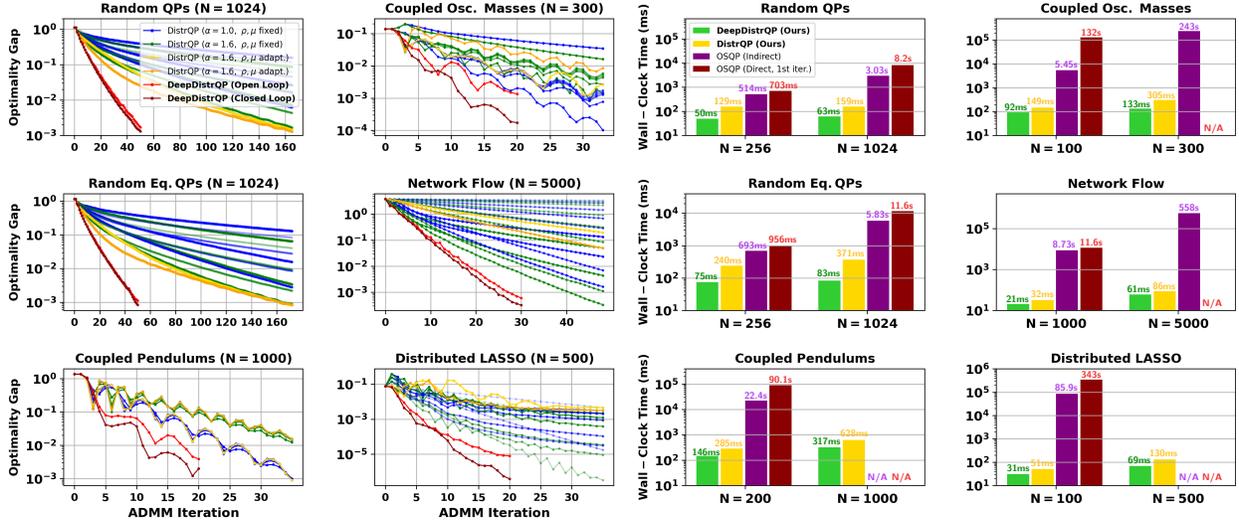
\begin{figure}[t]
\centering
\begin{tikzpicture}
\node[anchor=south west,inner sep=0] at (0,0){\includegraphics[width=1.0\textwidth, trim={0cm 0cm 0cm 0cm},clip]{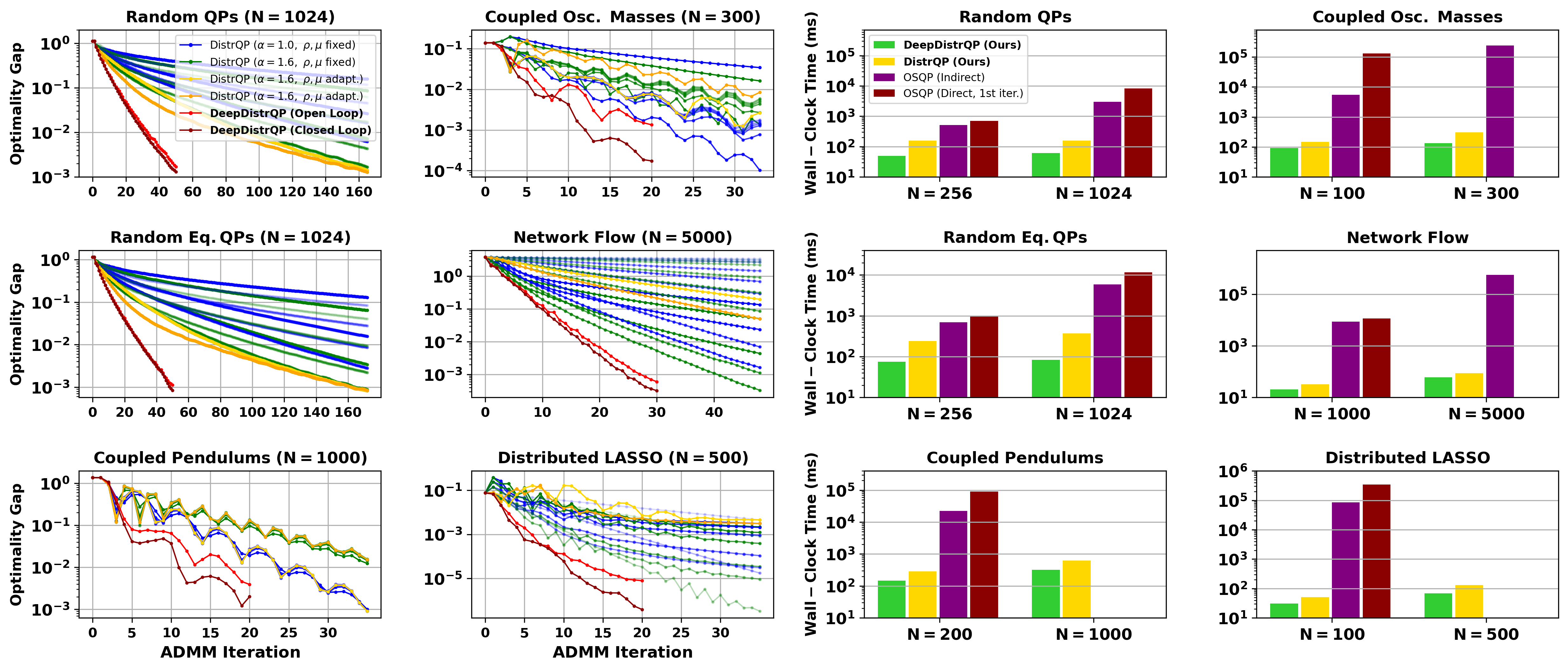}
};
\node[align=center, scale = 0.3, color = Green4] (c) at (9.38, 5.5) {\textbf{\textsf{50ms}}};
\node[align=center, scale = 0.3, color = Goldenrod1] (c) at (9.67, 5.67) {\textbf{\textsf{129ms}}};
\node[align=center, scale = 0.3, color = DarkOrchid2] (c) at (10.0, 5.83) {\textbf{\textsf{514ms}}};
\node[align=center, scale = 0.3, color = Firebrick2] (c) at (10.41, 5.87) {\textbf{\textsf{703ms}}};
\node[align=center, scale = 0.3, color = Green4] (c) at (11.01, 5.53) {\textbf{\textsf{63ms}}};
\node[align=center, scale = 0.3, color = Goldenrod1] (c) at (11.29, 5.67) {\textbf{\textsf{159ms}}};
\node[align=center, scale = 0.3, color = DarkOrchid2] (c) at (11.65, 6.08) {\textbf{\textsf{3.03s}}};
\node[align=center, scale = 0.3, color = Firebrick2] (c) at (12.0, 6.23) {\textbf{\textsf{8.2s}}};
\node[align=center, scale = 0.3, color = Green4] (c) at (9.38, 3.33) {\textbf{\textsf{75ms}}};
\node[align=center, scale = 0.3, color = Goldenrod1] (c) at (9.67, 3.55) {\textbf{\textsf{240ms}}};
\node[align=center, scale = 0.3, color = DarkOrchid2] (c) at (10.0, 3.75) {\textbf{\textsf{693ms}}};
\node[align=center, scale = 0.3, color = Firebrick2] (c) at (10.4, 3.82) {\textbf{\textsf{956ms}}};
\node[align=center, scale = 0.3, color = Green4] (c) at (11.01, 3.38) {\textbf{\textsf{83ms}}};
\node[align=center, scale = 0.3, color = Goldenrod1] (c) at (11.3, 3.64) {\textbf{\textsf{371ms}}};
\node[align=center, scale = 0.3, color = DarkOrchid2] (c) at (11.66, 4.16) {\textbf{\textsf{5.83s}}};
\node[align=center, scale = 0.3, color = Firebrick2] (c) at (11.98, 4.27) {\textbf{\textsf{11.6s}}};
\node[align=center, scale = 0.3, color = Green4] (c) at (9.38, 1.0) {\textbf{\textsf{146ms}}};
\node[align=center, scale = 0.3, color = Goldenrod1] (c) at (9.68, 1.12) {\textbf{\textsf{285ms}}};
\node[align=center, scale = 0.3, color = DarkOrchid2] (c) at (10.04, 1.75) {\textbf{\textsf{22.4s}}};
\node[align=center, scale = 0.3, color = Firebrick2] (c) at (10.37, 1.96) {\textbf{\textsf{90.1s}}};
\node[align=center, scale = 0.3, color = Green4] (c) at (10.98, 1.14) {\textbf{\textsf{317ms}}};
\node[align=center, scale = 0.3, color = Goldenrod1] (c) at (11.35, 1.24) {\textbf{\textsf{628ms}}};
\node[align=center, scale = 0.3, color = DarkOrchid2] (c) at (11.66, 0.65) {\textbf{\textsf{N/A}}};
\node[align=center, scale = 0.3, color = Firebrick2] (c) at (11.98, 0.65){\textbf{\textsf{N/A}}};
\node[align=center, scale = 0.3, color = Green4] (c) at (13.5, 5.59) {\textbf{\textsf{92ms}}};
\node[align=center, scale = 0.3, color = Goldenrod1] (c) at (13.82, 5.67) {\textbf{\textsf{149ms}}};
\node[align=center, scale = 0.3, color = DarkOrchid2] (c) at (14.17, 6.15) {\textbf{\textsf{5.45s}}};
\node[align=center, scale = 0.3, color = Firebrick2] (c) at (14.5, 6.575) {\textbf{\textsf{132s}}};
\node[align=center, scale = 0.3, color = Green4] (c) at (15.1, 5.63) {\textbf{\textsf{133ms}}};
\node[align=center, scale = 0.3, color = Goldenrod1] (c) at (15.43, 5.75) {\textbf{\textsf{305ms}}};
\node[align=center, scale = 0.3, color = DarkOrchid2] (c) at (15.8, 6.66) {\textbf{\textsf{243s}}};
\node[align=center, scale = 0.3, color = Firebrick2] (c) at (16.13, 5.3) {\textbf{\textsf{N/A}}};
\node[align=center, scale = 0.3, color = Green4] (c) at (13.52, 3.03) {\textbf{\textsf{21ms}}};
\node[align=center, scale = 0.3, color = Goldenrod1] (c) at (13.85, 3.1) {\textbf{\textsf{32ms}}};
\node[align=center, scale = 0.3, color = DarkOrchid2] (c) at (14.14, 3.76) {\textbf{\textsf{8.73s}}};
\node[align=center, scale = 0.3, color = Firebrick2] (c) at (14.5, 3.79) {\textbf{\textsf{11.6s}}};
\node[align=center, scale = 0.3, color = Green4] (c) at (15.13, 3.17) {\textbf{\textsf{61ms}}};
\node[align=center, scale = 0.3, color = Goldenrod1] (c) at (15.47, 3.21) {\textbf{\textsf{86ms}}};
\node[align=center, scale = 0.3, color = DarkOrchid2] (c) at (15.8, 4.25) {\textbf{\textsf{558s}}};
\node[align=center, scale = 0.3, color = Firebrick2] (c) at (16.11, 2.97) {\textbf{\textsf{N/A}}};
\node[align=center, scale = 0.3, color = Green4] (c) at (13.53, 0.78) {\textbf{\textsf{31ms}}};
\node[align=center, scale = 0.3, color = Goldenrod1] (c) at (13.85, 0.85) {\textbf{\textsf{51ms}}};
\node[align=center, scale = 0.3, color = DarkOrchid2] (c) at (14.17, 1.85) {\textbf{\textsf{85.9s}}};
\node[align=center, scale = 0.3, color = Firebrick2] (c) at (14.5, 2.025) {\textbf{\textsf{343s}}};
\node[align=center, scale = 0.3, color = Green4] (c) at (15.14, 0.89) {\textbf{\textsf{69ms}}};
\node[align=center, scale = 0.3, color = Goldenrod1] (c) at (15.46, 0.98) {\textbf{\textsf{130ms}}};
\node[align=center, scale = 0.3, color = DarkOrchid2] (c) at (15.8, 0.65) {\textbf{\textsf{N/A}}};
\node[align=center, scale = 0.3, color = Firebrick2] (c) at (16.11, 0.65) {\textbf{\textsf{N/A}}};
\end{tikzpicture}
\vspace{-0.6cm}
\caption{\textbf{Scaling DeepDistributedQP to high-dimensional problems.} Left: Comparison between DeepDistributedQP and its traditional optimization counterpart DistributedQP (same per-iteration complexity). Right: Total wall-clock time required by DeepDistributedQP, DistributedQP and OSQP (using indirect or direct method) to achieve the same accuracy.}
\vspace{-0.4cm}
\label{fig: distr results}
\end{figure}
%

\paragraph{Training on low-dimensional problems.} One of the key advantages of DeepDistributedQP is that it only requires using small-scale problems for training. The training dimensions for each problem are detailed in Table \ref{tab: training testing}. Both open-loop and closed-loop versions are trained using shared policies on datasets of size $H \in [500,1000]$. We employ the shared policies version of DeepDistributedQP to enable the same policies to be applied to larger problems during testing.

\begin{figure}
\includegraphics[width=\textwidth, trim={0.0cm 0cm 0.0cm 0.0cm},clip]{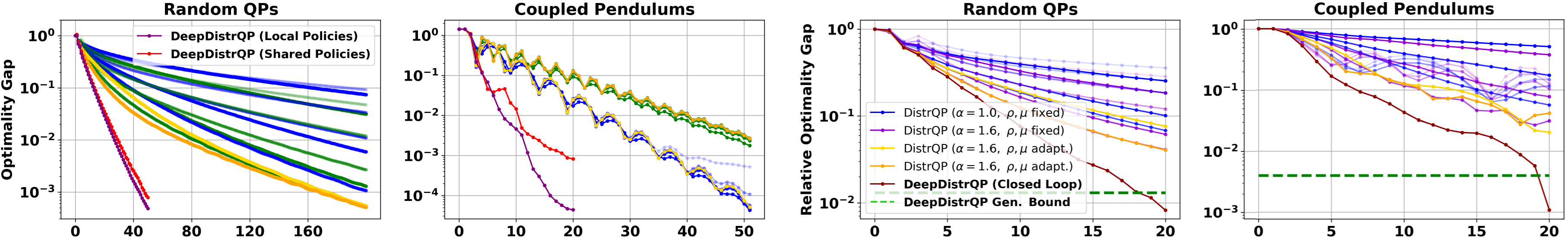}
\label{Case 1 fig b}
\caption{\textbf{Left: Local vs shared policies.} We showcase the advantage of learning local policies over shared ones. \textbf{Right: Performance guarantees.} The obtained generalization bounds guarantee the performance of DeepDistributedQP and its improvements over its standard optimization counterpart DistributedQP.}
\label{fig: local and gen bounds}
\end{figure}

\paragraph{Scaling to high-dimensional problems.} 
Subsequently, we evaluate DeepDistributedQP on problems with significantly larger scale than those used during training. The maximum problem dimensions tested are shown in Table \ref{tab: training testing}. On the left side of Fig. \ref{fig: distr results}, we highlight the superior performance of DeepDistributedQP over its standard optimization counterpart DistributedQP (same per-iteration complexity). In all cases, the learned algorithm achieves the same level of accuracy while requiring 1.5-3.5 times fewer iterations. Additionally, the right side of Fig. \ref{fig: distr results} compares the total wall-clock time between DeepDistributedQP, DistributedQP and OSQP (using indirect or direct method). For a complete illustration, we refer the reader to Table \ref{tab: wall-clock-appendix} in Appendix \ref{sec: appendix wall clock}. The provided results emphasize the superior scalability of the two proposed distributed methods against OSQP for large-scale QPs, as well as the advantage of our deep learning-aided approach over traditional optimization.
\paragraph{Local vs shared policies.} When applying a policy to a problem with the same dimensions as used during training, leveraging local policies instead of shared ones can be advantageous for better exploiting the structure of the problem. On the left side of Fig. \ref{fig: local and gen bounds} ,we compare the performance of local and shared policies on random QPs and coupled pendulums. For the coupled pendulums problem, which exhibits significant underlying structure, local policies demonstrate clear superiority. For the random QPs problem, where structural patterns are less pronounced, the advantage of local policies is smaller but still significant.
\paragraph{Performance guarantees.}
Next, we verify the guarantees of our framework for generalizing on unseen random QP $(N = 16)$ and coupled pendulums $(N = 10)$ problems. We switch from learning deterministic weights to learning stochastic ones and follow the procedure described in Appendix \ref{sec: opt for gen bounds} with $H = 15000$ training samples, $M = 30000$ sampled weights for the bounds evaluation, $\delta = 0.009$ and $\epsilon = 0.001$. The resulting generalization bounds, illustrated in Fig. \ref{fig: local and gen bounds} (right), are expressed in terms of the the expected final relative optimality gap - the progress metric used for deriving bounds in Section \ref{sec: gen bounds}, implying that with $99 \%$ probability the average performance of our framework will be bounded by this threshold. The bounds are observed to be tight compared to actual performance, underscoring their significance. Moreover, they outperform the standard optimizers, providing a strong guarantee of improved performance for DeepDistributedQP.

\begin{table}[t]
\begin{center}
\scalebox{0.87}{
\begin{tabular}{ccccccccc}
\toprule
  & \multicolumn{4}{c}{\textbf{Training}}  & \multicolumn{4}{c}{\textbf{Max Testing}}
\\[0.1cm]
\textbf{Problem Class} & $N$ & $n$ & $m$ & $\texttt{nnz}(\mQ, \mA)$ & $N$ & $n$ & $m$ & $\texttt{nnz}(\mQ, \mA)$
\\ \hline
Random QPs 
& 16 & 160 & 120 & 4,000 
& 1,024 & 10,240 & 9,920 & 300,800
\\ 
Random QPs w/ Eq. Constr. 
& 16 & 160 & 168 & 4,960
& 1,024 & 10,240 & 9,920 & 300,800
\\ 
Coupled Pendulums & 10 &  470 & 640 & 3,690 
& 1,000 & 47,000 & 64,000 & 380,880
\\ 
Coupled Osc. Masses & 10 & 470 & 1,580 & 4,590 & 300 & 28,200 & 47,400 & 141,180
\\ 
Network Flow & 20 & 100 & 140 & 600 & 5,000 & 25,000 & 35,000 & 150,000
\\ 
Distributed LASSO & 10 & 1,100 & 3,000 & 29,000 & 500 & 50,100 & 150,000 & 1,450,000
\\ \bottomrule
\end{tabular}
}
\end{center}
\caption{\textbf{Training and maximum testing dimensions for DeepDistributedQP.} The metric $\texttt{nnz}(\mQ, \mA)$ denotes the total number of non-zero elements in $\mQ$ and $\mA$.}
\label{tab: training testing}
\end{table}



\subsection{Discussion}
\label{sec: results disc}

%

\paragraph{In which cases can we use the direct method?} As illustrated in Fig. \ref{fig: distr results} and Table \ref{tab: wall-clock-appendix}, and further discussed in \cite{stellato2020osqp}, the indirect method is generally preferred for solving systems of the form (\ref{eq: DistrQP x update}) - or (\ref{eq: osqp - first update - direct}) for DeepQP/OSQP - once the problem reaches a certain scale. In this work, we adopt this approach both for training, due the memory and computational advantages outlined in Section \ref{sec: distr deep qp - backprop}, and evaluating DeepDistributedQP/DeepQP. However, it is worth considering whether the direct method might be advantageous during evaluation, a choice that depends on the problem scale and capabilities of the available hardware. Overall, the results of this work show that learning policies for the algorithm parameters is significantly beneficial in the context of both distributed and centralized QP assuming the indirect method is used. In future work, we wish to also explore schemes that adapt the parameters less frequently using the direct method and/or designing mechanisms to dynamically switch between the two approaches.

\paragraph{Limitations.} One limitation of the proposed framework is its reliance on a supervised training loss, requiring a dataset of pre-solved problems. In future work, we aim to explore training through directly minimizing the problem residuals rather than the optimality gaps. Furthermore, while PAC-Bayes theory provides an important probabilistic bound on average performance, stronger guarantees may be necessary for safety-critical applications to ensure reliability and robustness.

%% file: sections/sec7_conclusion.tex
\section{Conclusion and Future Work}
\vspace{-0.2cm}
In this work, we introduced DeepDistributedQP, a new deep learning-aided distributed optimization architecture for solving large-scale QP problems. The proposed method relies on unfolding the iterations of a novel optimizer named DistributedQP as layers of a supervised deep learning framework. The expected performance of our learned optimizer on unseen problems is also theoretically established through PAC-Bayes theory. DeepDistributedQP exhibits impressive scalability in effectively tackling large-scale optimization problems while being trained exclusively on much smaller ones. In addition, both DeepDistributedQP and Distributed significantly outperform OSQP in terms of required wall-clock time to reach the same accuracy as dimension increases. Furthermore, we showcase that the proposed PAC-Bayes bounds provide meaningful practical guarantees for the performance of DeepDistributedQP on new problems. 

In future work, we wish to extend the proposed framework to a semi-supervised version that relies less on pre-solved problems for training. In addition, we wish to explore incorporating more complex learnable components such as LSTMs for feedback within our framework. Finally, we wish to consider other classes of distributed constrained optimization methods outside of quadratic programming. 

%% file: sections/appendix.tex




\section{Complete Derivation of DistributedQP Algorithm}
\label{sec: appendix - distr_qp derivation}

\paragraph{Problem transformation and augmented Lagrangian.} Here, we present the detailed derivation of the DistributedQP algorithm presented in Section \ref{sec: distr qp - the algo}. We consider the over-relaxed version of ADMM \citep{boyd2011distributed} with $\alpha \in [1,2)$. First, let us rewrite problem (\ref{eq: distr qp problem}) as 
\begin{equation}
\min_{\vx} ~ \sum_{i \in \calV} \frac{1}{2} \vx_i^\top \mQ_i \vx_i + \vq_i^\top \vx_i
\quad \mathrm{s.t.} 
\quad \mA_i \vx_i = \vz_i,
~ \vz_i \leq \vb_i, 
~ \vx_i = \tilde{\vw}_i, 
\quad i \in \calV.
\label{eq: distr qp problem v2}
\end{equation}
where we have introduced the auxiliary variables $\vz_i$ for each $i = 1,\dots,N$. In addition, let us define the variables $\vs_i, \ i = 1, \dots, N$, and rewrite problem (\ref{eq: distr qp problem v2}) as
\begin{equation}
\begin{aligned}
& \min_{\vx} ~ \sum_{i \in \calV} \frac{1}{2} \vx_i^\top \mQ_i \vx_i 
+ \vq_i^\top \vx_i
+ \mathcal{I}_{\mA_i \vx_i = \vz_i} (\vx_i, \vz_i)
\\
& ~~~ \mathrm{s.t.} 
\quad 
\vz_i = \vs_i, 
~ \vs_i \leq \vb_i, 
~ \vx_i = \tilde{\vw}_i, 
\quad i \in \calV.
\end{aligned}
\label{eq: distr qp problem v3}
\end{equation}
The above splitting constitutes the problem suitable for being addressed with a two-block ADMM scheme, where the first block of variables consists of $\{\vx_i, \vz_i\}_{i = 1,\dots,N}$, and the second one contains the variables $\{\vs_i\}_{i = 1,\dots,N}$ and $\vw$. The (scaled) augmented Lagrangian (AL) for problem (\ref{eq: distr qp problem v3}) is given by
\begin{equation}
\begin{aligned}
\mathcal{L} 
& = \sum_{i \in \calV} 
\frac{1}{2} \vx_i^\top \mQ_i \vx_i 
+ \vq_i^\top \vx_i
+ \mathcal{I}_{\mA_i \vx_i = \vz_i} (\vx_i, \vz_i) 
+ \mathcal{I}_{\vs_i \leq \vb_i} (\vs_i) 
\\
& 
\qquad + \frac{\rho_i}{2} \left\| \vz_i - \vs_i + \frac{\blambda_i}{\rho_i} \right\|_2^2
+ \frac{\mu_i}{2} \left\| \vx_i - \tilde{\vw}_i + \frac{\vy_i}{\mu_i} \right\|_2^2.
\end{aligned}
\label{eq: AL distributed}
\end{equation}
\paragraph{First block of ADMM primal updates.} The first block of variables is updated through 
\begin{equation*}
\{\vx_i, \vz_i\}_{i \in \calV} = \argmin \mathcal{L}(\vx, \vz, \vs^k, \vw^k, \blambda^k, \vy^k).
\end{equation*}
This minimization can be decoupled to the following $N$ subproblems for each $i \in \calV$,
\begin{align*}
\{ \vx_i, \vz_i \} = \argmin
\frac{1}{2} \vx_i^\top \mQ_i \vx_i 
& + \vq_i^\top \vx_i
+ \frac{\rho_i}{2} \left\| \vz_i - \vs_i + \frac{\blambda_i}{\rho_i} \right\|_2^2
+ \frac{\mu_i}{2} \left\| \vx_i - \tilde{\vw}_i + \frac{\vy_i}{\mu_i} \right\|_2^2
\\[0.2cm]
& \quad \mathrm{s.t.} \quad \mA_i \vx_i = \vz_i,
\end{align*}
where we have temporarily dropped the superscript iteration indices for convenience.
Since these problems are equality constrained QPs, we can obtain a closed-form solution. The KKT conditions for each subproblem are given by
\begin{subequations}
\begin{align}
\mQ_i \vx_i + \vq_i
+ \mu_i ( \vx_i - \tilde{\vw}_i ) + \vy_i
+ \mA_i^\top \vnu_i & = \vzero,
\\
\rho_i ( \vz_i - \vz_i ) + \blambda_i - \vnu_i & = \vzero,
\\
\mA_i \vx_i - \vz_i & = \vzero,
\end{align}
\end{subequations}
where $\vnu_i$ is the Lagrange multiplier corresponding to the constraint $\mA_i \vx_i = \vz_i$. Eliminating $\vz_i$ leads to the following system of equations
\begin{equation}
\begin{bmatrix}
\mQ_i + \mu_i \mI & \mA_i^\top \\
\mA_i & - 1 / \rho_i \mI 
\end{bmatrix}
\begin{bmatrix}
\vx_i^{k+1} \\
\vnu_i^{k+1}
\end{bmatrix}
= 
\begin{bmatrix}
- \vq_i + \mu_i \tilde{\vw}_i^k - \vy_i^k
\\
\vz_i - 1 / \rho_i \blambda_i^k
\end{bmatrix},
\end{equation}
with $\vz_i^{k+1}$ given by
\begin{equation}
\vz_i^{k+1} =  \vs_i^k + \rho_i^{-1} (\vnu_i^{k+1} - \blambda_i^k).
\end{equation}
\paragraph{Second block of ADMM primal updates.}  The second block of updates is given by
\begin{equation*}
\{ \vs_i \}_{i \in \mathcal{V}}, \vw = \argmin \mathcal{L}(\vx^{k+1}, \vz^{k+1}, \vs, \vw, \blambda^k, \vy^k),
\end{equation*}
or more analytically
\begin{align*}
\{ \vs_i \}_{i \in \mathcal{V}}, \vw & = \argmin \sum_{i \in \calV} 
\frac{\rho_i}{2} \left\| \alpha \vz_i^{k+1} + (1 - \alpha) \vs_i^k - \vs_i + \frac{\blambda_i^k}{\rho_i} \right\|_2^2
\\
& \qquad \qquad \qquad \qquad \qquad + \frac{\mu_i}{2} \left\| \alpha \vx_i^{k+1} + (1 - \alpha) \tilde{\vw}_i^k - \tilde{\vw}_i + \frac{\vy_i^k}{\mu_i} \right\|_2^2
\quad \mathrm{s.t.} \quad 
\vs_i \leq \vb_i .
\end{align*}
Note that this minimization can be decoupled w.r.t. all $\vs_i$, $i \in \calV$, and $\vw$. In particular, each $\vs_i$ is updated in parallel through 
\begin{equation}
\vs_i^{k+1} = \Pi_{\vs_i \leq \vb_i}
\left( \alpha \vz_i^{k+1} + (1 - \alpha^k) \vs_i^k + \blambda_i^k / \rho_i \right).
\end{equation}
The global variable $\vw$ minimization can also be decoupled among its components $l = 1, \dots, n$, which gives
\begin{equation*}
\vw_l = \argmin \sum_{\calG(i,j) = l}
\frac{\mu_i}{2} \left\| \alpha [\vx_i^{k+1}]_j + (1 - \alpha) [\tilde{\vw}_i^k]_j - [\tilde{\vw}_i]_j + \frac{[\vy_i^k]_j}{\mu_i} \right\|_2^2.
\end{equation*}
Setting the gradient to be equal to zero gives
\begin{equation*}
\sum_{\calG(i,j) = l}
\mu_i \left[ \alpha [\vx_i^{k+1}]_j + (1 - \alpha) \vw_l^k - \vw_l^{k+1} + \frac{[\vy_i^k]_j}{\mu_i} \right]
= \vzero,
\end{equation*}
leading to
\begin{equation*}
\sum_{\calG(i,j) = l}
\mu_i \vw_l^{k+1}
= 
\sum_{\calG(i,j) = l}
\mu_i \left[ \alpha [\vx_i^{k+1}]_j + (1 - \alpha) \vw_l^k + \frac{[\vy_i^k]_j}{\mu_i} \right]
\end{equation*}
which eventually gives the update rule
\begin{equation}
\vw_l^{k+1} = \frac{\sum_{\mathcal{G}(i,j) = l} \alpha \mu_i [\vx_i^{k+1}]_j + [\vy_i^k]_j}{\sum_{\mathcal{G}(i,j) = l} \mu_i}
+ (1 - \alpha) \vw_l^k.
\end{equation}
\paragraph{ADMM dual updates.} Finally, the dual variables are updated through dual ascent steps as follows
\begin{align}
\blambda_i^{k+1} & = \blambda_i^k + \rho_i 
( \alpha \vz_i^{k+1} 
+ (1 - \alpha) \vs_i^k
- \vs_i^{k+1})
\\
\vy_i^{k+1} & = \vy_i^k + \mu_i 
( \alpha \vx_i^{k+1} 
+ (1 - \alpha) \tilde{\vw}_i^k
- \tilde{\vw}_i^{k+1}).
\end{align}
\paragraph{Simplifying the global update.} It is important to observe that after the first iteration, the global update can be simplified to 
\begin{equation}
\vw_l^{k+1} = \alpha \frac{\sum_{\mathcal{G}(i,j) = l} \mu_i [\vx_i^{k+1}]_j}{\sum_{\mathcal{G}(i,j) = l} \mu_i}
+ (1 - \alpha) \vw_l^k,
\end{equation}
since the summation 
\begin{align}
\sum_{\mathcal{G}(i,j) = l} [\vy_i^{k+1}]_j
& = 
\sum_{\mathcal{G}(i,j) = l} [\vy_i^k]_j
+ \mu_i ( \alpha [\vx_i^{k+1}]_j 
+ (1 - \alpha) [\tilde{\vw}_i^k]_j
- [\tilde{\vw}_i^{k+1}]_j)
\nonumber
\\
& = 
\sum_{\mathcal{G}(i,j) = l} [\vy_i^k]_j
+ \mu_i ( \alpha [\vx_i^{k+1}]_j 
+ (1 - \alpha) \vw_l^k
- \vw_l^{k+1})
\nonumber
\\
& = \sum_{\mathcal{G}(i,j) = l} [\vy_i^k]_j
+ \mu_i \bigg[ \alpha [\vx_i^{k+1}]_j 
+ \cancel{\purple{(1 - \alpha) \vw_l^k}}
\nonumber
\\
& \qquad \qquad - \frac{\sum_{\mathcal{G}(u,v) = l} \alpha \mu_u [\vx_u^{k+1}]_v + [\vy_u^k]_v}{\sum_{\mathcal{G}(u,v) = l} \mu_u}
- \cancel{\purple{(1 - \alpha) \vw_l^k}} \bigg]
\nonumber
\\
& 
= \sum_{\mathcal{G}(i,j) = l} [\vy_i^k]_j
+ \mu_i \bigg[ \alpha [\vx_i^{k+1}]_j 
- \frac{\sum_{\calG(u,v) = l} \alpha \mu_u [\vx_u^{k+1}]_v + [\vy_u^k]_v}{\sum_{\mathcal{G}(u,v) = l} \mu_u}
\bigg]
\nonumber
\\
& 
= \sum_{\mathcal{G}(i,j) = l} [\vy_i^k]_j
+ \alpha \mu_i [\vx_i^{k+1}]_j 
- \frac{\purple{\cancel{\sum_{\mathcal{G}(i,j) = l} \mu_i}} 
 \left[ \sum_{\mathcal{G}(u,v) = l} \alpha \mu_u [\vx_u^{k+1}]_v + [\vy_u^k]_v \right]}{\purple{\cancel{\sum_{\mathcal{G}(u,v) = l} \mu_u}}}
 \nonumber
\\
& 
= \sum_{\mathcal{G}(i,j) = l} [\vy_i^k]_j
+ \alpha \mu_i [\vx_i^{k+1}]_j  
- \sum_{\mathcal{G}(u,v) = l} \alpha \mu_u [\vx_u^{k+1}]_v + [\vy_u^k]_v = 0.
\label{eq: sum y zeros}
\end{align}

\section{Standard Convergence Guarantees for Simplified Version of DistributedQP}
\label{sec: appendix - distr_qp standard guarantees}

In the simplified case where $\rho_i^k = \rho$, $\mu_i^k = \mu$ for all $i \in \calV$ and for all $k$, as well as $\alpha^k = 1$ for all $k$, it would be straightforward to apply the classical convergence guarantees of two-block ADMM for convex optimization problems \citep{deng2016global} to ensure the convergence of DistributedQP. In the following, we show how DistributedQP fits under this setup.

Let us define the variables $\bar{\vx} = [\{\vx_i\}_{i \in \calV}; \{\vz_i\}_{i \in \calV}]$ and $\bar{\vz} = [\{ \vs_i\}_{i \in \calV}; \vw]$. Then, we can rewrite problem (\ref{eq: distr qp problem v3}) as 
\begin{equation}
\min f(\bar{\vx}) + g(\bar{\vz}) \quad \mathrm{s.t.} \quad 
\bar{\mA} \bar{\vx} + \bar{\mB} \bar{\vz} = \bar{\vc},
\end{equation}
where 
\begin{equation}
f(\bar{\vx}) = \sum_{i \in \calV} 
\frac{1}{2} \vx_i^\top \mQ_i \vx_i 
+ \vq_i^\top \vx_i
+ \mathcal{I}_{\mA_i \vx_i = \vz_i} (\vx_i, \vz_i), 
\quad
g(\bar{\vz}) = \sum_{i \in \calV} \mathcal{I}_{\vs_i \leq \vb_i} (\vs_i),
\end{equation}
and $\bar{\mA} = \mathrm{bdiag}(\mI, \mI)$, $\bar{\mB} = \mathrm{bdiag}(\mI, \mG)$ and $\vc = \vzero$, with $\mG \in \sR^{(\sum_i n_i) \times n}$ defined such that $\vx = \mG \vw$. In other words, $\mG$ is the matrix that represents the local-to-global variable components mapping, formally defined as $\mG = [\mG_1; \dots; \mG_N]$ with each submatrix $\mG_i \in \sR^{n_i \times n}$ given by
\begin{equation}
[\mG_i]_{u,v} = 
\begin{cases}
1,  & \text{if} ~ v = \calG(i, v)
\\ 
0, & \text{else}
\end{cases}.
\end{equation}
Given this representation, it becomes clear that our algorithm can be framed as a two-block ADMM. Now, note that $\mG$ is a full column rank matrix since all global variable components $\vg_l$ are mapped to at least one local variable component $[\vx_i]_j$. Then, since the functions $f,g$ are convex and the matrices $\bar{\mA}, \bar{\mB}$ are full column rank, it follows from \cite{deng2016global} that the algorithm is guaranteed to converge to the optimal solution.

Nevertheless, this analysis would have only been applicable to this simplified case of the proposed DistributedQP algorithm. In Appendix \ref{sec: appendix - distr_qp convergence proof}, we tackle the more complex case of iteration-varying relaxation and local penalty parameters. 

\section{Proof of DistributedQP Asymptotic Convergence}
\label{sec: appendix - distr_qp convergence proof}

In this section, we prove that DistributedQP is guaranteed to converge to optimality, even in the more challenging case of iteration-varying relaxation and local penalty parameters. The following analysis extends the theoretical results presented in \cite{xu2017adaptive}, where the convergence of an adaptive relaxed variant of two-block ADMM is provided. Nevertheless, this analysis is not directly applicable to our case which involves distinct local penalty parameters per computational node. 

\subsection{Sketch of Proof}
\label{sec: appendix - proof sketch}

To begin, we outline the following conventions. The points $\vx^*, \vz^*, \vs^{*}, \vw^{*}, \vy^{*}, \blambda^{*}$ are the KKT points of problem (\ref{eq: distr qp problem v3}). We refer to the notion of a distance function at any $(k+1)^{th}$ iteration to be representing a weighted squared norm of the difference between the variables $\vs^{k+1}, \vw^{k+1}, \vy^{k+1}, \blambda^{k+1}$ and their corresponding optimal values $\vs^{*}, \vw^{*}, \vy^{*}, \blambda^{*}$, indicating the distance from the optimal point. 

We prove the convergence of in the following steps:
\begin{itemize}
    \item First, we will derive a descent relation (\ref{Distr OSQP - descent relation final}), which establishes a relationship between the values of the distance function for consecutive iterations. To derive the descent relation in Lemma \ref{Lemma 4}, we first introduce the relations (\ref{R1})-(\ref{R8}) in Lemmas \ref{lemma 1}-\ref{lemma 3}. 
    \item Next, we use the derived descent relation to prove the convergence in Section \ref{Proof of DistrQP convergence} based on Assumption \ref{assumption: eta}.
\end{itemize}

\subsection{Necessary Lemmas}
\label{sec: appendix - proof lemmas}

Here, we present some necessary lemmas before proving the convergence of DistributedQP in Section \ref{Proof of DistrQP convergence}. For notational convenience, let us define
\begin{equation*}
f_i(\vx_i) = \frac{1}{2} \vx_i^\top \mQ_i \vx_i + \vq_i^\top \vx_i, \quad 
\calC_i = \{ \vs_i | \vs_i \leq \vb_i \}, \quad i \in \calV. 
\label{problem simplified notation}
\end{equation*}
%
%
\begin{lemma} \label{lemma 1}
For all $i \in \calV$, the following four relationships hold at every iteration $k$:
\begin{align}
& \sum_{i \in \calV} \mG_i^\top \vy_i^{k+1} = \vzero,
\tag{R1}
\label{R1}
\\
& \alpha^k \vx_i^{k+1} =
\frac{1}{\mu_i^k}(\vy_i^{k+1} - \vy_i^k) 
- (1 - \alpha^k) \mG_i \vw^k
+ \mG_i \vw^{k+1},
\tag{R2}
\label{R2}
\\
& \alpha^k \vz_i^{k+1} =
\frac{1}{\rho_i^k} (\blambda_i^{k+1}-\blambda_i^k)   
- (1 - \alpha^k) \vs_i^k
+ \vs_i^{k+1},
\tag{R3}
\label{R3}
\\
& \blambda_i^k{}^\top(\vt_1 - \vt_2) = 0, \quad \forall ~ \vt_1, \vt_2 \in \calC_i.
\tag{R4}
\label{R4}
\end{align}
\begin{proof}
Relationship (\ref{R1}) is equivalent with the argument proved in (\ref{eq: sum y zeros}). Indeed, if we observe that each matrix $\mG_i^\top \in \sR^{n \times n_i}$ indicates the mapping from local indices $(i,j)$ to global indices $l$ for a particular $i$, then we can write
\begin{equation}
\sum_{i \in \calV} \mG_i^\top \vy_i^{k+1} =  
\begin{bmatrix}
\sum_{\calG(i,j) = 1} [\vy_i^{k+1}]_j
\\
\vdots
\\
\sum_{\calG(i,j) = n} [\vy_i^{k+1}]_j
\end{bmatrix}
= 
\vzero,
\end{equation}
which yields (\ref{R1}).
Relationship (\ref{R2}) follows by rearranging the dual update (\ref{eq: DistrQP y update}) and replacing $\tilde{\vw}_i = \mG_i \vw$. Similarly, relationship (\ref{R3}) follows by rearranging the dual update (\ref{eq: DistrQP lambda update}). In the remaining, we focus on proving (\ref{R4}). Let us repeat the $\vs_i$-update (\ref{eq: DistrQP s update}) as
\begin{equation}
\vs_i^{k+1} = \Pi_{\calC_i}
\left( \alpha^k \vz_i^{k+1} + (1 - \alpha^k) \vs_i^k + \blambda_i^k / \rho_i^k \right).
\label{Distr OSQP l1 eq1}
\end{equation}
We define the closed convex cone $\bar{\calC}_i = \{ \vp | \; \vp \leq 0 \}$,
such that (\ref{Distr OSQP l1 eq1}) is rewritten as
%
\begin{equation}
    \vs_i^{k+1} = \Pi_{\bar{\calC}_i} \left( \hat{\vs}^{k+1}_i \right)
    + \vb_i,
    \label{Distr OSQP conv - s update relaxed}
\end{equation}
with 
\begin{equation}
    \hat{\vs}^{k+1}_i = \alpha^k \vz_i^{k+1} + (1 - \alpha^k) \vs_i^k + \blambda_i^k / \rho_i^k - \vb_i.
    \label{distr OSQP hat s definition}
\end{equation}
Next, let us rearrange the dual update (\ref{eq: DistrQP lambda update}) as
\begin{align}
    \blambda_i^{k+1} = \rho_i^k (\blambda_i^k/\rho_i^k + 
\alpha^k \vz_i^{k+1} 
+ (1 - \alpha^k) \vs_i^k
- \vs_i^{k+1}),
\end{align}
which can be rewritten through (\ref{distr OSQP hat s definition}) as
\begin{equation}
    \blambda_i^{k+1} = \rho_i^k ( \hat{\vs}^{k+1}_i + \vb_i
- \vs_i^{k+1}  )
\end{equation}
Substituting (\ref{Distr OSQP conv - s update relaxed}) in the above, we get
\begin{align}
    \blambda_i^{k+1} = \rho_i^k ( \hat{\vs}^{k+1}_i 
- \Pi_{\bar{\calC}_i} \left( \hat{\vs}^{k+1}_i \right)).
\label{Distr OSQP conv updates eq2}
\end{align}
For convenience, let us also repeat the definition of polar cones.
\begin{definition}[Polar cones]
Two cone sets $\calD$ and $\calD^o$ are called polar cones if for any $\vd \in \calD$ and $\bar{\vd} \in \calD^o$, if follows that $\vd^\top \bar{\vd} = 0$.  
\end{definition}
\noindent 
By Moreau's decomposition - refer to Theorems 1.1 and 1.2 from \cite{SOLTAN201945} - $\hat{\vs}^{k+1}_i$ can then be expressed as
\begin{equation}
    \hat{\vs}^{k+1}_i = \Pi_{\bar{\calC}_i} \left( \hat{\vs}^{k+1}_i \right)
    + \Pi_{\bar{\calC}_i^o} \left( \hat{\vs}^{k+1}_i \right),
    \label{Moreau decomposition}
\end{equation}
where $\bar{\calC}_i^o$ is a polar cone to $\bar{\calC}_i$. Thus, using (\ref{Distr OSQP conv updates eq2}) and (\ref{Moreau decomposition}), we get
$\blambda_i^{k+1} = 
\rho_i^k \Pi_{\bar{\calC}_i^o} \left( \hat{\vs}^{k+1}_i \right)$,
which implies that
$\blambda_i^{k+1}/\rho_i^k \in \bar{\calC}_i^o$.
Further, since $\bar{\calC}_i^o$ is a cone, and $\rho_i^k>0$, we get
\begin{align}
    \blambda_i^{k+1} \in \bar{\calC}_i^o.
    \label{lambda projection step}
\end{align}
Now, any vector $\vt \in \calC_i$ satisfies 
$\vt - \vb_i \in \bar{\calC}_i$.
Since $\bar{\calC}_i$ and $\bar{\calC}_i^o$ are polar cones, and using (\ref{lambda projection step}), the following relation holds true by the definition of polar cones,
\begin{equation*}
    \blambda_i^{k+1}{}^\top(\vt - \vb_i) = 0 \quad \text{for all } \vt \in \calC_i.
\end{equation*}
Thus, for any vectors $\vt_1, \vt_2 \in \calC_i$ and for all $k$, we have
\begin{align}
    \blambda_i^{k+1}{}^\top(\vt_1 - \vt_2)
    = \blambda_i^{k+1}{}^\top(\vt_1 - \vb_i - (\vt_2 - \vb_i))
    =0,
\end{align}
which proves (\ref{R4}).
\end{proof}
\end{lemma}

\begin{lemma} \label{lemma 2}
For all $i \in \calV$, the following two relationships hold at every iteration $k$:
\begin{align}
&
\big( \nabla f_i(\vx_i^*) +  \vy^*_i \big)^\top 
(\vx_i^* - \vx_i^{k+1}) + 
 \blambda_i^*{}^\top (\vz_i^* - \vz_i^{k+1}) = 0,
 \tag{R5}
 \label{R5}
\\[0.2cm]
&     
\bigg[ \nabla f_i(\vx_i^{k+1}) 
+ \vy_i^{k+1}
+ \mu_i^k \bigg(
(1 - \alpha^k)\vx_i^{k+1}
- (2 - \alpha^k)\mG_i \vw^k
+ \mG_i \vw^{k+1} \bigg) 
\bigg]^\top
(\vx_i^{k+1} - \vx_i^*)
\nonumber
\\
& \qquad 
+ \Bigg[ 
\blambda_i^{k+1} 
+ \rho_i^k \bigg( 
(1 - \alpha^k) \vz_i^{k+1} 
- (2 - \alpha^k) \vs_i^k
+ \vs_i^{k+1}
\bigg)
\bigg]^\top (\vz_i^{k+1} - \vz_i^{*})
= 0.
 \tag{R6}
 \label{R6}
\end{align}
\begin{proof}
We start with proving (\ref{R5}) using the KKT conditions for problem (\ref{eq: distr qp problem v3}). The point $(\vx^*, \vz^*, \vs^*, \vw^*)$ is the optimum of (\ref{eq: distr qp problem v3}) if and only if the following conditions are true:
\begin{subequations}
\begin{align}
    & \text{Optimality for $\vx_i$:}
    && \nabla f_i(\vx_i^*) + \mA_i^\top \vnu_i^* +  \vy^*_i = \vzero
    \label{Distr OSQP conv - optimal point KKT -xtilde}
    \\
    & \text{Optimality for $\vz_i$:} 
    && - \vnu_i^* + \blambda_i^* = \vzero
    \label{Distr OSQP conv - optimal point KKT -ztilde}
    \\
    & \text{Optimality for $\vs_i$:}
    && \blambda_i^* \in \calN_{\calC_i}(\vs_i^*)
    \; \Leftrightarrow \;
    \blambda_i^*{}^\top (\vs_i - \vs_i^*) \leq 0 \; \forall \;
    \vs_i \in \calC_i
    \label{Distr OSQP conv - optimal point KKT -z}
    \\
    & \text{Optimality for $\vw$:}
    && 
    \sum_{i \in \calV} \mG_i^\top \vy_i^{*} = \vzero
    \label{Distr OSQP conv - optimal point KKT - g}
    \\
    & \text{Constraints feasibility: }
    && \tilde{\vz}_i^* = \vs_i^* 
    \label{Distr OSQP conv - optimal point KKT z constraint}
    \\ &
    && \vx_i^* = \mG_i \vw^*
    \label{Distr OSQP conv - optimal point KKT g constraint} 
    \\ &
    && \mA_i \vx_i^* = \vz_i
    \label{Distr OSQP conv - optimal point KKT local constraint} 
    \\ &
    && \vs_i \in \calC_i 
    \label{Distr OSQP conv - optimal point KKT z in set}
\end{align}
\label{KKT conditions optimal}%
\end{subequations}
From (\ref{Distr OSQP conv - optimal point KKT -xtilde}), we have
\begin{equation}
\big( \nabla f_i(\vx_i^*) + \mA_i^\top \vnu_i^* + \vy^*_i \big)^\top 
(\vx_i^* - \vx_i^{k+1}) = 0,
\label{eq: lemma 2 proof first eq}
\end{equation}
and similarly from (\ref{Distr OSQP conv - optimal point KKT -ztilde}), we get
\begin{equation}
\big( - \vnu_i^* + \blambda_i^* \big)^\top (\vz_i^* - \vz_i^{k+1}) = 0.
\label{eq: lemma 2 proof second eq}
\end{equation}
Adding (\ref{eq: lemma 2 proof first eq}) and (\ref{eq: lemma 2 proof second eq}), we get
\begin{equation*}
\big( \nabla f_i(\vx_i^*) + \mA_i^\top \vnu_i^* + \vy^*_i \big)^\top 
(\vx_i^* - \vx_i^{k+1})
+ 
\big( - \vnu_i^* + \blambda_i^* \big)^\top (\vz_i^* - \vz_i^{k+1}) 
= 0,
\end{equation*}
which yields
\begin{equation}
\big( \nabla f_i(\vx_i^*) +  \vy^*_i \big)^\top 
(\vx_i^* - \vx_i^{k+1}) + 
\blambda_i^*{}^\top (\vz_i^* - \vz_i^{k+1})
+ \vnu_i^*{}^\top \big( \mA_i(\vx_i^* - \vx_i^{k+1}) - (\vz_i^* - \vz_i^{k+1}) \big)
= 0.
\label{eq: lemma 2 proof third eq}
\end{equation}
Using (\ref{Distr OSQP conv - optimal point KKT local constraint}) and the fact that $\mA_i \vx_i^{k+1} - \vz_i^{k+1} = \vzero$, we can then simplify (\ref{eq: lemma 2 proof third eq}) to
\begin{equation}
\big( \nabla f_i(\vx_i^*) + \vy^*_i \big)^\top 
(\vx_i^* - \vx_i^{k+1}) + 
\blambda_i^*{}^\top (\vz_i^* - \vz_i^{k+1}) = 0.
\end{equation}
which proves (\ref{R5}).

Subsequently, we proceed with proving relationship (\ref{R6}). The KKT conditions for the $(k+1)$-th update of $\vx_i, \vz_i$ are given by
\begin{subequations}
\begin{align}
    & \text{Optimality for $\vx_i$:}
    && \nabla f_i(\vx_i^{k+1}) 
    + \mA_i^\top \vnu_i^{k+1} 
    + \mu_i^k ( \vx_i^{k+1} - \mG_i \vw^k + \vy_i^k / \mu_i^k ) = \vzero
    \label{Distr OSQP conv - k+1 KKT -xtilde}
    \\
    & \text{Optimality for $\vz_i$:}
    && - \vnu_i^{k+1} 
    + \rho_i^k (  
    \vz_i^{k+1} - \vs_i^k + \blambda_i^k / \rho_i^k )= \vzero
    \label{Distr OSQP conv - k+1 KKT -ztilde}
    \\
    & \text{Constraints feasibility: }
    && \mA_i \vx_i^{k+1} = \vz_i^{k+1}
\end{align}
\label{KKT condition for k+1 subproblem}
\end{subequations}
From (\ref{Distr OSQP conv - k+1 KKT -xtilde}), we have
\begin{align}
\Big[ \nabla f_i(\vx_i^{k+1}) 
    + \mA_i^\top \vnu_i^{k+1} 
    + \mu_i^k ( \vx_i^{k+1} - \mG_i \vw^k + \vy_i^k / \mu_i^k )  
\Big]^\top
(\vx_i^{k+1} - \vx_i^*)
= 0.
\label{Distr OSQP conv - kkt relations eq2}
\end{align}
We rewrite the term $\mu_i^k ( \vx_i^{k+1} - \mG_i \vw^k + \vy_i^k / \mu_i^k )$ using (\ref{eq: DistrQP y update}) as follows
\begin{align*}
& \mu_i^k ( \vx_i^{k+1} - \mG_i \vw^k + \vy_i^k / \mu_i^k ) =
\nonumber
\\
& \qquad =
\mu_i^k \bigg( \vx_i^{k+1} - \mG_i \vw^k + \vy_i^{k+1} / \mu_i^k  -  
\Big( \alpha^k \vx_i^{k+1} 
+ (1 - \alpha^k) \mG_i \vw^k
- \mG_i \vw^{k+1} \Big) \bigg) 
\nonumber
\\
& \qquad =
\vy_i^{k+1}
+ \mu_i^k \Big( \vx_i^{k+1} - \mG_i \vw^k
- \alpha^k \vx_i^{k+1} 
- (1 - \alpha^k) \mG_i \vw^k
+ \mG_i \vw^{k+1} \Big) 
\nonumber
\\
& \qquad =
\vy_i^{k+1}
+ \mu_i^k \Big(
(1 - \alpha^k) \vx_i^{k+1}
- (2 - \alpha^k)\mG_i \vw^k
+ \mG_i \vw^{k+1} \Big)
\end{align*}
such that (\ref{Distr OSQP conv - kkt relations eq2}) is given as
\begin{align}
& \bigg[ \nabla f_i(\vx_i^{k+1}) 
+ \mA_i^\top \vnu_i^{k+1} 
+ \vy_i^{k+1} \nonumber \\
& \qquad \qquad \qquad
+ \mu_i^k \bigg(
(1 - \alpha^k) \vx_i^{k+1}
- (2 - \alpha^k)\mG_i \vw^k
+ \mG_i \vw^{k+1} \bigg) 
\bigg]^\top
(\vx_i^{k+1} - \vx_i^*)
= 0.
\label{Distr OSQP conv - kkt relations eq3}
\end{align}
Similarly, from (\ref{Distr OSQP conv - k+1 KKT -ztilde}), we get
\begin{align}
\Big[
- \vnu_i^{k+1} 
+ \rho_i^k (  
\vz_i^{k+1} - \vs_i^{k} + \blambda_i^k / \rho_i^k
)
\Big]^\top (\vz_i^{k+1} - \vz_i^*) = 0.
 \label{Distr OSQP conv - kkt relations eq4}
\end{align}
We rewrite the term 
$\rho_i^k ( \vz_i^{k+1} - \vs_i^{k} + \blambda_i^k / \rho_i^k)$
using (\ref{eq: DistrQP lambda update}) as follows
\begin{align*}
\rho_i^k ( \vz_i^{k+1} - \vs_i^{k} + \blambda_i^k / \rho_i^k)
& =
\rho_i^k \bigg(  
\vz_i^{k+1} - \vs_i^{k} + \blambda_i^{k+1} / \rho_i^k
- \Big( \alpha^k \vz_i^{k+1} 
+ (1 - \alpha^k) \vs_i^k
- \vs_i^{k+1} \Big) \bigg) 
\nonumber
\\
& =
\blambda_i^{k+1} 
+ \rho_i^k \big( 
\vz_i^{k+1} - \vs_i^k 
- \alpha^k \vz_i^{k+1}
- (1 - \alpha^k) \vs_i^k
+ \vs_i^{k+1}
\Big) 
\nonumber
\\
& =
\blambda_i^{k+1} 
+ \rho_i^k \Big( 
(1 - \alpha^k) \vz_i^{k+1} 
- (2 - \alpha^k) \vs_i^k
+ \vs_i^{k+1}
\Big),
\end{align*}
such that (\ref{Distr OSQP conv - kkt relations eq4}) is given as
\begin{equation}
\Bigg[
- \vnu_i^{k+1} 
+ \blambda_i^{k+1} 
+ \rho_i^k \bigg( 
(1 - \alpha^k) \vz_i^{k+1} 
- (2 - \alpha^k) \vs_i^k
+ \vs_i^{k+1} \bigg)
\Bigg]^\top (\vz_i^{k+1} - \vz_i^{*}) = 0.
\label{Distr OSQP conv - kkt relations eq5}
\end{equation}
Combining (\ref{Distr OSQP conv - kkt relations eq3}) and (\ref{Distr OSQP conv - kkt relations eq5}) and using (\ref{Distr OSQP conv - optimal point KKT local constraint}) and the fact that $\mA_i \vx_i^{k+1} - \vz_i^{k+1} = \vzero$, we obtain
\begin{equation}
\begin{aligned}
& \bigg[ \nabla f_i(\vx_i^{k+1}) 
+ \vy_i^{k+1}
+ \mu_i^k \bigg(
(1 - \alpha^k)\vx_i^{k+1}
- (2 - \alpha^k)\mG_i \vw^k
+ \mG_i \vw^{k+1} \bigg) 
\bigg]^\top 
(\vx_i^{k+1} - \vx_i^*)
\\
& \qquad 
+ \bigg[ 
\blambda_i^{k+1} 
+ \rho_i^k \bigg( 
(1 - \alpha^k) \vz_i^{k+1} 
- (2 - \alpha^k) \vs_i^k
+ \vs_i^{k+1}
\bigg)
\bigg]^\top (\vz_i^{k+1} - \vz_i^{*})
= 0
\end{aligned}
\end{equation}
which proves (\ref{R6}).
\end{proof}
\end{lemma}

\begin{lemma} \label{lemma 3}
For all $i \in \calV$, the following two relationships hold at every iteration $k$:
\begin{align}
    & \bigg( \vy_i^{k+1} - \vy^*_i
    + \mu_i^k \big(
    (1 - \alpha^k) \vx_i^{k+1}
    - (2 - \alpha^k)\mG_i \vw^k
    + \mG_i \vw^{k+1} \big) \bigg)^\top (\vx_i^{k+1} - \vx_i^{*}) \nonumber \\
    & \quad =
    \frac{1}{2 \alpha^k \mu_i^k}
    \big(
    \| \vy_i^{k+1} - \vy^*_i \|^2 
    - \| \vy_i^{k} - \vy^*_i \|^2 \big)
    + \frac{(2 - \alpha^k)}{ 2(\alpha^k)^2 \mu_i^k} \| \vy_i^{k+1} -\vy_i^k \|^2 \nonumber
    \\
    & \qquad
    + \frac{(2 - \alpha^k)\mu_i^k}{ 2(\alpha^k)^2 } \|\mG_i ( \vw^{k+1}  -  \vw^k ) \|^2 
        + \frac{\mu_i^k}{2 \alpha^k } (\|\mG_i ( \vw^{k+1}  -  \vw^* ) \|^2
    \nonumber
    \\
    & \qquad
        - \|\mG_i ( \vw^{k}  -  \vw^* ) \|^2 ) 
    + \frac{1}{\alpha^k} (\vy_i^{k+1} - \vy^*_i)^\top \mG_i \big(
        \vw^{k+1} 
        - (1 - \alpha^k) \vw^k - \alpha^k \vw^* \big) \nonumber
    \\
    & \qquad
    + \frac{1}{ (\alpha^k)^2} (\vy_i^{k+1} -\vy_i^k)^\top 
        \mG_i \big(
        (2 - \alpha^k) \vw^{k+1} 
        - (1+(1 - \alpha^k)^2) \vw^k
        - \alpha^k (1 - \alpha^k) \vw^* \big),
    \label{R7} 
    \tag{R7}
    \\
    & \bigg( \blambda_i^{k+1} - \blambda_i^*{} 
    + \rho_i^k \big( 
    (1 - \alpha^k) \vz_i^{k+1} 
    - (2 - \alpha^k) \vs_i^k
    + \vs_i^{k+1} \big) \bigg)^\top (\vz_i^{k+1} - \vz_i^{*}) \nonumber \\
    & \qquad =
    \frac{1}{2 \alpha^k \rho_i^k}
    \big(
    \| \blambda_i^{k+1} - \blambda_i^* \|^2 
    - \| \blambda_i^{k} - \blambda_i^* \|^2 \big)
    + \frac{(2 - \alpha^k)}{2 (\alpha^k)^2 \rho_i^k}
    \| \blambda_i^{k+1} - \blambda_i^k \|^2
    \nonumber \\
    & \qquad \qquad
    + \frac{\rho_i^k}{2 \alpha^k} \big(
     \|\vs_i^{k+1} - \vs_i^{*} \|^2 
    - \| \vs_i^{k} - \vs_i^{*} \|^2 \big)
    + \frac{(2 - \alpha^k)\rho_i^k}{2 (\alpha^k)^2}
    \| \vs_i^{k+1} - \vs_i^{k} \|^2
    \nonumber \\
    & \qquad \qquad
    + \frac{1}{\alpha^k}
    ( \blambda_i^{k+1} - \blambda_i^* )^\top  (- (1 - \alpha^k) \vs_i^{k}
    + \vs_i^{k+1} 
    - \alpha^k \vs_i^{*}).
    \label{R8}
    \tag{R8}
\end{align}
\end{lemma}
\begin{proof}
Let us first simplify the individual terms of the LHS of  (\ref{R7}). For that, we start by rewriting the term $\vx_i^{k+1} - \vx_i^{*}$ as follows using (\ref{R2}),
    \begin{equation*}
        \vx_i^{k+1} - \vx_i^{*} 
        = \frac{1}{\alpha^k}
     \bigg( \frac{1}{\mu_i^k}(\vy_i^{k+1} -\vy_i^k) 
    - (1 - \alpha^k) \mG_i \vw^k
    + \mG_i \vw^{k+1} 
    - \alpha^k \vx_i^{*} \bigg).
    \end{equation*}
    Using (\ref{Distr OSQP conv - optimal point KKT - g}), we can rewrite the above as
    \begin{equation}
    \vx_i^{k+1} - \vx_i^{*} 
    = \frac{1}{\alpha^k}
     \bigg( \frac{1}{\mu_i^k}(\vy_i^{k+1} -\vy_i^k) 
    - (1 - \alpha^k) \mG_i \vw^k
    + \mG_i \vw^{k+1} 
    - \alpha^k \mG_i \vw^* \bigg)
    \end{equation}
    which can be written in simplified form as 
    \begin{equation}
    \vx_i^{k+1} - \vx_i^{*} 
    = \frac{1}{\alpha^k \mu_i^k}
     (\vy_i^{k+1} -\vy_i^k) 
    + \frac{1}{\alpha^k} \mG_i \big(
     \vw^{k+1} 
    - (1 - \alpha^k) \vw^k
    - \alpha^k \vw^* \big).
    \label{R7 difference simplified eq1}
    \end{equation}
%
%
Let us now simplify the following term in the LHS of the relationship (\ref{R7})
\begin{equation}
    (1 - \alpha^k)\vx_i^{k+1}
    - (2 - \alpha^k)\mG_i \vw^k
    + \mG_i \vw^{k+1} = 
    (1 - \alpha^k)(\vx_i^{k+1} - \mG_i \vw^k)
    + \mG_i (\vw^{k+1} - \vw^k).
    \label{R7 eq4}
\end{equation}
We further simplify the term $(\vx_i^{k+1} - \mG_i \vw^k)$ using the relationship (\ref{R2}) as
\begin{equation*}
    \vx_i^{k+1} - \mG_i \vw^k = 
    \frac{1}{\alpha^k} \bigg( \frac{1}{\mu_i^k}(\vy_i^{k+1} - \vy_i^k) 
- (1 - \alpha^k) \mG_i \vw^k
+ \mG_i \vw^{k+1} \bigg) - \mG_i \vw^k,
\end{equation*}
which can be written in a simplified form as
\begin{equation}
    \vx_i^{k+1} - \mG_i \vw^k = 
    \frac{1}{\mu_i^k \alpha^k} (\vy_i^{k+1} - \vy_i^k) 
    + \frac{1}{\alpha^k} \mG_i ( \vw^{k+1}  -  \vw^k ).
    \label{R7 eq5}
\end{equation}
Substituting (\ref{R7 eq5}) in (\ref{R7 eq4}), we get
\begin{equation*}
(1 - \alpha^k)\vx_i^{k+1}
- (2 - \alpha^k)\mG_i \vw^k
+ \mG_i \vw^{k+1} 
= 
\frac{(1 - \alpha^k)}{\mu_i^k \alpha^k} (\vy_i^{k+1} - \vy_i^k) 
+ \frac{1}{\alpha^k} \mG_i ( \vw^{k+1}  -  \vw^k ).
\end{equation*}
Using the above result, we rewrite the following term on the LHS of (\ref{R7}) as
\begin{equation}
\begin{aligned}
    & \vy_i^{k+1} - \vy^*_i
    + \mu_i^k \big(
    (1 - \alpha^k) \vx_i^{k+1}
    - (2 - \alpha^k)\mG_i \vw^k
    + \mG_i \vw^{k+1} \big) \\
    & \qquad \qquad \qquad =
    \vy_i^{k+1} - \vy^*_i 
    +\frac{(1 - \alpha^k)}{ \alpha^k} (\vy_i^{k+1} - \vy_i^k) 
    + \frac{\mu_i^k}{\alpha^k} \mG_i ( \vw^{k+1}  -  \vw^k )
\end{aligned}
\label{R7 eq6}
\end{equation}
For notational simplicity, let us consider the LHS of (\ref{R7}) as LHS(R7).
Using (\ref{R7 eq6}) and (\ref{R7 difference simplified eq1}), we get
\begin{equation*}
\begin{aligned}
    \text{LHS(R7)} & = 
    \bigg(  \vy_i^{k+1} - \vy^*_i 
    +\frac{(1 - \alpha^k)}{ \alpha^k} (\vy_i^{k+1} - \vy_i^k) 
    + \frac{\mu_i^k}{\alpha^k} \mG_i ( \vw^{k+1}  -  \vw^k ) \bigg)^\top
    \\
    & \qquad \qquad
    \bigg( \frac{1}{\alpha^k \mu_i^k}
     (\vy_i^{k+1} -\vy_i^k) 
    + \frac{1}{\alpha^k} \mG_i \big(
     \vw^{k+1} 
    - (1 - \alpha^k) \vw^k
    - \alpha^k \vw^* \big) \bigg)
\end{aligned}
\end{equation*}
which can be further rewritten as
\begin{equation}
\begin{aligned}
    \text{LHS(R7)} & = 
    \frac{1}{\alpha^k \mu_i^k} (\vy_i^{k+1} - \vy^*_i)^\top (\vy_i^{k+1} -\vy_i^k)
    + \frac{1}{\alpha^k} (\vy_i^{k+1} - \vy^*_i)^\top \mG_i \big(
        \vw^{k+1} 
        - (1 - \alpha^k) \vw^k
    \\
    &~~~~
        - \alpha^k \vw^* \big) 
    + \frac{(1 - \alpha^k)}{ (\alpha^k)^2 \mu_i^k} \| \vy_i^{k+1} -\vy_i^k \|^2
    + \frac{(1 - \alpha^k)}{ (\alpha^k)^2} (\vy_i^{k+1} -\vy_i^k)^\top 
        \mG_i \big(
        \vw^{k+1} 
        \\
    &~~~~
        - (1 - \alpha^k) \vw^k
        - \alpha^k \vw^* \big)
    + \frac{1}{ (\alpha^k)^2} ( \vw^{k+1}  -  \vw^k )^\top \mG_i^\top 
        (\vy_i^{k+1} -\vy_i^k) 
    \\
    &~~~~
    + \frac{\mu_i^k}{(\alpha^k)^2} \big(\mG_i ( \vw^{k+1}  -  \vw^k ) \big)^\top \mG_i \big(
     \vw^{k+1} 
    - (1 - \alpha^k) \vw^k
    - \alpha^k \vw^* \big)
\end{aligned}
\label{R7 intermediate 1}
\end{equation}
Let us now simplify each term on the RHS of the above equation. We start with the terms including only the variables $\vy_i^{k+1}$, $\vy_i^k$ and $\vy_i^*$. 
Using the fact that $a^\top b  = \frac{1}{2} ( \|a\|^2 + \|b\|^2 - \|a-b\|^2)$, we get
\begin{equation*}
    \frac{1}{\alpha^k \mu_i^k} (\vy_i^{k+1} - \vy^*_i)^\top(\vy_i^{k+1} -\vy_i^k)
    =
    \frac{1}{2 \alpha^k \mu_i^k}
    \big(
    \| \vy_i^{k+1} - \vy^*_i \|^2 
    + \| \vy_i^{k+1} -\vy_i^k \|^2
    - \| \vy_i^{k} - \vy^*_i \|^2 \big).
\end{equation*}
Using the above result, we can write
\begin{align}
    & \frac{1}{\alpha^k \mu_i^k} (\vy_i^{k+1} - \vy^*_i)^\top (\vy_i^{k+1} -\vy_i^k) 
    + \frac{(1 - \alpha^k)}{ (\alpha^k)^2 \mu_i^k} \| \vy_i^{k+1} -\vy_i^k \|^2 
    \nonumber
    \\
    & \quad = 
    \frac{1}{2 \alpha^k \mu_i^k}
    \big(
    \| \vy_i^{k+1} - \vy^*_i \|^2 
    + \| \vy_i^{k+1} -\vy_i^k \|^2
    - \| \vy_i^{k} - \vy^*_i \|^2 \big)
    + \frac{(1 - \alpha^k)}{ (\alpha^k)^2 \mu_i^k} \| \vy_i^{k+1} -\vy_i^k \|^2 
    \nonumber \\
    & \quad =
    \frac{1}{2 \alpha^k \mu_i^k}
    \big(
    \| \vy_i^{k+1} - \vy^*_i \|^2 
    - \| \vy_i^{k} - \vy^*_i \|^2 \big)
    + \frac{(2 - \alpha^k)}{ 2(\alpha^k)^2 \mu_i^k} \| \vy_i^{k+1} -\vy_i^k \|^2
    \label{R7 y terms}
\end{align}
%
Next, we consider the following terms in the RHS of (\ref{R7 intermediate 1}) involving only the variables $\vw^{k+1}, \vw^{k}$ and $\vw^{*}$,
\begin{equation}
    \begin{aligned}
        & \frac{\mu_i^k}{(\alpha^k)^2} \big(\mG_i ( \vw^{k+1}  -  \vw^k ) \big)^\top \mG_i \big(
        \vw^{k+1} 
        - (1 - \alpha^k) \vw^k
        - \alpha^k \vw^* \big) 
        \\
        & =
        \frac{(1 - \alpha^k)\mu_i^k}{ (\alpha^k)^2 } \|\mG_i ( \vw^{k+1}  -  \vw^k ) \|^2 
        + \frac{\mu_i^k}{\alpha^k } \big(\mG_i ( \vw^{k+1}  -  \vw^k ) \big)^\top \big( \mG_i 
     (\vw^{k+1} - \vw^*) \big)
    \end{aligned}
    \label{R7 wterms 1}
\end{equation}
Using a similar approach as used to derive (\ref{R7 y terms}), we obtain
\begin{equation}
    \begin{aligned}
        &
        \frac{(1 - \alpha^k)\mu_i^k}{ (\alpha^k)^2 } \|\mG_i ( \vw^{k+1}  -  \vw^k ) \|^2 
        + \frac{\mu_i^k}{\alpha^k } \big(\mG_i ( \vw^{k+1}  -  \vw^k ) \big)^\top \big( \mG_i 
     (\vw^{k+1} - \vw^*) \big) 
     \\
     & =
     \frac{(2 - \alpha^k)\mu_i^k}{ 2(\alpha^k)^2 } \|\mG_i ( \vw^{k+1}  -  \vw^k ) \|^2 
        + \frac{\mu_i^k}{2 \alpha^k } (\|\mG_i ( \vw^{k+1}  -  \vw^* ) \|^2
        - \|\mG_i ( \vw^{k}  -  \vw^* ) \|^2 ).
    \end{aligned}
    \label{R7 w terms 2}
\end{equation}
%
%
Now, let us consider the following terms from the RHS of (\ref{R7 intermediate 1}),
\begin{align}
    & \frac{(1 - \alpha^k)}{ (\alpha^k)^2} (\vy_i^{k+1} -\vy_i^k)^\top 
        \mG_i \big(
        \vw^{k+1} 
        - (1 - \alpha^k) \vw^k
        - \alpha^k \vw^* \big) \nonumber
    \\ & \qquad \qquad \qquad \qquad \qquad \qquad
    + \frac{1}{ (\alpha^k)^2} ( \vw^{k+1}  -  \vw^k )^\top \mG_i^\top 
        (\vy_i^{k+1} -\vy_i^k) \nonumber
    \\
    &~~~ =
    \frac{1}{ (\alpha^k)^2} (\vy_i^{k+1} -\vy_i^k)^\top 
        \mG_i \big(
        (1 - \alpha^k) \vw^{k+1} 
        - (1 - \alpha^k)^2 \vw^k
        - \alpha^k (1 - \alpha^k) \vw^* 
    + \vw^{k+1}  -  \vw^k \big)
    \nonumber
    \\
    &~~~ =
    \frac{1}{ (\alpha^k)^2} (\vy_i^{k+1} -\vy_i^k)^\top 
        \mG_i \big(
        (2 - \alpha^k) \vw^{k+1} 
        - (1+(1 - \alpha^k)^2) \vw^k
        - \alpha^k (1 - \alpha^k) \vw^* \big).
    \label{R7 innerproduct terms relation}
\end{align}
%
Substituting (\ref{R7 y terms}), (\ref{R7 wterms 1}), (\ref{R7 w terms 2}), and (\ref{R7 innerproduct terms relation}) into (\ref{R7 intermediate 1}), we get
\begin{equation}
    \begin{aligned}
   \text{LHS(R7)} & = 
    \frac{1}{2 \alpha^k \mu_i^k}
    \big(
    \| \vy_i^{k+1} - \vy^*_i \|^2 
    - \| \vy_i^{k} - \vy^*_i \|^2 \big)
    + \frac{(2 - \alpha^k)}{ 2(\alpha^k)^2 \mu_i^k} \| \vy_i^{k+1} -\vy_i^k \|^2
    \\
    &~~~~
    + \frac{(2 - \alpha^k)\mu_i^k}{ 2(\alpha^k)^2 } \|\mG_i ( \vw^{k+1}  -  \vw^k ) \|^2 
        + \frac{\mu_i^k}{2 \alpha^k } (\|\mG_i ( \vw^{k+1}  -  \vw^* ) \|^2
        \\
    &~~~~
        - \|\mG_i ( \vw^{k}  -  \vw^* ) \|^2 )
    + \frac{1}{\alpha^k} (\vy_i^{k+1} - \vy^*_i)^\top \mG_i \big(
        \vw^{k+1} 
        - (1 - \alpha^k) \vw^k - \alpha^k \vw^* \big) 
    \\
    &~~~~
    + \frac{1}{ (\alpha^k)^2} (\vy_i^{k+1} -\vy_i^k)^\top 
        \mG_i \big(
        (2 - \alpha^k) \vw^{k+1} 
        - (1+(1 - \alpha^k)^2) \vw^k
        - \alpha^k (1 - \alpha^k) \vw^* \big)
    \end{aligned}
\end{equation}
which proves (\ref{R7}).

Subsequently, we prove relationship (\ref{R8}). Using similar steps as for (\ref{R7}), we get
\begin{equation}
\begin{aligned}
    & \bigg( \blambda_i^{k+1} - \blambda_i^*{} 
    + \rho_i^k \big( 
    (1 - \alpha^k) \vz_i^{k+1} 
    - (2 - \alpha^k) \vs_i^k
    + \vs_i^{k+1} \big) \bigg)^\top (\vz_i^{k+1} - \vz_i^{*}) \\
    & \qquad =
    \frac{1}{2 \alpha^k \rho_i^k}
    \big(
    \| \blambda_i^{k+1} - \blambda_i^* \|^2 
    - \| \blambda_i^{k} - \blambda_i^* \|^2 \big)
    + \frac{(2 - \alpha^k)}{2 (\alpha^k)^2 \rho_i^k}
    \| \blambda_i^{k+1} - \blambda_i^k \|^2
    \\
    & \qquad \qquad
    + \frac{\rho_i^k}{2 \alpha^k} \big(
     \|\vs_i^{k+1} - \vs_i^{*} \|^2 
    - \| \vs_i^{k} - \vs_i^{*} \|^2 \big)
    + \frac{(2 - \alpha^k)\rho_i^k}{2 (\alpha^k)^2}
    \| \vs_i^{k+1} - \vs_i^{k} \|^2
    \\
    & \qquad \qquad
    + \frac{1}{\alpha^k}
    ( \blambda_i^{k+1} - \blambda_i^* )^\top  (\vs_i^{k+1} - (1 - \alpha^k) \vs_i^{k}
    - \alpha^k \vs_i^{*}) 
    \\
    & \qquad \qquad
    + \frac{1}{ (\alpha^k)^2} (\blambda_i^{k+1} - \blambda_i^k)^\top 
        \big(
        (2 - \alpha^k) \vs_i^{k+1} 
        - (1+(1 - \alpha^k)^2) \vs_i^k
        - \alpha^k (1 - \alpha^k) \vs_i^* \big)
\end{aligned}
\label{R8 intermediate}
\end{equation}
Let us now simplify the last term of the RHS of the above equation as follows
\begin{equation}
    \begin{aligned}
        & (\blambda_i^{k+1} - \blambda_i^k)^\top 
        \big(
        (2 - \alpha^k) \vs_i^{k+1} 
        - (1+(1 - \alpha^k)^2) \vs_i^k
        - \alpha^k (1 - \alpha^k) \vs_i^* \big) \\
        & =
        (1+(1 - \alpha^k)^2) (\blambda_i^{k+1} - \blambda_i^k)^\top 
        (\vs_i^{k+1} - \vs_i^k)
        + \alpha^k (1 - \alpha^k) (\blambda_i^{k+1} - \blambda_i^k)^\top 
        (\vs_i^{k+1} - \vs_i^*).
    \end{aligned}
\end{equation}
From (\ref{eq: DistrQP s update}) and (\ref{Distr OSQP conv - optimal point KKT z in set}), we have that the vectors $\vs_i^{k}, \vs_i^{k+1}, \vs_i^* \in \calC_i$. Using (\ref{R4}), the above equation gives us
\begin{align*}
    & (\blambda_i^{k+1} - \blambda_i^k)^\top 
        \big(
        (2 - \alpha^k) \vs_i^{k+1} 
        - (1+(1 - \alpha^k)^2) \vs_i^k
        - \alpha^k (1 - \alpha^k) \vs_i^* \big) 
        \nonumber \\
    & \qquad = 
    (\blambda_i^{k+1} - \blambda_i^k)^\top 
        \big(
        (2 - \alpha^k) \vs_i^{k+1}
        - (2 + (\alpha^k)^2 - 2 \alpha^k) \vs_i^k
        + (-\alpha^k + (\alpha^k)^2) \vs_i^* \big) = 0.
\end{align*}
It follows that (\ref{R8 intermediate}) simplifies to
\begin{equation*}
\begin{aligned}
    & \bigg( \blambda_i^{k+1} - \blambda_i^*{} 
    + \rho_i^k \big( 
    (1 - \alpha^k) \vz_i^{k+1} 
    - (2 - \alpha^k) \vs_i^k
    + \vs_i^{k+1} \big) \bigg)^\top(\vz_i^{k+1} - \vz_i^{*}) \\
    & \qquad =
    \frac{1}{2 \alpha^k \rho_i^k}
    \big(
    \| \blambda_i^{k+1} - \blambda_i^* \|^2 
    - \| \blambda_i^{k} - \blambda_i^* \|^2 \big)
    + \frac{(2 - \alpha^k)}{2 (\alpha^k)^2 \rho_i^k}
    \| \blambda_i^{k+1} - \blambda_i^k \|^2
    \\
    & \qquad \qquad
    + \frac{\rho_i^k}{2 \alpha^k} \big(
     \|\vs_i^{k+1} - \vs_i^{*} \|^2 
    - \| \vs_i^{k} - \vs_i^{*} \|^2 \big)
    + \frac{(2 - \alpha^k)\rho_i^k}{2 (\alpha^k)^2}
    \| \vs_i^{k+1} - \vs_i^{k} \|^2
    \\
    & \qquad \qquad
    + \frac{1}{\alpha^k}
    ( \blambda_i^{k+1} - \blambda_i^* )^\top  (\vs_i^{k+1} - (1 - \alpha^k) \vs_i^{k}
    - \alpha^k \vs_i^{*}),
\end{aligned}
\end{equation*}
which proves (\ref{R8}).
\end{proof}
\begin{lemma} \label{Lemma 4}
The following inequality holds true at every iteration $k$:
\begin{equation}
    \begin{aligned}
    & \sum_{i \in \calV} \bigg( \frac{1}{\mu_i^k}
    \big(
    \| \vy_i^{k+1} - \vy^*_i \|^2 
    - \| \vy_i^{k} - \vy^*_i \|^2 \big) 
    + \mu_i^k (\|\mG_i ( \vw^{k+1}  -  \vw^* ) \|^2
        - \|\mG_i ( \vw^{k}  -  \vw^* ) \|^2 ) 
    \\
     & \qquad \qquad
     + \frac{1}{\rho_i^k}
    \big(
    \| \blambda_i^{k+1} - \blambda_i^* \|^2 
    - \| \blambda_i^{k} - \blambda_i^* \|^2 \big)
    + \rho_i^k \big(
     \|\vs_i^{k+1} - \vs_i^{*} \|^2 
    - \| \vs_i^{k} - \vs_i^{*} \|^2 \big) \bigg) \\
    & \leq
    - \frac{(2 - \alpha^k)}{\alpha^k} 
    \sum_{i \in \calV} \bigg( \frac{1}{\mu_i^k} \| \vy_i^{k+1} -\vy_i^k \|^2
    + \mu_i^k \|\mG_i (\vw^{k+1}-\vw^k) \|^2 
    + \frac{1}{\rho_i^k}
            \| \blambda_i^{k+1} - \blambda_i^k \|^2
    \\
    & \qquad \qquad \qquad \qquad \qquad
    + \rho_i^k
    \| \vs_i^{k+1} - \vs_i^{k} \|^2 \bigg).
    \end{aligned}
    \label{Distr OSQP - descent relation final}
    \end{equation}
\end{lemma}
\begin{proof}
    We start by combining the relationships (\ref{R5}) and (\ref{R6}) to get
    \begin{equation}
        \begin{aligned}
            & \bigg( \vy_i^{k+1} - \vy^*_i
            + \mu_i^k \big(
            (1 - \alpha^k) \vx_i^{k+1}
            - (2 - \alpha^k)\mG_i \vw^k
            + \mG_i \vw^{k+1} \big) \bigg)^\top (\vx_i^{k+1} - \vx_i^{*}) \\
            & \qquad
            + \bigg( \blambda_i^{k+1} - \blambda_i^*{} 
            + \rho_i^k \big( 
            (1 - \alpha^k) \vz_i^{k+1} 
            - (2 - \alpha^k) \vs_i^k
            + \vs_i^{k+1} \big) \bigg)^\top (\vz_i^{k+1} - \vz_i^{*}) \\
            & 
            = - ( \nabla f_i(\vx_i^{k+1}) - \nabla f_i(\vx_i^*) )^\top (\vx_i^{k+1} - \vx_i^{*}).
        \end{aligned}
        \label{descent eq1}
    \end{equation}
    Since $f_i$ is convex, then we have
        $( \nabla f_i(\vx_i^{k+1}) - \nabla f_i(\vx_i^*) )^\top (\vx_i^{k+1} - \vx_i^{*}) \geq 0$,
    which gives
    \begin{equation}
        \begin{aligned}
            & \bigg( \vy_i^{k+1} - \vy^*_i
            + \mu_i^k \big(
            (1 - \alpha^k) \vx_i^{k+1}
            - (2 - \alpha^k)\mG_i \vw^k
            + \mG_i \vw^{k+1} \big) \bigg)^\top (\vx_i^{k+1} - \vx_i^{*}) \\
            & \qquad
            + \bigg( \blambda_i^{k+1} - \blambda_i^*{} 
            + \rho_i^k \big( 
            (1 - \alpha^k) \vz_i^{k+1} 
            - (2 - \alpha^k) \vs_i^k
            + \vs_i^{k+1} \big) \bigg)^\top (\vz_i^{k+1} - \vz_i^{*})
            \leq 0.
        \end{aligned}
        \label{lemma 4 proof: agent wise eq}
    \end{equation}
    Summing (\ref{lemma 4 proof: agent wise eq}) over all $i \in \calV$, we get
    \begin{equation}
        \begin{aligned}
            & \sum_{i \in \calV} \bigg( \vy_i^{k+1} - \vy^*_i
            + \mu_i^k \big(
            (1 - \alpha^k) \vx_i^{k+1}
            - (2 - \alpha^k)\mG_i \vw^k
            + \mG_i \vw^{k+1} \big) \bigg)^\top (\vx_i^{k+1} - \vx_i^{*}) \\
            & \qquad
            +  \sum_{i \in \calV}  \bigg( \blambda_i^{k+1} - \blambda_i^*{} 
            + \rho_i^k \big( 
            (1 - \alpha^k) \vz_i^{k+1} 
            - (2 - \alpha^k) \vs_i^k
            + \vs_i^{k+1} \big) \bigg)^\top (\vz_i^{k+1} - \vz_i^{*})
            \leq 0.
        \end{aligned}
    \end{equation}
    Now, we use the relationships (\ref{R7}) and (\ref{R8}) to rewrite the above inequality as
    \begin{equation}
    \begin{aligned}
    0 \geq & \sum_{i \in \calV} \bigg( \frac{1}{2 \alpha^k \mu_i^k}
    \big(
    \| \vy_i^{k+1} - \vy^*_i \|^2 
    - \| \vy_i^{k} - \vy^*_i \|^2 \big)
    + \frac{(2 - \alpha^k)}{ 2(\alpha^k)^2 \mu_i^k} \| \vy_i^{k+1} -\vy_i^k \|^2 
    \\
    & \qquad
    + \frac{(2 - \alpha^k)\mu_i^k}{ 2(\alpha^k)^2 } \|\mG_i ( \vw^{k+1}  -  \vw^k ) \|^2 
        + \frac{\mu_i^k}{2 \alpha^k } (\|\mG_i ( \vw^{k+1}  -  \vw^* ) \|^2
    \\
    & \qquad
        - \|\mG_i ( \vw^{k}  -  \vw^* ) \|^2 ) 
    + \frac{1}{\alpha^k} (\vy_i^{k+1} - \vy^*_i)^\top \mG_i \big(
        \vw^{k+1} 
        - (1 - \alpha^k) \vw^k - \alpha^k \vw^* \big) 
    \\
    & \qquad
    + \frac{1}{ (\alpha^k)^2} (\vy_i^{k+1} -\vy_i^k)^\top 
        \mG_i \big(
        (2 - \alpha^k) \vw^{k+1} 
        - (1+(1 - \alpha^k)^2) \vw^k
        - \alpha^k (1 - \alpha^k) \vw^* \big)\\
     & \qquad
     + \frac{1}{2 \alpha^k \rho_i^k}
    \big(
    \| \blambda_i^{k+1} - \blambda_i^* \|^2 
    - \| \blambda_i^{k} - \blambda_i^* \|^2 \big)
    + \frac{(2 - \alpha^k)}{2 (\alpha^k)^2 \rho_i^k}
    \| \blambda_i^{k+1} - \blambda_i^k \|^2
    \\
    & \qquad
    + \frac{\rho_i^k}{2 \alpha^k} \big(
     \|\vs_i^{k+1} - \vs_i^{*} \|^2 
    - \| \vs_i^{k} - \vs_i^{*} \|^2 \big)
    + \frac{(2 - \alpha^k)\rho_i^k}{2 (\alpha^k)^2}
    \| \vs_i^{k+1} - \vs_i^{k} \|^2
    \\
    & \qquad 
    + \frac{1}{\alpha^k}
    ( \blambda_i^{k+1} - \blambda_i^* )^\top  (- (1 - \alpha^k) \vs_i^{k}
    + \vs_i^{k+1} 
    - \alpha^k \vs_i^{*}) \bigg).
    \end{aligned}
    \label{descent eq2}
    \end{equation}
    Let us now further simplify the terms on the RHS of the above equation. 
    For that, let us start with the last term on the RHS.
    We have
    \begin{equation}
        \begin{aligned}
            ( \blambda_i^{k+1} - \blambda_i^* )^\top  (- (1 - \alpha^k) \vs_i^{k}
        + \vs_i^{k+1} 
        - \alpha^k \vs_i^{*}) 
        & =
        ( \blambda_i^{k+1} - \blambda_i^* )^\top(\vs_i^{k+1} - \vs_i^{*}) 
        \\
        & \qquad
        - (1 - \alpha^k) ( \blambda_i^{k+1} - \blambda_i^* )^\top(\vs_i^{k} - \vs_i^{*})
        \end{aligned}
        \label{descent mixed lambda}
    \end{equation}
    Using (\ref{R4}), (\ref{Distr OSQP conv - optimal point KKT -z}), and the fact that $\vs_i^{k}, \vs_i^{k+1}, \vs_i^* \in \calC_i$, we get
    \begin{align}
        & ( \blambda_i^{k+1} - \blambda_i^* )^\top(\vs_i^{k+1} - \vs_i^{*}) \geq 0, \label{descent lambda 1}\\
        & ( \blambda_i^{k+1} - \blambda_i^* )^\top(\vs_i^{k} - \vs_i^{*}) \geq  0. \label{descent lambda 2}
    \end{align}
   Thus, for $\alpha^k \geq 1$, combining (\ref{descent mixed lambda}), (\ref{descent lambda 1}), and (\ref{descent lambda 2}), we get
   \begin{equation}
       ( \blambda_i^{k+1} - \blambda_i^* )^\top  (- (1 - \alpha^k) \vs_i^{k}
        + \vs_i^{k+1} 
        - \alpha^k \vs_i^{*}) \geq 0.
        \label{descent eq 8}
   \end{equation}
   Now, the following results hold based on the relationship (\ref{R1}) and (\ref{Distr OSQP conv - optimal point KKT - g}).
    \begin{align}
        & \sum_{i \in \calV} (\vy_i^{k+1} - \vy^*_i)^\top \mG_i = 0,
       \quad
        \sum_{i \in \calV} (\vy_i^{k+1} -\vy_i^k)^\top 
        \mG_i = 0.
        \label{descent eq 7}
    \end{align}
    %
    By substituting (\ref{descent eq 8}) and (\ref{descent eq 7}) in (\ref{descent eq2}), and by rearranging the terms, we get
    \begin{equation}
    \begin{aligned}
    & \sum_{i \in \calV} \bigg( \frac{1}{2 \alpha^k \mu_i^k}
    \big(
    \| \vy_i^{k+1} - \vy^*_i \|^2 
    - \| \vy_i^{k} - \vy^*_i \|^2 \big) 
    + \frac{\mu_i^k}{2 \alpha^k } (\|\mG_i ( \vw^{k+1}  -  \vw^* ) \|^2
        - \|\mG_i ( \vw^{k}  -  \vw^* ) \|^2 ) 
    \\
     & \qquad \qquad
     + \frac{1}{2 \alpha^k \rho_i^k}
    \big(
    \| \blambda_i^{k+1} - \blambda_i^* \|^2 
    - \| \blambda_i^{k} - \blambda_i^* \|^2 \big)
    + \frac{\rho_i^k}{2 \alpha^k} \big(
     \|\vs_i^{k+1} - \vs_i^{*} \|^2 
    - \| \vs_i^{k} - \vs_i^{*} \|^2 \big) \bigg) \\
    & \leq
    - \sum_{i \in \calV} \bigg( \frac{(2 - \alpha^k)}{                  2(\alpha^k)^2 \mu_i^k} \| \vy_i^{k+1} -\vy_i^k \|^2
    + \frac{(2 - \alpha^k)\mu_i^k}{ 2(\alpha^k)^2 } \|\mG_i (           \vw^{k+1}  -  \vw^k ) \|^2 
    \\
    & \qquad \qquad \qquad
    + \frac{(2 - \alpha^k)}{2 (\alpha^k)^2 \rho_i^k}
            \| \blambda_i^{k+1} - \blambda_i^k \|^2
    + \frac{(2 - \alpha^k)\rho_i^k}{2 (\alpha^k)^2}
    \| \vs_i^{k+1} - \vs_i^{k} \|^2 \bigg).
    \end{aligned}
    \end{equation}
    Since, $\alpha^k \geq 1$, we can multiply the above equation with $2 \alpha^k$ to obtain
    \begin{equation}
    \begin{aligned}
    & \sum_{i \in \calV} \bigg( \frac{1}{\mu_i^k}
    \big(
    \| \vy_i^{k+1} - \vy^*_i \|^2 
    - \| \vy_i^{k} - \vy^*_i \|^2 \big) 
    + \mu_i^k (\|\mG_i ( \vw^{k+1}  -  \vw^* ) \|^2
        - \|\mG_i ( \vw^{k}  -  \vw^* ) \|^2 ) 
    \\
     & \qquad \qquad
     + \frac{1}{\rho_i^k}
    \big(
    \| \blambda_i^{k+1} - \blambda_i^* \|^2 
    - \| \blambda_i^{k} - \blambda_i^* \|^2 \big)
    + \rho_i^k \big(
     \|\vs_i^{k+1} - \vs_i^{*} \|^2 
    - \| \vs_i^{k} - \vs_i^{*} \|^2 \big) \bigg) \\
    & \leq
    - \frac{(2 - \alpha^k)}{\alpha^k} 
    \sum_{i \in \calV} \bigg( \frac{1}{\mu_i^k} \| \vy_i^{k+1} -\vy_i^k \|^2
    + \mu_i^k \|\mG_i (\vw^{k+1}-\vw^k) \|^2 
    + \frac{1}{\rho_i^k}
            \| \blambda_i^{k+1} - \blambda_i^k \|^2
    \\
    & \qquad \qquad \qquad \qquad \qquad
    + \rho_i^k
    \| \vs_i^{k+1} - \vs_i^{k} \|^2 \bigg).
    \end{aligned}
    \end{equation}
\end{proof}
\subsection{Proof of Theorem \ref{theorem: DistrQP convergence}} \label{Proof of DistrQP convergence}
Let us first rewrite the relation (\ref{Distr OSQP - descent relation final}) derived in Lemma \ref{Lemma 4} for $\alpha^k \in [1,2)$, as
\begin{equation*}
    \begin{aligned}
    & \sum_{i \in \calV} \bigg( \frac{1}{\mu_i^k}
    \big(
    \| \vy_i^{k+1} - \vy^*_i \|^2 
    - \| \vy_i^{k} - \vy^*_i \|^2 \big) 
    + \mu_i^k (\|\mG_i ( \vw^{k+1}  -  \vw^* ) \|^2
        - \|\mG_i ( \vw^{k}  -  \vw^* ) \|^2 ) 
    \\
     & \qquad \qquad
     + \frac{1}{\rho_i^k}
    \big(
    \| \blambda_i^{k+1} - \blambda_i^* \|^2 
    - \| \blambda_i^{k} - \blambda_i^* \|^2 \big)
    + \rho_i^k \big(
     \|\vs_i^{k+1} - \vs_i^{*} \|^2 
    - \| \vs_i^{k} - \vs_i^{*} \|^2 \big) \bigg) \\
    & \leq
    - \frac{(2 - \alpha^k)}{\alpha^k} 
    \sum_{i \in \calV} \bigg( \frac{1}{\mu_i^k} \| \vy_i^{k+1} -\vy_i^k \|^2
    + \mu_i^k \|\mG_i (\vw^{k+1}-\vw^k) \|^2 
    + \frac{1}{\rho_i^k}
            \| \blambda_i^{k+1} - \blambda_i^k \|^2
    \\
    & \qquad \qquad \qquad \qquad \qquad
    + \rho_i^k
    \| \vs_i^{k+1} - \vs_i^{k} \|^2 \bigg)
    \end{aligned}
\end{equation*}
which can be rearranged as follows
\begin{equation*}
    \begin{aligned}
    &  
    \frac{(2 - \alpha^k)}{\alpha^k} 
    \sum_{i \in \calV} \bigg( \frac{1}{\mu_i^k} \| \vy_i^{k+1} -\vy_i^k \|^2
    + \mu_i^k \|\mG_i (\vw^{k+1}-\vw^k) \|^2 
    + \frac{1}{\rho_i^k}
            \| \blambda_i^{k+1} - \blambda_i^k \|^2
    + \rho_i^k
    \| \vs_i^{k+1} - \vs_i^{k} \|^2 \bigg)
    \\
    & \leq
    \sum_{i \in \calV} \bigg( \frac{1}{\mu_i^k}
    \big(
    \| \vy_i^{k} - \vy^*_i \|^2 
    - \| \vy_i^{k+1} - \vy^*_i \|^2 \big) 
    + \mu_i^k (\|\mG_i ( \vw^{k}  -  \vw^* ) \|^2
        - \|\mG_i ( \vw^{k+1}  -  \vw^* ) \|^2 ) 
    \\
     & \qquad \qquad
     + \frac{1}{\rho_i^k}
    \big(
    \| \blambda_i^{k} - \blambda_i^* \|^2 
    - \| \blambda_i^{k+1} - \blambda_i^* \|^2 \big)
    + \rho_i^k \big(
     \|\vs_i^{k} - \vs_i^{*} \|^2 
    - \| \vs_i^{k+1} - \vs_i^{*} \|^2 \big) \bigg).
    %
    \end{aligned}
    \label{main conv - descent relation}
\end{equation*}
For convenience, let us define for each iteration $k$, the terms $\eta_i^k$, $i \in \calV$, and $\eta^k$ such that 
\begin{equation*}
\eta_i^k + 1 = \max \left( 
\frac{\rho_i^k}{\rho_i^{k-1}}, 
\frac{\rho_i^{k-1}}{\rho_i^k},
\frac{\mu_i^k}{\mu_i^{k-1}}, 
\frac{\mu_i^{k-1}}{\mu_i^k} 
\right), 
\quad
\eta^k_{\text{max}} = \max_{i \in \calV} \eta_i^k,
\label{eta definition}
\end{equation*}
and the term $V^{k}$ as 
\begin{equation*}
    \begin{aligned}
    V^{k} = \sum_{i \in \calV} \bigg( \frac{1}{\mu_i^{k-1}} \| \vy_i^{k} - \vy^*_i \|^2 
    + \mu_i^{k-1} \|\mG_i ( \vw^{k}  -  \vw^* ) \|^2
    + \frac{1}{\rho_i^{k-1}} \| \blambda_i^{k} - \blambda_i^* \|^2 
    + \rho_i^{k-1} \|\vs_i^{k} - \vs_i^{*} \|^2 \bigg).
    \end{aligned}
\end{equation*}
Based on the definition of $\eta_i^k$, we can write
\begin{equation*}
    \begin{aligned}
    & \frac{1}{\mu_i^k} \| \vy_i^{k} - \vy^*_i \|^2 
    + \mu_i^k \|\mG_i ( \vw^{k}  -  \vw^* ) \|^2
    + \frac{1}{\rho_i^k} \| \blambda_i^{k} - \blambda_i^* \|^2 
    + \rho_i^k \|\vs_i^{k} - \vs_i^{*} \|^2 \\
    &
    \quad \leq 
    (\eta^k_i + 1) \bigg( 
    \frac{1}{\mu_i^{k-1}} \| \vy_i^{k} - \vy^*_i \|^2 
    + \mu_i^{k-1} \|\mG_i ( \vw^{k}  -  \vw^* ) \|^2
    + \frac{1}{\rho_i^{k-1}} \| \blambda_i^{k} - \blambda_i^* \|^2  
    + \rho_i^{k-1} \|\vs_i^{k} - \vs_i^{*} \|^2 \bigg).
    \end{aligned}
\end{equation*}
By adding the above result over all $i \in \calV$, and using the fact that $\eta^k_{\text{max}} \geq \eta^k_i$ for all $i$, we get
\begin{equation}
\begin{aligned}
    & \sum_{i \in \calV} \bigg( \frac{1}{\mu_i^k} \| \vy_i^{k} - \vy^*_i \|^2 
    + \mu_i^k \|\mG_i ( \vw^{k}  -  \vw^* ) \|^2
    + \frac{1}{\rho_i^k} \| \blambda_i^{k} - \blambda_i^* \|^2 
    + \rho_i^k \|\vs_i^{k} - \vs_i^{*} \|^2 \bigg) \\
    &
    \leq 
    \sum_{i \in \calV} (\eta^k_i + 1) \bigg( 
    \frac{1}{\mu_i^{k-1}} \| \vy_i^{k} - \vy^*_i \|^2 
    + \mu_i^{k-1} \|\mG_i ( \vw^{k}  -  \vw^* ) \|^2
    + \frac{1}{\rho_i^{k-1}} \| \blambda_i^{k} - \blambda_i^* \|^2 
    \\
    & \qquad \qquad 
    + \rho_i^{k-1} \|\vs_i^{k} - \vs_i^{*} \|^2 \bigg) 
    \\
    & \quad \leq 
    (\eta^k_{\text{max}} + 1) \sum_{i \in \calV} \bigg( 
    \frac{1}{\mu_i^{k-1}} \| \vy_i^{k} - \vy^*_i \|^2 
    + \mu_i^{k-1} \|\mG_i ( \vw^{k}  -  \vw^* ) \|^2
    + \frac{1}{\rho_i^{k-1}} \| \blambda_i^{k} - \blambda_i^* \|^2 
    \\
    & \qquad \qquad 
    + \rho_i^{k-1} \|\vs_i^{k} - \vs_i^{*} \|^2 \bigg) \\
    & \quad \quad 
    = (\eta^k_{\text{max}} + 1) V^k.
\end{aligned}
\end{equation}
Substituting the above result in (\ref{main conv - descent relation}), we get
\begin{equation}
\begin{aligned} 
    \frac{(2 - \alpha^k)}{\alpha^k} 
    \sum_{i \in \calV} \bigg( \frac{1}{\mu_i^k} \| \vy_i^{k+1} -\vy_i^k \|^2
    & + \mu_i^k \|\mG_i (\vw^{k+1}-\vw^k) \|^2 
    + \frac{1}{\rho_i^k}
            \| \blambda_i^{k+1} - \blambda_i^k \|^2
    \\
    & \qquad
    + \rho_i^k
    \| \vs_i^{k+1} - \vs_i^{k} \|^2 \bigg)
    \leq
    (\eta^k_{\text{max}} + 1) V^k - V^{k+1}.
\end{aligned}
\label{main conv - descent 2}
\end{equation}
Now that we have derived the above relation, we need to next prove that $V^k$ is bounded. By the definition of $V^k$, we have that $V^k$ is lower bounded by zero. Thus, we now prove that $V^k$ is upper bounded. From (\ref{main conv - descent 2}), we have
\begin{equation}
    V^{k+1} \leq
    (\eta^k_{\text{max}} + 1) V^k,
\end{equation}
which leads to the following relationship 
\begin{equation}
    V^{k+1} \leq \prod_{l = 1}^{k} (\eta^l_{\text{max}} + 1) V^1.
    \label{main conv - v upper bound}
\end{equation}
It should be noted that based on Assumption \ref{assumption: eta}, we have $(\eta^k_{\text{max}}+1) \rightarrow 1$, as $k \rightarrow \infty$.
Therefore, (\ref{main conv - v upper bound}) implies that $V^{k+1}$ is upper bounded for all $k$, and there exists $V_{\text{max}}$ such that  
\begin{equation}
    V^k \leq V_{\text{max}} < \infty, \quad \text{for all } k.
    \label{main conv - v bound}
\end{equation}
Let us now consider adding the result (\ref{main conv - descent 2}) over $k$ as follows
\begin{equation}
    \begin{aligned}
    \sum_{k=1}^{\infty} \frac{(2 - \alpha^k)}{\alpha^k} 
    \sum_{i \in \calV} \bigg( \frac{1}{\mu_i^k} \| \vy_i^{k+1} -\vy_i^k \|^2
    & + \mu_i^k \|\mG_i (\vw^{k+1}-\vw^k) \|^2 
    + \frac{1}{\rho_i^k}
            \| \blambda_i^{k+1} - \blambda_i^k \|^2
    \\
    & \qquad
    + \rho_i^k
    \| \vs_i^{k+1} - \vs_i^{k} \|^2 \bigg)
    \leq
    \sum_{k =1}^{\infty} (\eta^k_{\text{max}} + 1) V^k - V^{k+1}.
    \end{aligned}
    \label{main conv - descent 3}
\end{equation}
The term on the RHS of the above equation can be further simplified as follows
\begin{equation*}
    \sum_{k =1}^{\infty} (\eta^k_{\text{max}} + 1) V^k - V^{k+1} 
    = \sum_{k =1}^{\infty} \eta^k_{\text{max}} V^k 
        + \sum_{k =1}^{\infty} (V^k- V^{k+1})
    = V^1 - V^{\infty} + \sum_{k =1}^{\infty} \eta^k_{\text{max}} V^k.
\end{equation*}
Based on Assumption \ref{assumption: eta}, we have $\eta^k_{\text{max}} \rightarrow 0$ as $k \rightarrow \infty$, which implies that
\begin{equation}
     \sum_{k =1}^{\infty} \eta^k_{\text{max}} < \infty.
\end{equation}
Using the above fact and (\ref{main conv - v bound}), we can upper bound $\sum_{k =1}^{\infty} \eta^k_{\text{max}} V^k$ as follows
\begin{equation}
     \sum_{k =1}^{\infty} \eta^k_{\text{max}} V^k \leq 
      \bigg( \sum_{k =1}^{\infty} \eta^k_{\text{max}} \bigg) V_{\text{max}} 
      < \infty.
\end{equation}
Using the facts that $V^1$ is upper bounded, and $V^{\infty}$ is lower bounded by zero, and using the above equation, we get
\begin{equation*}
    V^1 - V^{\infty} + \sum_{k =1}^{\infty} \eta^k_{\text{max}} V^k
    \leq
    V^1 + \sum_{k =1}^{\infty} \eta^k_{\text{max}} V^k
    < \infty.
\end{equation*}
Thus, we can rewrite (\ref{main conv - descent 3}) as
\begin{equation}
    \begin{aligned}
    \sum_{k=1}^{\infty} \frac{(2 - \alpha^k)}{\alpha^k} 
    \sum_{i \in \calV} \bigg( \frac{1}{\mu_i^k} \| \vy_i^{k+1} -\vy_i^k \|^2
    & + \mu_i^k \|\mG_i (\vw^{k+1}-\vw^k) \|^2 
    + \frac{1}{\rho_i^k}
            \| \blambda_i^{k+1} - \blambda_i^k \|^2
    \\
    & \qquad
    + \rho_i^k
    \| \vs_i^{k+1} - \vs_i^{k} \|^2 \bigg)
    < \infty.
    \end{aligned}
    \label{main conv - descent 4}
\end{equation}
Since $\alpha^k \in [1,2)$, we have $\frac{(2 - \alpha^k)}{\alpha^k} > 0$ for all $k$. Further, we have $0< \mu_i^k, \rho_i^k < \infty$ for all $k$. Thus, (\ref{main conv - descent 4}) implies that as $k \rightarrow \infty$, 
\begin{equation}
(\vy_i^{k+1} -\vy_i^k) \rightarrow \vzero, \quad
\mG_i (\vw^{k+1}-\vw^k) \rightarrow \vzero, \quad 
(\blambda_i^{k+1} - \blambda_i^k) \rightarrow \vzero, \quad
\vs_i^{k+1} - \vs_i^{k} \rightarrow \vzero,
\label{convergence result 1}
\end{equation}
for all $i \in \calV$.
This proves the convergence of the variables $\vy_i, \blambda_i$ and $\vs_i$. Further, it follows that $\mG (\vw^{k+1}-\vw^k) \rightarrow \vzero$. Since $\mG$ is full column rank, this implies that as $k \rightarrow \infty$, 
\begin{equation}
    (\vw^{k+1}-\vw^k) \rightarrow \vzero,
\end{equation}
which proves the convergence of the global variable $\vw$. Subsequently, combining (\ref{R2}), (\ref{R3}), and the convergence result (\ref{convergence result 1}), we also obtain that as $k \rightarrow \infty$,
\begin{equation}
(\vx_i^{k+1} -\vx_i^k) \rightarrow \vzero, \quad
(\vz_i^{k+1}-\vz_i^k) \rightarrow \vzero,
\end{equation}
for all $i \in \calV$. Hence, we have proved the convergence of the DistributedQP algorithm. 

Now that we have proved the convergence of all variables, we proceed with verifying that the limit point of convergence is the optimal solution to problem (\ref{eq: distr qp problem v3}). For that, we need to check if the limit point satisfies the KKT condition (\ref{KKT conditions optimal}) for the problem \ref{eq: distr qp problem v3}. The convergence of the dual variables $\vy_i$ and $\blambda_i$, and the update steps verify that the limit points have constraint feasibility (\ref{Distr OSQP conv - optimal point KKT z constraint} - \ref{Distr OSQP conv - optimal point KKT z in set}). The constraint feasibility of the limit points and the optimality conditions of $(k+1)$-th update of $\vx_i, \vz_i$ (\ref{KKT condition for k+1 subproblem}) imply that the limit points satisfy the optimality conditions (\ref{Distr OSQP conv - optimal point KKT -xtilde} - \ref{Distr OSQP conv - optimal point KKT -ztilde}). Further, using relations (\ref{R1}) and (\ref{R4}), we can prove that the limit points also satisfy (\ref{Distr OSQP conv - optimal point KKT -z} - \ref{Distr OSQP conv - optimal point KKT - g}). 

\section{Details on DeepDistributedQP Feedback Policies}
\label{sec: appendix - feedback residuals}

In DeepDistributedQP, the penalty parameters are given by 
\begin{equation}
\rho_i^k = \mathrm{SoftPlus}\Big( \bar{\rho}_i^k 
+ \underbrace{\pi_{i,\rho}^k (r_{i, \rho}^k, s_{i, \rho}^k; \theta_{i, \rho}^k)}_{\hat{\rho}_i^k} \Big), \quad 
\mu_i^k = \mathrm{SoftPlus}\Big( \bar{\mu}_i^k + \underbrace{\pi_{i,\mu}^k (r_{i, \mu}^k, s_{i, \mu}^k; \theta_{i, \mu}^k)}_{\hat{\mu}_i^k}
\Big)
\end{equation}
where $\bar{\rho}_i^k$, $\bar{\mu}_i^k$ are learnable feed-forward parameters and $\hat{\rho}_i^k$, $\hat{\mu}_i^k$ and the feedback parts. The latter are obtain through the learnable policies $\pi_{i, \cdot}^k (r_{i, \cdot}^k, s_{i, \cdot}^k; \theta_{i, \cdot}^k)$ 
parameterized by fully-connected neural network layers with inputs $r_{i, \cdot}^k, s_{i, \cdot}^k$ and weights $\theta_{i,\cdot}^k$. The analytical expressions for $r_{i, \cdot}^k, s_{i, \cdot}^k$ are provided as follows:
%
%
\begin{subequations}
\begin{alignat}{2}
& r_{i, \rho}^k = 
\begin{bmatrix}
\| \vz_i^k - \vs_i^k \|_2 
\\[0.1cm]
\| \mA_i \vx_i^k - \vs_i^k \|_2
\end{bmatrix}, 
\quad
&& s_{i, \rho}^k =
\begin{bmatrix}
\| \vs_i^k - \vs_i^{k-1} \|_2
\\[0.1cm]
\| \mQ_i \vx_i^k + \vq_i + \mA_i^\top \blambda_i^k \|_2 
\end{bmatrix} 
\\
& r_{i, \mu}^k = \| \vx_i^k - \tilde{\vw}_i^k \|_2, 
\quad
&& s_{i, \mu}^k = \| \tilde{\vw}_i^k - \tilde{\vw}_i^{k-1} \|_2,
\end{alignat}
\end{subequations}
being motivated by the primal and dual residuals of ADMM \citep[Section 3]{boyd2011distributed} and the ones used in the OSQP algorithm \citep{stellato2020osqp}.

\section{The Centralized Version: DeepQP}
\label{sec: appendix - non distr deep qp}
The centralized version of DeepDistributedQP boils down to simply unfolding the iterates of the standard OSQP algorithm for solving centralized QPs (\ref{eq: general qp problem}), while applying the same principles as in Section \ref{sec: distr deep qp - main arch} for DeepDistributedQP. 

For convenience, we repeat the OSQP updates from \cite{stellato2020osqp} here:

\begin{enumerate}
\item \textit{Update for $(\vx, \vz)$:} Solve linear system
\begin{equation}
\begin{bmatrix}
\mQ + \sigma \mI & \mA^\top \\
\mA & - 1/\rho^k \mI 
\end{bmatrix}
\begin{bmatrix}
\vx^{k+1} \\
\vnu^{k+1}
\end{bmatrix}
= 
\begin{bmatrix}
\sigma \vt^k - \vq
\\
\vs^k - 1/\rho^k \blambda^k
\end{bmatrix}
\label{eq: osqp - first update - direct}
\end{equation}

and update 
\begin{equation}
\vz^{k+1} =  \vs^k + 1/\rho^k (\vnu^{k+1} - \blambda^k).
\end{equation}
As explained in \cite{stellato2020osqp}, as the scale of the system )\ref{eq: osqp - first update - direct}) increases, it is often preferable to solve the following system instead,
\begin{equation}
(\mQ + \sigma \mI + \rho^k \mA^\top \mA) \vx^{k+1} = 
\sigma \vx^k - \vq + \mA^\top (\rho^k \vz^k - \vy^k),
\label{eq: osqp - first update - indirect}
\end{equation}
using a method such as CG.
\item \textit{Update for $(\vt, \vs)$:}
\begin{subequations}    
\begin{align}
\vt^{k+1} & = \alpha^k \vx^{k+1} + (1 - \alpha^k) \vt^k
\\
\vs^{k+1} & = \Pi_{\calC} \left( \alpha^k \vz^{k+1} + (1 - \alpha^k) \vs^k + \blambda^k / \rho^k \right)
\end{align}
\end{subequations}
\item \textit{Dual update for $\blambda$:}
\begin{equation}
\blambda^{k+1} = \blambda^k + \rho^k (\alpha^k \vz^{k+1} + (1-\alpha^k) \vs^k - \vs^{k+1})
\end{equation}
\end{enumerate}
The DeepQP framework then emerges through unfolding the OSQP updates following the same methodology as in DeepDistributedQP. In particular, its iterations are unrolled for a prescribed amount of $K$ iterations as shown in Fig. \ref{fig: DeepQP}.

\section{Proof of Indirect Method Implicit Differentiation}
\label{sec: appendix - IFT proof}

We start by restating the implicit function theorem, whose proof can be found in \citep{krantz2002implicit}.
\begin{lemma}[Implicit Function Theorem] \label{lemma:IFT}
    Let $r: \sR^n \times \sR^m \to \sR^n$ be a continuously differentiable function. Let $(\vx_0, \btheta_0)$ be a point such that $r(\vx_0, \btheta_0) = 0$. If the Jacobian matrix $\frac{\partial r}{\partial \vx}(\vx_0, \btheta_0)$ is invertible, then there exists a function $\vx^*(\cdot)$ defined in a neighborhood of $\btheta_0$ such that $\vx^*(\btheta_0) = \vx_0$, and
    \begin{align}
        \frac{\partial \vx^*}{\partial \btheta} (\btheta) = - \left( \frac{\partial r}{\partial \vx}(\vx^*(\btheta), \btheta) \right)^{-1} \frac{\partial r}{\partial \btheta}(\vx^*(\btheta), \btheta) .
    \end{align}
\end{lemma}

\begin{proof}[Proof of Theorem \ref{thm:implicit_differentiation}]
Let $\btheta = (\bar{\mQ}_i^k, \bar{\vb}_i^k)$ be the concatenation of all the parameters in (\ref{eq:indirect_method_linear_system}). $\bar{\mQ}_i^k$ is always positive definite since $\mQ_i$ is positive definite and the penalty parameters are always non-negative. Therefore, (\ref{eq:indirect_method_linear_system}) has a unique solution $\vx_i^{k + 1}$ satisfying $r(\vx_i^{k + 1}, \btheta) \vcentcolon= \bar{\mQ}_i^k \vx_i^{k + 1} - \bar{\vb}_i^k = 0$. Applying (\ref{lemma:IFT}) to this residual function yields the relationship $\frac{\partial \vx_i^{k + 1}}{\partial \btheta}(\btheta) = -(\bar{\mQ}_i^k)^{-1} \frac{\partial r}{\partial \btheta}(\vx_i^{k + 1}(\btheta), \btheta)$.

Now, for any downstream loss function $L(\vx_i^{k + 1}(\btheta))$, we have that
\begin{align}
    \nabla_{\btheta} L(\vx_i^{k + 1}(\btheta)) &= \frac{\partial \vx_i^{k + 1}}{\partial \btheta}(\btheta) \nabla_{\vx} L(\vx_i^{k + 1}(\btheta)) \nonumber \\
    &= -\frac{\partial r}{\partial \btheta}(\vx_i^{k + 1}(\btheta), \btheta)^\top (\bar{\mQ}_i^k)^{-1} \nabla_{\vx} L(\vx_i^{k + 1}(\btheta)) \nonumber \\
    &= \frac{\partial r}{\partial \btheta}(\vx_i^{k + 1}(\btheta), \btheta)^\top d\vx_i^{k + 1}, \label{eq:loss_grad_final}
\end{align}
where $d\vx_i^{k + 1}$ is the unique solution to the linear system
\begin{equation*}
    \bar{\mQ}_i^k d\vx_i^{k + 1} = -\nabla_{\vx} L(\vx_i^{k + 1}(\btheta)).
\end{equation*}
Expanding the matrix multiplication in (\ref{eq:loss_grad_final}) yields
\begin{align*}
    \nabla_{\bar{\mQ}_i^k} L &= \frac{1}{2} (\vx_i^{k + 1} \otimes d\vx_i^{k + 1} + d\vx_i^{k + 1} \otimes \vx_i^{k + 1}), \\
    \nabla_{\bar{\vb}_i^k} L &= -d\vx_i^{k + 1}.
\end{align*}
\end{proof}


\section{Background on PAC-Bayes Theory}
\label{ssec:pac bayes background}

Here, we provide a brief overview of PAC-Bayes theory (\cite{alquier2024user}). Consider a bounded loss function $\ell(\zeta; \theta)$. Without loss of generality, we assume that this loss is uniformly bounded between $0$ and $1$. PAC-Bayes theory aims to providing a probabilistic bound for the true expected loss
\begin{equation}
\ell_{\calD}(\calP) = \E_{\zeta \sim \calD} ~ \E_{\theta \sim \calP} \left[ \ell (\zeta; \theta) \right],
\end{equation}
where $\calD$ is the data distribution --- in our case, this is the distribution optimization problems are drawn from. The empirical expected loss is given by,
\begin{equation}
\ell_\mathcal{\mathcal{S}}(\mathcal{P}) = \E_{\theta \sim \calP}\left[\frac{1}{H} \sum_{j = 1}^H (\zeta^j; \theta)\right],
\end{equation}
where $\calS = \{\zeta^j\}_{j = 1}^H$ is the training dataset consisting of $H$ problem instances.

The PAC-Bayes framework operates by forming a bound that holds in high probability on the true loss $\ell_{\calD}(\calP)$ in terms of the empirical loss and a the deviation between the learned policy $\calP$ and a prior policy $\calP_0$ used to as an initial guess for $\calP$. This deviation is measured using the KL divergence. Importantly, $\calP_0$ need not be a Bayesian prior but can be any distribution independent of the data used to train $\calP$ and evaluate the sample loss. Moreover, $\ell(\zeta; \theta)$ need not be the loss used to train $\calP$, but can be any bounded function. This observation is useful because, both in the literature and in the sequel, it is common to use a loss function modified for practicality during training before evaluating the bound using the loss function of interest.

Specifically, the following PAC-Bayes bounds hold with probability $ 1 - \delta$,
\begin{equation}
\ell_{\calD}(\calP)
\leq 
\sD_{\text{KL}}^{-1}    
\left(
\ell_{\calS} (\calP) \| 
\frac{\sD_{\text{KL}}(\calP \| \calP_0) + \log \frac{2 \sqrt{H}}{\delta}}{H}
\right)
\leq
\ell_\calS(\calP)
+ \sqrt{
\frac{\sD_{\text{KL}}(\calP \| \calP_0) + \log \frac{2 \sqrt{H}}{\delta}}{2 H}
}, \label{eq:pac bayes}
\end{equation}
where the $\sD_{\text{KL}}^{-1}(p \| c)$ is the \emph{inverse of the KL divergence} for Bernoulli random variables $\mathcal{B}(p), \mathcal{B}(q)$:
\begin{align}
\sD_{\text{KL}}^{-1}(p \| c) = \sup \{q \in [0, 1]\ |\ \sD_{\text{KL}}(\mathcal{B}(p) \| \mathcal{B}(q)) \leq c\}.
\end{align}

The probability $\delta$ captures the failure case that the data set $\calS$ is not sufficiently representative of the data distribution $\calD$. In the sequel, both of the above inequalities will be used. As the first bound is tighter, it is used to evaluate the generalization capabilities of the learned optimizer. The benefit of the second, loser, bound is that its form is convenient to use during training as a regularizer. Using both bounds in this manner is a common technique in the PAC-Bayes literature (\cite{majumdar2021pac}, \cite{dziugaite2017computing}).

\section{Optimizing and Evaluating Generalization Bound}
\label{sec: opt for gen bounds}

Two important requirements for establishing a tight PAC-Bayes bound are selecting an informative prior and optimizing the PAC-Bayes bounds in (\ref{eq:pac bayes}) instead of simply minimizing the loss function. The choice of prior $\calP_0$ is particularly important because the KL divergence is unbounded and can produce a vacuous result \cite{dziugaite2021role}. While the distribution $\calP_0$ need not be a Bayesian prior, it  must be selected independently from the data used to optimize $\calP$ and evaluate the bound. To select $\calP_0$, we follow a common approach in the literature and split our training set $\calS$ into two disjoint subsets $\calS_0, \calS_1$. The prior $\calP_0$ is first trained using the data set $\calS_0$ and the loss $\ell(\calD; \Theta)$ discussed in (\ref{sec: distr deep qp}). 

Subsequently, the posterior $\calP$ is trained by minimizing the looser (i.e., rightmost) PAC-Bayes bound in (\ref{eq:pac bayes}). This bound is used for training because it is straightforward to evaluate in comparison to computing the inverse of the KL divergence, and this objective is easily interpreted as minimizing an expected loss function with a regularizer. To evaluate the loss function in the PAC-Bayes bound, parameters are sampled from $\calP$ using the current network weights and an empirical average is used. Once training is complete, the PAC-Bayes bound is evaluated as described in Theorem \ref{thm:gen bound}, i.e., by using the tighter PAC-Bayes bound in (\ref{eq:pac bayes}) and the sample convergence bound in (\ref{eq:sample convergence}).

\section{Details on Experiments}
\label{sec: appendix - exp details}

This section provides further details on the problems considered in the experiments, the training of the learned optimizers, as well as the evaluation of both learned and traditional methods.

\subsection{Problem Types in Centralized Experiments}

\paragraph{Random QPs.} 
We consider randomly generated problems of the following form
\begin{equation}
\min_{\vx} ~ \frac{1}{2} \vx^\top \mQ \vx + \vq^\top \vx
\quad \mathrm{s.t.} 
\quad \mA \vx \leq \vb,
~~ \mC \vx = \vd. 
\end{equation}
For each generated problem, the cost Hessian is given by $\mQ = \mF^\top \mF + \gamma \mI$, where each element of $\mF \in \sR^{n \times n}$ is sampled through $\mF_{ij} \sim \calN(0,1)$ and $\gamma = 1.0$. The coefficients of $\vq$ are also sampled as $\vq_i \sim \calN(0,1)$. The elements of the inequality constraints matrix $\mA \in \mathbb{R}^{m \times n}$ are given by $\mA_{ij} \sim \calN(0,1)$, while $\vb = \mA \btheta$, where each element of $\btheta \in \sR^n$ is sampled through $\btheta_i \sim \calN(0,1)$. Similarly, the elements of the equality constraints matrix $\mC \in \mathbb{R}^{p \times n}$ are given by $\mC_{ij} \sim \calN(0,1)$, while $\vd = \mC \bxi$, where each element of $\bxi \in \sR^n$ is $\bxi_i \sim \calN(0,1)$. For random QPs without equality constraints, we set $n = 50$, $m = 40$ and $p = 0$. For random QPs with equality constraints, we set $n = 50$, $m = 25$ and $p = 20$.

\paragraph{Optimal control.} We consider linear optimal control problems of the following form
\begin{subequations}
\begin{align}
\min_{\vx, \vu} ~ & \sum_{t = 0}^{T-1} \vx_t^\top \mQ \vx_t + \vu_t^\top \mR \vu_t + \vx_T^\top \mQ_T \vx_T
\\
\text{s.t.} 
\quad & \vx_{t+1} = \mA_{\text{d}} \vx_t + \mB_{\text{d}} \vu_t, \quad t = 0, \dots, T-1,
\\
& \mA_u \vu_t \leq \vb_u, \quad \mA_x \vx_t \leq \vb_x, \quad t = 0, \dots, T,
\quad 
\\
& \vx_0 = \bar{\vx}_0.
\end{align}
\end{subequations}
where $\vx = \{ \vx_0, \dots, \vx_T \}$ is the state trajectory, $\vu = \{ \vu_0, \dots, \vu_{T-1} \}$ is the control trajectory, $\bar{\vx}_0$ is the given initial state condition, $\mQ$ and $\mR$ are the running state and control cost matrices, $\mQ_T$ is the terminal state cost matrix, $\mA_{\text{d}}$ and $\mB_{\text{d}}$ are the dynamics matrices, and finally $\mA_u, \vb_u$ and $\mA_x, \vb_x$ are the control and state constraints coefficients, respectively. 

Both the double integrator and the mass-spring problem setups are drawn from \cite{chen2022large}. For the double integrator system, we have $x_t \in \sR^2$ and $u_t \in \sR$, with time horizon $T = 20.$ The dynamics matrices are given by
\begin{equation}
\mA_{\text{d}} = \begin{bmatrix}
1 & 1 \\
0 & 1
\end{bmatrix}, 
\quad 
\mB_{\text{d}} = \begin{bmatrix}
0.5 \\ 0.1 
\end{bmatrix}
\end{equation}
The cost matrices are $\mQ = \mQ_T = \mI_2$ and $R = 1.0$. The state and control constraint coefficients are given by
\begin{equation}
\mA_x = 
\begin{bmatrix}
\mI_2 \\ - \mI_2
\end{bmatrix}, 
\quad
\vb_x = 
\begin{bmatrix}
5 & 1 & 5 & 1
\end{bmatrix}^\top, 
\quad
\mA_u = 
\begin{bmatrix}
  1 \\ -1
\end{bmatrix}, 
\quad
\vb_u = 
\begin{bmatrix}
0.1 & 0.1
\end{bmatrix}^\top.
\end{equation}
Finally, the initial state conditions are sampled from the uniform distribution $\calU[ [-1; -0.3], [1; 0.3]]$.

For the oscillating masses, we have $x_t \in \sR^{12}$ and $u_t \in \sR^3$, with time horizon $T = 10.$ The discrete-time dynamics matrices are obtained from the continuous-time ones through Euler discretization,
\begin{equation}
\mA_{\text{d}} = \mI + \mA_{\text{c}} \Delta t, \quad 
\mB_{\text{d}}  = \mA_{\text{c}} \Delta t.
\end{equation}
The continuous-time dynamics matrices are given by
\begin{equation}
\mA_{\text{c}} = \begin{bmatrix}
\vzero_{6 \times 6} & \mI_6 \\
a \mI_6 + c \mL_6 + c \mL_6^\top & b \mI_6 + d \mL_6 + d \mL_6^\top
\end{bmatrix}, 
\quad 
\mB_{\text{c}} = \begin{bmatrix}
\vzero_{6 \times 3} \\ \mF
\end{bmatrix}
\end{equation}
with $c = 1.0$, $d = 0.1$, $a = -2 c$, $b = -2.0.$ $\mL_6$ is the $6 \times 6$ lower shift matrix and 
\begin{equation}
\mF = 
\begin{bmatrix}
\ve_1 & - \ve_1 & \ve_2 & \ve_3 & -\ve_2 & \ve_3
\end{bmatrix}^\top
\end{equation}
where $\ve_1, \ve_2, \ve_3$ are the standard basis vectors in $\sR^3$.

The timestep is set as $\Delta t = 0.5$. The cost matrices are $\mQ = \mQ_T = \mI_{12}$ and $\mR = \mI_3$. The state and control constraints are defined through
\begin{equation}
\mA_x = 
\begin{bmatrix}
\mI_{12} \\ - \mI_{12}
\end{bmatrix}, 
\quad
\vb_x = 
4 \cdot \mathbf{1}_{24}, 
\quad
\mA_u = 
\begin{bmatrix}
\mI_3 \\ - \mI_3
\end{bmatrix}, 
\quad
\vb_u = 
0.5 \cdot \mathbf{1}_6.
\end{equation}

The initial conditions $\bar{\vx}_0$ are sampled from $\calU \big[ [-1, 1]^{12} \big]$.

\paragraph{Portfolio optimization.} We consider the same portfolio optimization problem setup as in \cite{stellato2020osqp}. For completeness, we briefly repeat it here,
\begin{equation}
\max_{\vx} ~ \bmu^\top \vx - \gamma (\vx^\top \bSigma \vx)
\quad \mathrm{s.t.} 
\quad x_1 + \dots + x_n = 1, 
\quad \vx \geq \vzero,
\label{eq: portfolio problem}
\end{equation}
where $\vx \in \sR^n$ is the assets allocation vector, $\bmu \in \sR^n$ is the expected returns vector, $\bSigma \in \sR_+^N$ is the risk covariance matrix and $\gamma > 0$ is the risk aversion parameter. The matrix $\bSigma$ is of the form $\bSigma = \mF \mF^\top + \mD$ with $\mF \in \sR^{d \times n}$ is the factors matrix and $\mD \in \sR^{n \times n}$ is a diagonal matrix involving individual asset risks. Using an auxiliary variable $\vt = \mF^\top \vx$, then problem \eqref{eq: portfolio problem} is rewritten as
\begin{equation}
\min_{\vx, \vt} ~ \vx^\top \mD \vx + \vt^\top \vt - \frac{1}{\gamma} \bmu^\top \vx
\quad \mathrm{s.t.} 
\quad \vt = \mF^\top \vx,
\quad \mathbf{1}^\top \vx = 1, 
\quad \vx \geq \vzero.
\end{equation}
For the problems we are generating, we use $n = 250$, $k = 25$ and $\gamma = 1.0$. Each element of the expected return vector $\bmu$ is sampled through $\mu_i \sim \calN(0,1)$. The matrix $\mF$ consists of $50 \%$ non-zero elements sampled through $F_{ij} \sim \calN(0,1)$. Finally, the diagonal elements of $\mD$ are sampled with $\calD_{ii} \sim \calU[0, \sqrt{k}]$.

\paragraph{LASSO.} The least absolute shrinkage and selection operator (LASSO) is a linear regression technique with an added $\ell_1$-norm regularization term to promote sparsity in the parameters \citep{tibshirani1996regression}. We again consider the same problem setup as in \cite{stellato2020osqp}, where the initial optimization problem 
\begin{equation}
\min_{\vx} ~ \| \mA \vx - \vb \|_2^2 + \lambda \| \vx \|_1
\end{equation}
is rewritten as
\begin{equation}
\min_{\vx, \vt} ~ (\mA \vx - \vb)^\top (\mA \vx - \vb) + \lambda \mathbf{1}^\top \vt
\quad \mathrm{s.t.} 
\quad -\vt \leq \vx \leq \vt, 
\end{equation}
where $\vx \in \sR^n$ is the vector of parameters, $\mA \in \sR^{m \times n}$ is the data matrix, $\lambda$ is the weighting parameter, and $\vt \in \sR^n$ are newly introduced variables. The matrix $\mA$ consists of 15\% non-zero elements sampled through $\mA_{ij} \sim \calN(0,1).$ The true sparse vector $\vv \in \sR^n$ to be learned consists of 50\% non-zero elements sampled through $\vv_i \sim \calN(0, 1/n)$. We then construct $\vb = \mA \vv + \bxi$ where $\bxi_i \sim \calN(0,1)$ represents noise in the data. Finally, we set $\lambda = (1/5) \| \mA^\top \vb \|_\infty$. For the problems we are generating, we set $n = 100$ and $m = 10^4$.

\subsection{Problem Types in Distributed Experiments}

\paragraph{Random Networked QPs.} In this family of problems, we generate random QPs with an underlying network structure. Consider an undirected graph $\calG(\calV, \calE)$, where $\calV$ and $\calE$ are the nodes and edges sets, respectively. Each node $i$ is associated with a decision variable $\vx_i \in \sR^{n_i}$. Then, we generate problems of the following form 
\begin{subequations}
\begin{align}
& ~~~~~~~~~~~~~~~~~~~ \min_{\{\vx_i\}_{i \in \calV}} ~  
\sum_{i \in \calV}
\frac{1}{2} \vx_i^\top \mQ_i \vx_i + \vq_i^\top \vx_i
\\[0.1cm]
& \mathrm{s.t.} \quad 
\mA_{ij} 
\begin{bmatrix}
\vx_i
\\
\vx_j
\end{bmatrix}
\leq \vb_{ij}, 
\quad 
\mC_{ij} 
\begin{bmatrix}
\vx_i
\\
\vx_j
\end{bmatrix}
= \vd_{ij},
\quad (i,j) \in \calE.
\end{align}
\label{eq: network random qp}
\end{subequations}
For each generated problem, a cost Hessian is constructed as $\mQ_i = \mF_i^\top \mF_i + \gamma \mI$, where each element of $\mF_i \in \sR^{n_i \times n_i}$ is sampled through $\mF_{i}^{kl} \sim \calN(0,1)$ and $\gamma = 1.0$. The elements of the cost coefficients vectors $\vq_i$ are also sampled through $\vq_i^k \sim \calN(0,1)$. The elements of the inequality constraints matrix $\mA_{ij} \in \sR^{m_{ij} \times (n_i + n_j)}$ are given by $\mA_{ij}^{kl} \sim \calN(0,1)$. The vectors $\vb_{ij} \in \sR^{m_{ij}}$ are obtained through $\vb_{ij} = \mA_{ij} \btheta_{ij}$, where each element of $\btheta_{ij} \in \sR^{n_i + n_j}$ is sampled through $\btheta_{ij}^k \sim \calN(0,1)$. In a similar manner, the elements of the equality constraints matrices $\mC_{ij} \in \sR^{p_{ij} \times (n_i + n_j)}$ are generated through $\mC_{ij}^{kl} \sim \calN(0,1)$, while the vectors $\vd_{ij} \in \sR^{p_{ij}}$ are acquired through $\vd_{ij} = \mC_{ij} \bxi_{ij}$, where each element of $\bxi_{ij} \in \sR^{n_i + n_j}$ is generated with $\bxi_{ij}^k \sim \calN(0,1)$.

It is straightforward to observe that problems of the form (\ref{eq: network random qp}) can be cast in the form (\ref{eq: distr qp problem}) by introducing the augmented node variables $\vx_i^{aug} = [x_i, \{x_j\}_{j \in \calN_i}]^\top$. The problem data can then be augmented based on this new $\vx_i^{aug}$ to yield the desired problem structure. Most notably, the constraints can be rewritten as $\mA_i^{aug} \vx_i^{aug} \leq b_i^{aug}$ and $\mC_i^{aug} x_i^{aug} = d_i^{aug}$, respectively. In our experiments, the underlying graph structure is a square grid. For random QPs without equality constraints, we set $n_i = 10$, $m_{ij} = 5$, and $p_{ij} = 0$. For random QPs with equality constraints, we set $n_i = 10$, $m_{ij} = 3$, and $p_{ij} = 2$ for the $N = 16$ training experiment and $p_{ij} = 1$ for the rest of the testing experiments until $N = 1,024$.

\paragraph{Multi-agent optimal control.} 
We adapt the distributed MPC problem from \citep{conte2012computational,conte2012distributed}, which generalizes to different systems based on the choice of dynamics matrices, as described below. The optimization problem is given as
\begin{subequations} \label{eq:distributed_mpc}
\begin{align}
    \min_{\vx, \vu} ~ & \sum_{i \in V} \sum_{t = 0}^{T - 1} (\vx_i^t)^\top \mQ_i \vx_i^t + (\vu_i^t)^\top \mR_i \vu_i^t + (\vx_i^T)^\top \mP_i \vx_i^T, \label{eq:distr_mpc_obj} \\
    \mathrm{s.t.} \quad & \vx_i^{t + 1} = \mA_{ii} \vx_i^t + \mB_i \vu_i^t + \sum_{j \in \calN_i} \mA_{ij} \vx_j^t, \quad t = 0, \ldots, T - 1, \quad i \in \calV \label{eq:distr_mpc_dynamics} \\
    & \mG_x^i \vx_i^t \leq \vf_x^i, \, \mG_u^i \vu_i^t \leq \vf_u^i, \quad t = 0, \ldots, T, \quad i \in \calV \label{eq:distr_mpc_local_constraints} \\
    & \vx_i^0 = \bar{\vx}_i^0, \quad i \in \calV, \label{eq:distr_mpc_init_cond}
\end{align}
\end{subequations}
where $\vx_i^t$ and $\vu_i^t$ are the state and control for agent $i$ at time $t$. (\ref{eq:distr_mpc_dynamics}) describes the dynamics and the coupling between the agents, (\ref{eq:distr_mpc_local_constraints}) describe local inequality constraints, and (\ref{eq:distr_mpc_init_cond}) describes the initial condition for each of the agents.

For the coupled pendulums, the individual state $\vx_i^t \in \sR^2$ for each agent consists of the angle and angular velocity of the pendulum and the control $\vu_i^t \in \sR^1$ is the torque. The dynamics matrices are given as
\begin{align*}
    \mA_{ii} = \begin{bmatrix}
        1 & dt \\
        -(\frac{g}{\ell} + \frac{\text{nn}(i) k}{m}) dt & 1 - \frac{\text{nn}(i) c}{m}dt
    \end{bmatrix}, \quad \mA_{ij} = \begin{bmatrix}
        0 & 0 \\
        \frac{k}{m} dt & \frac{c}{m} dt
    \end{bmatrix}, \quad \mB_i = \begin{bmatrix}
        0 \\ \frac{1}{m \ell^2} dt
    \end{bmatrix},
\end{align*}
where $dt = 0.1$ is the discretization step size, $g = 9.81$ is the gravitational constant, $m = 1.0$ is the mass of each pendulum, $\ell = 0.5$ is the length of each pendulum, $\text{nn}(i)$ is the number of neighbors of agent $i$, $k = 0.1$ is the spring constant between each pendulum, and $c = 0.1$ is the damping constant between each pendulum. We have used the small angle assumption $\sin \theta \approx \theta$ so the dynamics are linear and therefore the optimization is convex. There are no inequality constraints for the coupled pendulums. The initial states are sampled uniformly from $\calU[-\pi, \pi]$. Finally, we considered $N = 10$ and $T = 30$.

For the coupled oscillating masses, we adapt the same benchmark system from \citet{chen2022large} used in the non-distributed experiments. The individual state $\vx_i^t \in \sR^2$ for each agent consists of the displacement and velocity of the mass and the control $\vu_i^t \in \sR^1$ is the force acting on the mass. The dynamics matrices are
\begin{align*}
    \mA_{ii} = \begin{bmatrix}
        1 & dt \\
        -\frac{2k}{m} dt & 1 - \frac{2 c}{m}dt
    \end{bmatrix}, \quad \mA_{ij} = \begin{bmatrix}
        0 & 0 \\
        \frac{k}{m} dt & \frac{c}{m} dt
    \end{bmatrix}, \quad \mB_i = \begin{bmatrix}
        0 \\ \frac{1}{m} dt
    \end{bmatrix},
\end{align*}
where $dt = 0.5$ is the discretization step size, $m = 1.0$ is the mass, $k = 0.4$ is the spring constant between each mass, and $c = 0.1$ is the damping constant between each mass. The initial states are sampled uniformly from $\calU[-2.0, 2.0]$. Inequality constraints $-4 \leq \vx_i^t \leq 4$ and $-0.5 \leq \vu_i^t \leq 0.5$ are represented as
\begin{align*}
    \mG_x^i = \begin{bmatrix}
        \mI_2 \\
        -\mI_2
    \end{bmatrix}, \quad \vf_x^i = 4 \cdot \mathbf{1}_4, \quad \mG_u^i = \begin{bmatrix}
        1 \\
        -1
    \end{bmatrix}, \quad \vf_u^i = 0.5 \cdot \mathbf{1}_2,
\end{align*}

For both the distributed MPC problems described above, the cost matrices are taken to be identity matrices: $\mQ_i = \mI_2$, $\mR_i = \mI_1$, and $\mP_i = \mI_2$, for all $i \in \calV$.

The optimization (\ref{eq:distributed_mpc}) can be expressed in the form of (\ref{eq: distr qp problem}) by defining an augmented vector consisting of the individual agent's states and controls, as well as the states and controls of its neighbors. Letting $\vz_i = [ \vx_i^0, \vu_i^0, \ldots, \vx_i^T ]^\top$, the augmented optimization vector for each agent $i$ is given as $\vx_i^\text{aug} = [\vz_i, \{ \vz_j \}_{j \in \calN_i}]^\top$. The cost, dynamics, and constraint matrices can be augmented straightforwardly based on this new $\vx_i^\text{aug}$. For all problems, we considered $T = 15$.

\paragraph{Network flow.}

The network flow problem is adapted from \citet{mota2013communication,mota2014distributed}. We consider a directed regular graph with 200 nodes and 1000 directed edges $x_{ij} \in \calE$. Each edge has an associated quadratic cost function $\phi_{ij}(x_{ij}) = \frac{1}{2} (x_{ij} - a_{ij})^2$, where $a_{ij}$ is sampled from $[1.0, 2.0, 3.0, 4.0, 5.0, 10.0]$ with probabilities $[0.2, 0.2, 0.2, 0.2, 0.1, 0.1]$. The objective is to optimize the flow through the graph subject to equality constraints on the flow into and out of each node. Namely, the flow into each node should be equal to the flow out of the node. For node $i$, the flow conservation constraint is $\sum_{j \in \calE_i^-} x_{ji} = \sum_{k \in \calE_i^+} x_{ik}$, where $\calE_i^-$ is the set of all incoming edges to node $i$, and similarly $\calE_i^+$ is the set of all outgoing edges from node $i$.
100 nodes are randomly selected and injected with an external flow $f_k$ sampled identically to $a_{ij}$. For each of these nodes, a reachable descendant is randomly selected and an equivalent amount of flow $f_k$ is removed from those nodes.

This problem is straightforward to express in the form (\ref{eq: distr qp problem}) by considering each node as an individual agent and defining the local state vector for each agent as
\begin{align}
    \vx_i = \begin{bmatrix}
        \{ x_{ji} \}_{j \in \calE_i^-} \\
        \{ x_{ik} \}_{k \in \calE_i^+} 
    \end{bmatrix},
\end{align}
consisting of all the incoming and outgoing edges for node $i$. Each agent is responsible for its own flow constraint defined by
\begin{align}
    \mA_i = \begin{bmatrix}
        \{ 1 \}_{j \in \calE_i^-} & \{ -1 \}_{k \in \calE_i^+} \\
        \{ -1 \}_{j \in \calE_i^-} & \{ 1 \}_{k \in \calE_i^+}
    \end{bmatrix}, \quad \vb_i = \vzero,
\end{align}
where $\vb_i$ might instead contain the external injected or removed flow $f_i$ for that node $i$. The augmented cost matrix $\mQ_i$ is zero for all incoming edges and has entries $1/2$ on the diagonal of the outgoing edges. The augmented cost vector $\vq_i$ contains each of the quadratic cost offsets $a_{ik}$:
\begin{align}
    \mQ_i = \begin{bmatrix}
        \{0\}_{j \in \calE_i^-} & \\
        & \{ \frac{1}{2} \}_{k \in \calE_i^+}
    \end{bmatrix}, \quad \vq_i = \begin{bmatrix}
        \{0\}_{j \in \calE_i^-} \\
        \{ -a_{ik} \}_{k \in \calE_i^+}
    \end{bmatrix}.
\end{align}
Finally, we impose the constraint $-f_{\max} \cdot \mathbf{1} \leq \vx_i \leq f_{\max} \cdot \mathbf{1}$ on the maximum allowed flow of all edges, with $f_{\max} = 5$.

\paragraph{Distributed LASSO.} Distributed LASSO \citep{mateos2010distributed} extends LASSO to situations where the training data are distributed across different agents and agents cannot share training data with each other. It can be formulated as
\begin{equation}
\min_{\{\vx_i\}_{i = 1}^N, \vw} ~  
\sum_{i = 1}^N
\| \mA_i \vx_i - \vb_i \|_2^2 + \frac{\lambda}{N} \| \vx_i \|_1
\quad \mathrm{s.t.} 
\quad \vx_i = \vw,
\quad i = 1, ..., N
\label{eq: first distributed LASSO formulation}
\end{equation}
where $\vw \in \sR^{n_i}$ is a global vector of regression coefficients, $\vx_i \in \sR^{n_i}$ is a local copy of $\vw$, $\mA_i \in \sR^{m_i \times n_i}$ and $\vb \in \sR^{m_i}$ are the training data available to agent $i$, and $\lambda$ is the weighting parameter. Similarly to non-distributed LASSO, this formulation is rewritten as 
\begin{subequations}
\begin{align}
\min ~ & 
\sum_{i = 1}^N
(\mA_i \vx_i - \vb_i)^\top (\mA_i \vx_i - \vb_i) + \frac{\lambda}{N} \mathbf{1}^\top \vt_i \\
\mathrm{s.t.} 
\quad & \vt_i \leq \vx_i \leq \vt_i, \quad \vx_i = \vw,
\quad \vt_i = \vg,
\quad i = 1, ..., N
\label{eq: second distributed LASSO formulation}
\end{align}
\end{subequations}
where $\vt_i \in \sR^{n_i}$ are newly-introduced variables and $\vg$ is the global copy of $\vt_i$.

The matrix $\mA_i$ consists of 15\% non-zero elements sampled through $\mA_i^{kl} \sim \calN(0, 1).$ The true sparse vector $\vv \in \sR^n$ to be learned consists of 50\% non-zero elements sampled through $\vv_i \sim \calN(0, 1/n)$. We then construct $\vb = \mA \vv + \bxi$ where $\bxi_i \sim \calN(0,1)$ represents noise in the data. Finally, we set $\lambda = (1/5) \max_{i} ( \| \mA_i^\top \vb_i \|_\infty)$. For the problems, we have $n_i = 50$ and $m_i = 5 \cdot 10^3.$

\begin{table}[t]
\begin{center}
\footnotesize{
\scalebox{1.0}{
\begin{tabular}{cccccc}
\toprule
Problem Class & Layers $K$ & Train samples & Epochs & Train time & Test samples
\\ \hline
Random QPs 
& 30 & 2,000 & 125 & 21min & 1,000
\\ 
Random QPs with Eq. Constraints
& 30 & 2,000 & 125 & 23min & 1,000
\\ 
Double Integrator
& 30 & 500 & 300 & 28min & 1,000
\\ 
Osc. Masses
& 15 & 500 & 300 & 48min & 1,000
\\ 
Portfolio Optimization
& 30 & 500 & 300 & 1h 14min & 1,000
\\ 
LASSO
& 10 & 500 & 300 & 20min & 1,000
\\ \bottomrule
\end{tabular}
}
}
\end{center}
\caption{Training and testing details for DeepQP.}
\label{tab: DeepQP train details}
\end{table}

\begin{table}[t]
\begin{center}
\footnotesize{
\scalebox{1.0}{
\begin{tabular}{cccccc}
\toprule
Problem Class & Layers $K$ & Training samples & Epochs & Train time & Test samples
\\ \hline
Random QPs 
& 50 & 1,000 & 300 & 3h 21min & 500
\\ 
Random QPs with Eq. Constraints
& 50 & 500 & 600 & 3h 29min & 500
\\ 
Coupled Pendulums
& 20 & 500 & 400 & 1h 49min & 500
\\ 
Coupled Osc. Masses
& 20 & 500 & 600 & 2h 29min & 500
\\ 
Network Flow
& 30 & 500 & 600 & 2h 8min & 500
\\ 
Distributed LASSO
& 20 & 500 & 600 & 56min & 500
\\ \bottomrule
\end{tabular}
}
}
\end{center}
\caption{Training and testing details for DeepDistributedQP.}
\label{tab: DeepDistrQP train details}
\end{table}

\subsection{Details on Training and Testing}
\label{sec: appendix - train details}

Here, we discuss details regarding the training and testing of DeepQP and DeepDistributedQP in the presented experiments.

\paragraph{Centralized experiments.} 
Table \ref{tab: DeepQP train details} shows the number of layers $K$, training dataset size, number of epochs, total training time and testing dataset size for DeepQP in every centralized problem. The increased dataset size and number of epochs for RandomQPs is motivated by the fact that the structure in these problems is less clear; learning policies that exploit this structure therefore requires more examples and takes longer. In all experiments, DeepQP was trained with a batch size of $50$ using the Adam optimizer with learning rate $10^{-3}$. The feedback layers are set as $2 \times 16$ MLPs. DeepQP and OSQP always start with zero initializations in all comparisons. The weights of the training loss were set to $\gamma_k = \exp\left(\left(k - K\right)/5\right)$ in all experiments. Both the training and testing datasets are contructed after letting OSQP running until optimality.

\begin{table}[t]
\begin{center}
\footnotesize{
\scalebox{1.0}{
\begin{tabular}{cc}
\toprule
Problem Class & List of penalty parameters $\rho$
\\ \hline
Random QPs 
& 0.1, 0.3, \dots, 3, 10
\\ 
Random QPs with Eq. Constraints
& 0.1, 0.3, \dots, 3, 10
\\ 
Double Integrator
& 3, 5, \dots, 100, 300
\\ 
Osc. Masses
& 0.1, 0.3, \dots, 3, 10
\\ 
Portfolio Optimization
& 3, 5, \dots, 100, 300
\\ 
LASSO
& 30, 50, \dots, 1000, 3000
\\ \bottomrule
\end{tabular}
}
}
\end{center}
\caption{List of OSQP penalty parameters used in centralized experiments.}
\label{tab: OSQP pen params}
\end{table}

\begin{table}[t]
\begin{center}
\footnotesize{
\scalebox{1.0}{
\begin{tabular}{cc}
\toprule
Problem Class & List of penalty parameters $\rho$ and $\mu$
\\ \hline
Random QPs 
& 0.1, 0.3, \dots, 3, 10
\\ 
Random QPs with Eq. Constraints
& 0.1, 0.3, \dots, 3, 10
\\ 
Coupled Pendulums
& 0.1, 0.3, \dots, 3, 10
\\ 
Coupled Osc. Masses
& 0.1, 0.3, \dots, 3, 10
\\ 
Network Flow
& 0.1, 0.3, \dots, 3, 10
\\ 
Distributed LASSO
& 30, 50, \dots, 1000, 3000
\\ \bottomrule
\end{tabular}
}
}
\end{center}
\caption{List of DistributedQP penalty parameters used in distributed experiments}
\label{tab: DistrQP pen params}
\end{table}

\paragraph{Distributed experiments.} Table \ref{tab: DeepDistrQP train details} shows the number of layers $K$, training dataset size, number of epochs, total training time and testing dataset size for DeepDistributedQP in every distributed problem. In all experiments, DeepDistributedQP was trained with a batch size of $50$ using the Adam optimizer with learning rate $10^{-3}$. The feedback layers are set as $2 \times 16$ MLPs. DeepDistributedQP and DistributedQP always start with zero initializations in all comparisons. In all experiments, the weights of the training loss were set to $\gamma_k = \exp\left(\left(k - K\right)/5\right)$. For the low-dimensional testing datasets, these datasets are constructed using OSQP. For larger scales, the testing dataset is constructed with DistributedQP instead as it is much faster (see Table \ref{tab: wall-clock-appendix}), after ensuring convergence to optimality.


\paragraph{Generalization bounds experiments.} These experiments were performed on a networked random QPs problem with $N = 16, n_i = 10, m_{ij} = 5, p_{ij} = 0$ and on a coupled pendulums problem with $N = 10$ and the same parameters as described in the previous section. The prior was obtained through training on a small separate dataset of $500$ problems for $50$ epochs. The posterior was then acquired through optimizing for the generalization bound with a dataset of $15,000$ problems for $100$ epochs.

\begin{table}[t]
\begin{center}
\scalebox{0.69}{
\begin{tabular}{cccc>{\columncolor{green!20}}c>{\columncolor{green!20}}c>{\columncolor{yellow!20}}c>{\columncolor{yellow!20}}c>{\columncolor{magenta!20}}c>{\columncolor{magenta!20}}c>{\columncolor{red!20}}c>{\columncolor{red!20}}c}
\hline
  & & & 
& \multicolumn{2}{>{\columncolor{green!20}}c}{\textbf{DeepDistrQP (ours)}}  
& \multicolumn{2}{>{\columncolor{yellow!20}}c}{\textbf{DistrQP (ours)}}
& \multicolumn{2}{>{\columncolor{magenta!20}}c}{\textbf{OSQP (Indirect)}}
& \multicolumn{2}{>{\columncolor{red!20}}c}{\textbf{OSQP (Direct)}}
\\
\toprule
\multicolumn{12}{c}{\textbf{Networked Random QPs}}  
\\
\hline
$N$ & $n$ & $m$ & $\texttt{nnz}(\mQ, \mA)$ 
& Time & Iters
& Time & Iters
& Time & Iters
& Time (1st iter.) & Iters
\\ \hline
16 & 160 & 120 & 4,000 
& 33.05 ms & 50 & 141.9 ms & 208 
& 46.16 ms & 29 & \textbf{0.86 ms} & 29
\\ 
64 & 640 & 560 & 17,600 
& 39.11 ms & 50 & 129.2 ms & 192
& 185.1 ms & 28 & \textbf{23.8 ms} & 28
\\ 
256 & 2,560 & 2,400 & 73,600
& \textbf{50.21 ms} & 50 & 128.8 ms & 168 
& 514 ms & 23 & 703.5 ms & 23
\\ 
1,024 & 10,240 & 9,920 & 300,800
& \textbf{62.68 ms} & 50 & 158.9 ms & 165 
& 3.03s & 23 & 8.20 s & 23 
\\
\toprule
\multicolumn{12}{c}{\textbf{Networked Random QPs with Equality Constraints}}  
\\
\hline
$N$ & $n$ & $m$ & $\texttt{nnz}(\mQ, \mA)$ 
& Time & Iters
& Time & Iters
& Time & Iters
& Time (1st iter.) & Iters
\\ \hline
16 & 160 & 168 & 4,960
& 37.21 ms & 50 & 138.9 ms  & 170
& 36.52 ms & 19 & \textbf{0.76 ms} & 19
\\ 
64 & 640 & 560 & 17,600
& 57.76 ms & 50 & 238.1 ms & 172
& 109.0 ms & 17 & \textbf{26.9 ms} & 17
\\ 
256 & 2,560 & 2,400 & 73,600
& \textbf{74.54 ms} & 50 & 239.5 ms & 164
& 692.5 ms & 17 & 956.0 ms  & 17
\\ 
1,024 & 10,240 & 9,920 & 300,800
& \textbf{82.55 ms} & 50 & 371.0 ms & 172
& 5.83 s & 16 & 11.60 s & 16
\\
\toprule
\multicolumn{12}{c}{\textbf{Coupled Pendulums Optimal Control}}  
\\
\hline
$N$ & $n$ & $m$ & $\texttt{nnz}(\mQ, \mA)$ 
& Time & Iters
& Time & Iters
& Time & Iters
& Time (1st iter.) & Iters
\\ \hline
10 & 470 & 640 & 3,690
& 50.99 ms & 20 & 89.81 ms & 35
& 49.46 ms & 8 & \textbf{4.95 ms} & 8
\\ 
20 & 940 & 1,200 & 7,500
& \textbf{66.44 ms} & 20 & 116.7 ms & 35
& 372.0 ms & 8 & 199.7 ms & 8
\\ 
50 & 2,350 & 3,200 & 18,930
& \textbf{75.9 ms} & 20 & 142.1 ms & 34
& 948.8 ms & 8 & 4.38 s & 8
\\ 
100 & 4,700 & 6,400 & 37,980
& \textbf{101.9 ms} & 20 & 201.9 ms & 35
& 3.97 s & 9 & 19.91 s & 9
\\ 
200 & 9,400 & 12,800 & 76,080
& \textbf{146.0 ms} & 20 & 284.8 ms & 34
& 22.41 s & 8 & 90.07 s & 8
\\ 
500 & 23,500 & 32,000 & 190,380
& \textbf{204.3 ms} & 20 & 379.8 ms & 36
& 112.9 s & 9 & \multicolumn{2}{>{\columncolor{red!20}}c}{Out of memory} 
\\ 
1,000 & 47,000 & 64,000 & 380,880
& \textbf{317.2 ms} & 20 & 628.2 ms & 34
& \multicolumn{2}{>{\columncolor{magenta!20}}c}{Out of memory} & \multicolumn{2}{>{\columncolor{red!20}}c}{Out of memory}  
\\
\toprule
\multicolumn{12}{c}{\textbf{Coupled Oscillating Masses Optimal Control}}  
\\
\hline
$N$ & $n$ & $m$ & $\texttt{nnz}(\mQ, \mA)$ 
& Time & Iters
& Time & Iters
& Time & Iters
& Time (1st iter.) & Iters
\\ \hline
10 & 470 & 1,580 & 4,590 
& \textbf{48.22 ms} & 20 & 73.58 ms & 33
& 79.1 ms & 9 & 178.4 ms & 9
\\ 
20 & 940 & 3,160 & 9,300
& \textbf{67.93 ms} & 20 & 91.53 ms & 33
& 641.9 ms & 9 & 2.37 s & 9
\\ 
50 & 2,350 & 7,900 & 23,430 
& \textbf{73.92 ms} & 20 & 97.34 ms & 32
& 1.07 s & 8 & 28.1 s & 8
\\ 
100 & 4,700 & 15,800 & 46,980 
& \textbf{91.93 ms} & 20 & 148.8 ms & 33
& 5.45 s & 8 & 132 s & 8
\\ 
200 & 9,400 & 31,600 & 94,080 
& \textbf{109.4 ms} & 20 & 194.4 ms & 34
& 31.8 s & 8 & 614 s & 8
\\ 
300 & 28,200 & 47,400 & 141,180 
& \textbf{132.8 ms} & 20 & 304.8 ms & 33
& 243 s & 8 & \multicolumn{2}{>{\columncolor{red!20}}c}{Out of memory} 
\\
\toprule
\multicolumn{12}{c}{\textbf{Network Flow}}  
\\
\hline
$N$ & $n$ & $m$ & $\texttt{nnz}(\mQ, \mA)$ 
& Time & Iters
& Time & Iters
& Time & Iters
& Time (1st iter.) & Iters
\\ \hline
20 & 100 & 140 & 600
& 6.80 ms & 30 & 10.68 ms & 50
& 9.51 ms & 15 & \textbf{0.59 ms} & 15
\\ 
50 & 250 & 350 & 1,500
& 7.81 ms & 30 & 13.17 ms & 48
& 14.81 ms & 16 & \textbf{1.30 ms} & 16
\\ 
200 & 1,000 & 1,400 & 6,000
& \textbf{12.08 ms} & 30 & 17.61 ms & 42
& 208.19 ms & 17 & 61.93 ms & 17 
\\ 
500 & 2,500 & 3,500 & 15,000
& \textbf{13.63 ms} & 30 & 19.73 ms & 40
& 425.7 ms & 17 & 745.2 ms & 17
\\ 
1,000 & 5,000 & 7,000 & 30,000
& \textbf{20.51 ms} & 30 &  31.62 ms & 40
& 8.73 s & 18 & 11.59 s & 18
\\ 
2,000 & 10,000 & 14,000 & 60,000
& \textbf{29.86 ms} & 30 &  47.22 ms & 40
& 51.6 s & 18 & 73.9 s & 18 
\\ 
5,000 & 25,000 & 35,000 & 150,000
& \textbf{61.23 ms} & 30 & 85.99 ms & 39
& 558 s & 18 & \multicolumn{2}{>{\columncolor{red!20}}c}{Out of memory} 
\\
\toprule
\multicolumn{12}{c}{\textbf{Distributed LASSO}}  
\\
\hline
$N$ & $n$ & $m$ & $\texttt{nnz}(\mQ, \mA)$ 
& Time & Iters
& Time & Iters
& Time & Iters
& Time (1st iter.) & Iters
\\ \hline
10 & 1,100 & 3,000 & 29,000
& \textbf{15.06 ms} & 20 & 28.57 ms & 37
& 2.04 s & 33 & 148.2 ms & 33
\\ 
50 & 5,500 & 15,000 & 145,000
& \textbf{24.92 ms} & 20 & 44.27 ms & 38
& 13.74 s & 31 & 49.21 s & 31
\\ 
100 & 10,100 & 30,000 & 290,000
& \textbf{30.51 ms} & 20 & 51.44 ms & 35
& 85.92 s & 32 & 342.9 s & 32
\\ 
200 & 20,100 & 60,000 & 580,000
& \textbf{40.88 ms} & 20 & 76.21 ms & 36
& 418.9 s & 32 & \multicolumn{2}{>{\columncolor{red!20}}c}{Out of memory} 
\\ 
500 & 50,100 & 150,000 & 1,450,000
& \textbf{69.19 ms} & 20 & 130.24 ms & 35
& \multicolumn{2}{>{\columncolor{magenta!20}}c}{Out of memory} & \multicolumn{2}{>{\columncolor{red!20}}c}{Out of memory} 
\\ \hline
\end{tabular}
}
\end{center}
\caption{\textbf{Wall-clock times and iterations for DeepDistributedQP, DistributedQP, OSQP (indirect) and OSQP (direct).} This comparison shows the total wall-clock times for DistributedQP and OSQP (indirect or direct method) required to reach the same accuracy as DeepDistributedQP. For OSQP with direct method, we only report the time for the first iteration, assuming the best-case scenario in which the factorized KKT matrix can be reused for all subsequent iterations. Both DeepDistributedQP and DistributedQP demonstrate orders-of-magnitude improvements compared to OSQP as scale increases. In additon, DeepDistributedQP maintains a significant advantage over its standard optimization counterpart in all cases.}
\label{tab: wall-clock-appendix}
\end{table}

\subsection{Details on Standard Optimizers}

\paragraph{Details on OSQP.}
When comparing with OSQP using fixed penalty parameters, we selected the best-performing subsequence of $\{..., 0.1, 0.3, 0.5, 1.0, 3.0, 5.0, ...\}$ as the penalty parameters to plot against. Table \ref{tab: OSQP pen params} shows these parameters for every centralized problem in our experiments. For equality constraints, we scaled $\rho$ by $10^3$, as in \cite{stellato2020osqp}. For the adaptive version, we prefered the standard heuristic adaptation rule shown in \cite{boyd2011distributed} with $\tau = 2.0$ and $\mu = 10.0$, instead of the OSQP adaptation scheme \citep{stellato2020osqp}, as it performed better in our problem instances. We hypothesize that this might be due to the fact that as scale increases the infinity norm is ignoring more information that the 2-norm. The initial $\rho^0$ was initialized as the median of the range of fixed penalty parameters. 


\paragraph{Details on DistributedQP.} The range of fixed penalty parameters to compare with was chosen using the same methodology as with OSQP. Table \ref{tab: DistrQP pen params} shows these parameters for every distributed problem in our experiments. For the adaptive version, we used the standard heuristic adaptation rule shown in \cite{boyd2011distributed} with $\tau = 2.0$ and $\mu = 10.0$. The initial value was again always chosen as the median value of the above lists.

\subsection{Details on Wall-Clock Times}
\label{sec: appendix wall clock}

In Table \ref{tab: wall-clock-appendix}, we list the observed wall-clock times for DeepDistributedQP (ours), DistributedQP (ours) and OSQP using either the indirect or the direct method. The table presents all six studied problems with an increasing dimension. As clearly observed, DeepDistributedQP and DistributedQP demonstrate a substantially more favorable scalability than OSQP. In fact, the two algorithms can efficiently solve problems that OSQP cannot handle due to memory overflow on our system. Finally, DeepDistributedQP also maintains a clear advantage over its standard optimization counterpart DistributedQP across all experiments which signifies the importance of learning policies for the algorithm parameters.